\def\shortversion{January 3, 2001}
\def\longversion{Submitted September 5, 1999; Revised \shortversion}
\documentclass[reqno,12pt]{amsart}
\usepackage{amssymb}
\newcommand{\old}[1]{}
\renewcommand{\Pr}{{\mathbb P}}
\newcommand{\SAT}{\operatorname{SAT}}
\newcommand{\UNSAT}{\operatorname{UNSAT}}
\newcommand{\TRUE}{\operatorname{TRUE}}
\newcommand{\FALSE}{\operatorname{FALSE}}
\newcommand{\Var}{\operatorname{Var}}
\newcommand{\Bin}{\operatorname{Binomial}}
\newcommand{\Proof}[1]{\noindent{\it #1} \hspace*{.1em}}
\newenvironment{proof_of}[1]{
    \Proof{Proof of #1.}}{
    \qed}
\newenvironment{proof_left}[1]{
    \Proof{#1 in the subcritical regime.}}{
    \lqed}
\newcommand{\lqed}{\hspace*{\fill} $\underset{\displaystyle\square}\leftarrow$}
\newcommand{\rqed}{\hspace*{\fill} $\underset{\displaystyle\square}\rightarrow$}
\renewcommand{\qed}{\hspace*{\fill} $\square$}
\renewcommand{\o}{o_{\varepsilon,\lambda_n}}
\newenvironment{proof_right}[1]{
    \Proof{#1 in the supercritical regime.}}{
    \rqed}

\newtheorem{theorem}{Theorem}[section]
\newtheorem{lemma}[theorem]{Lemma}
\newtheorem{corollary}[theorem]{Corollary}

\theoremstyle{definition}
\newtheorem{definition}[theorem]{Definition}
\theoremstyle{remark}
\newtheorem{remark}[theorem]{\it Remark}

\numberwithin{equation}{section}

\begin{document}

\begin{quote}
\phantom{.}
\vspace{-15mm}
\rightline{\small Technical Report}
\rightline{\small MSR-TR-99-41}
\end{quote}
\vspace{10mm}

\title[Scaling Window of the 2-SAT Transition]%
{The Scaling Window of the 2-SAT Transition}

\author[B. Bollob\'as, C. Borgs, J. T. Chayes, J. H. Kim, D. B. Wilson]{}

\address{B\'ela Bollob\'as\newline
Department of Mathematical Sciences\newline
University of Memphis\newline
Memphis, TN 38152\newline
and Trinity College\newline
Cambridge CB2 1TQ, England}
\email{bollobas@msci.memphis.edu\newline
\hglue3cm{b.bollobas@dpmms.cam.ac.uk}}

\address{Christian Borgs\newline
Microsoft Research\newline
One Microsoft Way\newline
Redmond, WA 98052}
\email{borgs@microsoft.com}

\address{Jennifer Tour Chayes \newline
Microsoft Research \newline
One Microsoft Way \newline
Redmond, WA 98052}
\email{jchayes@microsoft.com}

\address{Jeong Han Kim\newline
Microsoft Research\newline
One Microsoft Way\newline
Redmond, WA 98052}
\email{jehkim@microsoft.com}

\address{David Bruce Wilson\newline
Microsoft Research\newline
One Microsoft Way\newline
Redmond, WA 98052}
\email{dbwilson@microsoft.com}

\maketitle
\thispagestyle{empty}

\vspace{1mm}

\renewcommand{\thefootnote}{\roman{footnote}}

\centerline{
B\'ela Bollob\'as,$^{2,3}$ 
Christian Borgs,$^{1}$
Jennifer T. Chayes$,^{1}$
}
\centerline{
Jeong Han Kim,$^1$ and
David B. Wilson$^1$
}

\vspace{6mm}

\begin{centering}
{\it
$^1$Microsoft Research, Redmond, Washington\\
$^2$Department of Mathematical Sciences, University of Memphis\\
$^3$Trinity College, Cambridge, England\\
}
\end{centering}

\vspace{6mm}

\begin{center}
\longversion
\end{center}
\vspace{20mm}

{\small
\noindent{\bf Abstract.}
We consider the random 2-satisfiability problem, in which each
instance is a formula that is the conjunction of $m$ clauses
of the form $x\vee y$,
chosen uniformly at random from among all 2-clauses on
$n$ Boolean variables and
their negations.  As $m$ and $n$ tend to infinity in the
ratio $m/n\rightarrow\alpha$, the problem is known to have a phase transition
at $\alpha_c = 1$, below which the probability that the formula
is satisfiable tends
to one and above which it tends to zero.  We determine the finite-size
scaling about this transition, namely the scaling of the maximal window
$W(n,\delta) = (\alpha_-(n,\delta),
\alpha_+(n,\delta))$ such that the probability of satisfiability
is greater than $1-\delta$ for $\alpha < \alpha_-$ and is
less than $\delta$ for $\alpha > \alpha_+$.  We show that
$$
W(n,\delta)=(1-\Theta(n^{-1/3}),
1+\Theta(n^{-1/3})),
$$
where the constants implicit in $\Theta$ depend on $\delta$.  We
also determine the rates at which the probability of satisfiability
approaches one and zero at the boundaries of the window.
Namely, for $m=(1+\varepsilon)n$, where $\varepsilon$ may depend on
$n$ as long as $|\varepsilon|$ is sufficiently small and
$|\varepsilon| n^{1/3}$ is sufficiently large, we show that
the probability
of satisfiability decays
like
$
\exp\left(-\Theta\left({n\varepsilon^3}\right)\right)
$
above the window, and goes to one like
$
1-\Theta\left( n^{-1}|\varepsilon|^{-3}\right)
$
below the window.
 We prove
these results by defining an order parameter for the transition
and establishing its scaling behavior in $n$ both inside and outside
the window.  Using this order parameter,
we prove that the 2-SAT phase transition
is continuous with an order parameter critical exponent of 1.
We also determine the values of two other critical exponents,
showing that the exponents of 2-SAT are identical to those of
the random graph.

\vspace{3mm}


\vspace{6mm}

}


\newpage

\setcounter{page}{1}
\setcounter{footnote}{0}
\renewcommand{\thefootnote}{\arabic{footnote}}

\section{Introduction and Statement of Results}
\label{sec:intro}

There has recently been interest in
a new field emerging at the intersection of
statistical physics, discrete mathematics, and
theoretical computer science.   The field is
characterized by the study of phase transitions
in combinatorial structures arising in problems from
theoretical computer science.

Perhaps the most interesting phenomena in
statistical physics are phase transitions.
These transitions occur in systems with infinitely
many degrees of freedom, i.e.\ systems specified by
infinitely many random variables.  Physically, the transitions
represent changes in the state of the system; mathematically,
the transitions are manifested as nonanalyticities in
relevant functions of an external control parameter,
such as the temperature.  In systems with a large but
finite number of degrees of freedom, one can study
the approach to nonanalytic behavior.  This study
is called finite-size scaling.  In systems with
continuous phase transitions characterized by
critical exponents, the form of the finite-size
scaling turns out to be related to these exponents.

Discrete mathematics often focuses on the study
of large combinatorial structures.  Random versions
of these structures (with respect to natural distributions)
are discrete systems with large but finite numbers
of degrees of freedom.  In the limit of an infinite
number of degrees of freedom, these systems can and
often do undergo phase transitions.  The study
of threshold phenomena emerging in these large combinatorial
structures is therefore analogous to finite-size
scaling in statistical physics.

The theory of complexity focuses on the difficulty
of solving certain combinatorial problems which
arise naturally in theoretical computer science.
The complexity of a given problem is determined
by the difficulty of solving any instance
of the problem (i.e., in the worst case).  Researchers
have also studied
randomly chosen
 instances of certain problems, and determined average-
or typical-case complexity.  However, even when it
is determined that a problem is easy or hard on
average, it is not clear what properties characterize
the hard instances.

The convergence of these three disciplines is a consequence
of the recent observation that one can define control
parameters in terms of which certain theoretical
computer science problems undergo phase transitions,
and the even more interesting observation that the
hardest instances of these problems seem to be concentrated
at the phase transition point.  The problem for which this
phenomenon has been studied most extensively is the
$k$-satis\-fiability problem.  Our work is the first
complete, rigorous
analysis of finite-size scaling for a satis\-fiability
problem.

The $k$-satisfiability ($k$-SAT) problem is
a canonical constraint satisfaction problem
in theoretical computer science.  Instances of
the problem are formulae in conjunctive normal
form:  a $k$-SAT formula is a conjunction of $m$
clauses, each of which is a disjunction of
length $k$.  The $k$ elements of
each clause are chosen from among $n$ Boolean
variables and their negations.  Given a formula,
the decision version of the problem is whether
there exists an assignment of the $n$ variables
satisfying the formula.

It is known that the $k$-satisfiability problem
behaves very differently for $k = 2$ and $k \geq 3$
\cite{Coo71}.  For $k=2$, the problem is in P \cite{Coo71};
indeed, it can be solved by a linear time
algorithm \cite{APT79}.  For $k \geq 3$, the problem
is NP-complete \cite{Coo71}, so that in the worst case
it is difficult to
determine whether a $k$-SAT
formula is satisfiable or not --- assuming P $\neq$ NP.
Note,
however, that even for $k=2$, variants of the $k$-SAT
problem are difficult.  For example, the
MAX-2-SAT problem, in which one determines
whether the maximum number of satisfiable
clauses in a 2-SAT formula is bounded by
a given integer, is an NP-complete problem \cite{GJS76} (see also
\cite{GJ79}), and even approximating it to within a factor of
$4/3-\varepsilon$ is NP-hard \cite{Has97}.

More recently, it has been realized that---rather
than focusing on worst-case instances---it is
often useful to study typical
instances of
the fixed-$k$ problem
as a function of the parameter $\alpha=m/n$.
Consider the random $k$-SAT problem, in which
formulae are generated by choosing
uniformly at random from among all possible clauses.
As $m$ and $n$ tend to infinity with limiting ratio
$m/n\rightarrow\alpha$, considerable
empirical evidence suggests
that the random
$k$-SAT problem undergoes a {\em phase transition}
at some value
$\alpha_c(k)$ of the parameter $\alpha$
(\cite{MSL92},
\cite{CA93}, \cite{LT93}, \cite{KS94}):  For
$\alpha < \alpha_c$, a random formula is
satisfiable with probability tending to
one as $m$ and $n$ tend to infinity in the
fixed ratio $\alpha=m/n$, while if
$\alpha > \alpha_c$, a random formula is
unsatisfiable with probability tending to
one as $m$ and $n$ tend to infinity, again
with $m/n\rightarrow\alpha$.

Existence of the phase transition is on a
different footing for $k=2$ and $k \ge 3$.
For $k=2$, it was shown by Goerdt
(\cite{Goe92}, \cite{Goe96}),
Chv\'atal and Reed \cite{CR92}, and
Fernandez de la Vega \cite{Fer92}
that a transition occurs at
$\alpha_c(2) = 1$.
For $k \ge 3$, it may not be
possible to locate the exact value of the
transition point.  However,
there has been considerable work bounding
the value of the presumed 3-SAT threshold from
below and above.  Using a succession of
increasingly sophisticated and clever algorithms
for finding SAT solutions with high probability,
lower bounds on $\alpha_c(3)$ were improved
from 1 (\cite{CF86}, \cite{CF90}, \cite{CR92})
to 1.63 \cite{BFU93} to
3.003 \cite{FS96} to 3.145 \cite{Ach00} to 3.26 \cite{AS00}.
Bounding the probability
of finding a solution by the expected number
of solutions gave an upper bound on
$\alpha_c(3)$ of 5.191 \cite{FP83};
increasingly sophisticated
counting arguments gave
a succession of improved upper bounds on
$\alpha_c(3)$ from 5.08 \cite{EF95} to 4.758 \cite{KMPS95}
to 4.643 \cite{DB97} to 4.602 \cite{KKK96} to 4.596 \cite{JSV00}.
More recently a bound of 4.506 \cite{DBM99} has been announced.
Although these bounds are relatively tight,
they nevertheless allow for the possibility of a non-sharp
transition.
However, motivated by the empirical evidence,
Friedgut and later Bourgain showed that indeed there is a sharp
transition \cite{FB99} (although they did not prove that the
probability of satisfiability approaches
a limit).  These proofs were based
on a general argument which shows that
global, as opposed to local, phenomena
lead to sharp transitions.  However, the existence
of a limiting threshold is still an open problem.
\addtolength{\textheight}{1\baselineskip}

Having established the sharpness of the transition,
the next step is to analyze some of its properties.
{\em Finite-size scaling} is the study
of changes in the transition behavior due
to finite-size effects, in particular,
broadening of the transition
region for finite $n$.
To be precise, for $0 <\delta < 1$, let
$\alpha_-(n,\delta)$ be the supremum over
$\alpha$ such that for $m=\alpha n$, the
probability of a random formula being
satisfiable is at least $1-\delta$.
Similarly, let
$\alpha_+(n,\delta)$ be the infimum over
$\alpha$ such that for $m=\alpha n$, the
probability of a random formula being
satisfiable is at most $\delta$.
Then, for $\alpha$  within the {\em scaling window}
\begin{equation}
W(n,\delta) = (\alpha_-(n,\delta),
\alpha_+(n,\delta)),
\label{window}
\end{equation}
the probability of
a random formula being satisfiable is
between $\delta$ and $1-\delta$.
Since, by \cite{FB99}, for all $\delta$,
$|\alpha_+(n,\delta) - \alpha_-(n,\delta)| \to 0$
as $n\to \infty$,
we
say that the scaling window represents
the broadening of the transition due
to finite-size effects.  Sometimes we shall
omit the explicit $\delta$ dependence of
$\alpha_\pm(n,\delta)$ and $W(n,\delta)$, writing
instead $\alpha_\pm(n)$ and $W(n)$.  In these cases,
the power laws we quote will be uniform in $\delta$,
but the implicit constants may depend on $\delta$.

The first model for which such broadening was established rigorously
is the random graph model.  The phase transition for this model, namely
the sudden emergence of a giant component, was already proved by
Erd\H os and R\'enyi
(\cite{ER60}, \cite{ER61}).  But the characteristic width of the
transition was (correctly) investigated only 24 years later
by Bollob\'as \cite{Bol84} (see also \cite{Bol85} and the references therein).
In particular, this work showed that the width
of the scaling window $W(n)$ is $n^{-1/3+o(1)}$;
the precise growth rate was later shown
to be $\Theta(n^{-1/3})$ by \L uczak \cite{Luc90}.
Many additional properties of the phase transition were then determined
using generating functions
\cite{LPW94} \cite{JKLP94}.
For the finite-dimensional analogue of the
random graph problem, namely percolation on
a low-dimensional hypercubic lattice, the
broadening was established by
Borgs, Chayes, Kesten and Spencer
(\cite{BCKS98a}, \cite{BCKS98b}), who also
related the power law form
of $\alpha_{\pm}(n)$ to the critical exponents
of the percolation model.

The question of finite-size scaling in
the $k$-SAT model was first addressed
by Kirkpatrick and Selman \cite{KS94},
who presented both a
heuristic framework
and empirical evidence for analysis of
the problem.  There has also been
subsequent empirical (\cite{SK96},
\cite{MZKST99}) and theoretical
(\cite{MZ96}, \cite{MZ97}, \cite{MZKST99})
work, the
latter using the replica method familiar
from the study of disordered, frustrated
models in condensed matter physics
(see \cite{MPV87} and references therein).
Although the theoretical
work has yielded a good deal of insight,
the empirical work on finite-size scaling
has been misleading \cite{Wil00}, and
rigorous progress on finite-size scaling
in $k$-SAT has been quite limited.

\addtolength{\textheight}{-1\baselineskip}
In this work, we address the question
of finite-size scaling in the 2-SAT problem;
in particular, we obtain the power law form of
the
scaling window
$W(n) = (\alpha_-(n),
\alpha_+(n))$,
together with the rates of convergence at
the boundaries of the window.
Previous work on 2-SAT by Goerdt \cite{Goe99}
has shown that
$\alpha_-(n) \ge 1 - O(1/\sqrt{\log n})$,
while Verhoeven \cite{Ver99} has
recently obtained the result
$\alpha_+(n) \le 1+ O(n^{-1/4})$.
Numerical work on the scaling window for
2-SAT is somewhat controversial: While
earlier simulations \cite{MZKST99}
suggested that the
window scales like
$W(n) = (1 - \Theta (n^{-1/2.8}), 1+\Theta (n^{-1/2.8}))$,
recent simulations by Wilson \cite{Wil98} indicate that
the 2-SAT formulae considered in \cite{MZKST99} are not long
enough to reach the asymptotic regime.\footnote{As
usual, $f=\Theta(g)$ means that there exist
positive, finite
constants $c_1$ and $c_2$ such that $c_1 \leq f/g \leq c_2$.
Unless noted otherwise, these constants are universal.  In fact,
in the above formulae for $W(n)$, the constants depend on
$\delta$.}
Indeed,
we shall prove in this paper that
$W(n) = (1 - \Theta (n^{-1/3}), 1+\Theta (n^{-1/3}))$,
as conjectured earlier by
Bollob\'as, Borgs, Chayes and Kim \cite{BBCK98}
and predicted numerically in \cite{Wil98}.
We also show how the probability of satisfiability tends
to $1$ and $0$ at the edges of the window.

In order to state our results precisely, we
need a little notation.  Let $x_1, \dots, x_n$
denote
$n$ Boolean variables.  Writing $\overline x$ for the
{\em negation} of $x$, our
$n$ variables give $2n$ {\em literals}
$x_1, \dots, x_n,
\overline x_1, \dots, \overline x_n$.
Two literals $x$ and $y$ are said to be
{\em strictly distinct} if neither
$x=y$ nor $x=\overline y$.  A
{\em k-clause} is a
disjunction $C=u_1\vee \cdots \vee u_k$
of $k$
strictly distinct literals, and
a $k$-SAT formula is a conjunction
$F=C_1\wedge\dots\wedge C_m$ of
$k$-clauses
$C_1$, $\dots$, $C_m$.
We say that $H$ is a subformula
of $F$ if it can be obtained from
$F$ by deleting some of its clauses.
A $k$-SAT formula
$F=F(x_1,\dots,x_n)$ is said to be satisfiable,
or SAT, if there exists a truth assignment
$\eta_i\in\{0,1\}$, $i=1,\dots, n$,
such that $F(\eta_1,\dots,\eta_n)=1$.  Here, as usual,
$0$ stands for the logical value FALSE,
and $1$ is the logical value TRUE.
We write ``$F$ is SAT'' if the formula $F$
is satisfiable, and ``$F$ is UNSAT'' if the formula
$F$ is not satisfiable.  We also sometimes use the alternative
notation $\SAT(F)$ and $\UNSAT(F)$ to denote these two
cases.

We consider the random 2-SAT problem in two
essentially equivalent forms, given by
{\em a priori} different probability distributions
of random 2-SAT formulae on $x_1$, $\dots$, $x_n$.
First, we consider the probability space of formulae
$F_{n,m}$ chosen
uniformly at random from all 2-SAT formulae with
exactly $m$ different clauses.  (Here $x \vee y$
is considered to be the same as $y \vee x$, but
different from e.g.\ $x \vee \overline y$.)  Second,
we consider the space of formulae $F_{n,p}$
with 2-clauses on
$x_{_1}, \dots, x_{_n}$
chosen independently with probability $p$.  In this
introduction, we shall state theorems in terms of
the $F_{n,m}$; the equivalent theorems for the $F_{n,p}$
will be given in Section \ref{sec:strat}.   The
conversion between the two formulations of the
problem is given
in Appendix A.  In both cases, we use $\Pr(A)$
to denote the probability of an event $A$.

As usual in  2-SAT, it is convenient to study the
phase transition in terms of the parameter $\varepsilon$
representing the deviation of $\alpha$ from its
critical value:
\begin{equation}
\label{def-epsilon}
m = (1 + \varepsilon) n.
\end{equation}
When studying finite-size effects,
we shall take the parameter $\varepsilon$ to depend on $n$.
Our analysis shows that the appropriate scaling of
$\varepsilon$ is
$n^{-1/3}$, so that it is natural
to define yet another parameter
$\lambda= \lambda_n$ according to
\begin{equation}
\label{def-lambda}
\varepsilon = \lambda_n n^{-1/3},
\end{equation}
and distinguish the cases $\lambda_n$ bounded,
$\lambda_n \rightarrow \infty$ and $\lambda_n \rightarrow
-\infty$.

Our main result is the following theorem.

\begin{theorem}
\label{sat_{n,m}}
There are constants $\varepsilon_0$
and $\lambda_0$, $0<\varepsilon_0 < 1$,
$0<\lambda_0<\infty$, such that
\begin{equation}
\label{new1.1}
\Pr(F_{n,m} \; \text{\rm is SAT})=
\begin{cases}
1-\Theta\Big(\frac{1}{|\lambda_n|^3}\Big)
&\qquad\text{if $ -\varepsilon_0 n^{1/3} \leq \lambda_n \leq -\lambda_0$},
\\
\Theta(1)
&\qquad\text{if $ -\lambda_0 \leq \lambda_n \leq \lambda_0$},
\\
\exp\big(-\Theta\big({\lambda_n^3}\big)\big)
&\qquad\text{if $\lambda_0 \leq \lambda_n \leq \varepsilon_0 n^{1/3}$}.
\end{cases}
\end{equation}
\end{theorem}
\noindent
Note that the behaviors for $\lambda_n <0$
and $\lambda_n >0$ can be cast in the same form
by writing
$\Pr(F_{n,m} \; \text{\rm is SAT})
= 1-\Theta(|\lambda_n|^{-3})
 = \exp(-\Theta({|\lambda_n|^{-3}}))$.

Theorem \ref{sat_{n,m}} gives us the exact form of
the scaling window:

\begin{corollary}
\label{cor:window}
For all sufficiently small $\delta > 0$,
the scaling window (\ref{window}) is of the form
\begin{equation*}
W(n,\delta)=(1-\Theta(n^{-1/3}),
1+\Theta(n^{-1/3})),
\end{equation*}
where the constants implicit in the definition of $\Theta$
depend on $\delta$, and are easily calculated from equation
(\ref{new1.1}).
\end{corollary}

\noindent
Of course, the theorem gives us more than
the boundaries of the window; it also gives
us the rates of approach of the probability of satisfiability
at these boundaries.  As an easy special
case of the rate result at the upper boundary, note that
if $\varepsilon$ is positive and independent of
$n$, then our result for $\lambda > \lambda_0$
gives that
\begin{equation}
\Pr(F_{n,m} \; \text{\rm is SAT})
=\exp(-\Theta({\varepsilon^3}n)).
\end{equation}
This
strengthens both the result of Fernandez de la Vega \cite{Fer98}
that
$\Pr(F_{n,m} \; \text{is SAT})
=O(\exp(-{f(\varepsilon)\sqrt{n}}))$
and the recent improvement of Achlioptas and Molloy \cite{AM98}
that $\Pr(F_{n,m} \; \text{is SAT})
=O(\exp(-{f(\varepsilon)n}))$ for
some $f(\varepsilon) > 0$.

We remark that in the random graph model, the existence of a complex
component, i.e.\ a connected component with more than one cycle, is
roughly analogous to the existence of a contradiction in a random
2-SAT formula.  When there are $m=\frac12 n(1+\lambda n^{-1/3})$ edges,
Britikov \cite{Bri89} showed
\begin{equation}
\Pr\left(\parbox{1.4in}{$G_{n,m}$ contains no
\\ complex component}\right)=
(1+o_n(1))\begin{cases}
1-\frac{5+o_\lambda(1)}{24}\frac{1}{|\lambda|^3}
&\qquad\text{if $-\omega(n)\leq \lambda \leq -\lambda_0$},
\\
P(\lambda)
&\qquad\text{if $ -\lambda_0 \leq \lambda \leq \lambda_0$},
\\
\frac{\sqrt{2\pi}+o_\lambda(1)}{2^{1/4}\Gamma(1/4)}\frac{e^{-\lambda^3/6}}{\lambda^{3/4}}
&\qquad\text{if $\lambda_0 \leq \lambda \leq \omega(n)$},
\end{cases}
\end{equation}
where $P(\lambda)$ is an explicit power series in $\lambda$, and
$\omega(n)$ is an unspecified slowly growing function of $n$.  Later
this formula was shown to be valid for $|\lambda|\ll n^{1/12}$
\cite[pg.\ 289]{JKLP94}.  Observe that in comparison, the analogous
bounds for random 2-SAT in Theorem~\ref{sat_{n,m}} are not so precise,
but they hold for the full range of $|\lambda|\ll n^{1/3}$.

The key to our analysis is the introduction of an
{\em order parameter} for the 2-SAT phase transition.
As usual in statistical physics,
an order parameter is a function which
vanishes on one side of the transition and becomes
non-zero on the other side.  Control of the growth of the
order parameter was the key to Bollob\'as' analysis of
finite-size scaling in the random graph \cite{Bol84}, and to
Borgs, Chayes, Kesten and Spencer's
analysis of finite-size scaling in percolation
\cite{BCKS98b}.
Our order parameter for satisfiability is the
average density of the
{\em spine} of
a  Boolean formula, which we define as follows.
Given a formula $F$ in conjunctive normal form,
we define the spine $S(F)$ as the
set of literals
$x$ such that there is a satisfiable
subformula $H$ of $F$ for which $H \wedge x$
is not satisfiable,
\begin{equation}
S(F)=\{x\mid \exists H\subset F,
\text{$H$ is SAT
and $H\wedge x$ is UNSAT}\} .
\label{defspine}
\end{equation}

Our notion of the spine was motivated by the insightful concept
of the backbone, $B(F)$, introduced by
Monasson and Zecchina \cite{MZ96} ---
where the backbone density $|B(F)|/n$ was
originally called ``the fraction of frozen
variables.''
The backbone $B(F)$ is the set of literals that are required to
be FALSE in any assignment that minimizes the number
of unsatisfied clauses in $F$.  It is easy to see that
$B(F) \subset S(F)$, and in particular $B(F) = S(F)$
if $F$ is satisfiable.  One of the principal differences
between the spine and the backbone is that the spine is
monotone in the sense that adding clauses to a formula
only enlarges its spine.  It is the monotonicity which
enables us to achieve analytical control of the spine.
In addition, we have found that the spine is easier to
simulate than the backbone \cite{Wil98}.  We believe that
the spine will become an important tool in the analysis
of satisfiability problems.

Consider now a satisfiable 2-SAT formula $F$.
It is not hard to see that the addition of the
 2-clause
$C=x\vee y$ makes $F$ (or, more precisely, makes $F\wedge C$)
unsatisfiable if and only if both $x$ and $y$ lie in the
spine.  Building a random 2-SAT formula
by adding clauses one by one at random to an initially empty (and
hence satisfiable) formula, we can therefore control
the probability that a formula is satisfiable if we
have sufficient control of the spine
in each step.  This is the strategy we shall
follow to prove Theorem~\ref{sat_{n,m}}.

In the course of proving
Theorem~\ref{sat_{n,m}}, we obtain detailed estimates
on the expectation and variance of the size of the spine
inside the scaling window
$m\in [n-\Theta(n^{2/3}), n+\Theta(n^{2/3})]$,
i.e.\ the finite-size scaling of the spine.  Before stating
these results, however, let
us give the behavior of
the size of $S(F_{n,m})$ on the scale $n$.  To this end,
let
$
\vartheta:(0,\infty)\to(0,1)$
be the function satisfying
\begin{equation}\label{def-theta}
1-\vartheta(\varepsilon)
= \exp[-(1+\varepsilon)\vartheta(\varepsilon)],
\end{equation}
i.e.
\begin{equation}
\label{alternative-theta}
\vartheta(\varepsilon) =
1 - \sum_{k=1}^\infty
\frac{k^{k-1}}{k!}(1+\varepsilon)^{k-1} e^{-(1+\varepsilon)k}.
\end{equation}
The $k$th term in \eqref{alternative-theta} is the probability
that a Poisson birth-and-death process with birth rate $1+\varepsilon$
will have size $k$, while $\vartheta(\varepsilon)$ is the probability
that it is infinite.
Note that $\vartheta(\varepsilon) =2\varepsilon+O(\varepsilon^2)$
for positive $\varepsilon$ sufficiently small.
The size of the spine is given by:

\begin{theorem}
\label{thm1.2}
For any fixed $\varepsilon\in(-\varepsilon_0,\varepsilon_0)$,
where $\varepsilon_0$ is the constant from
Theorem~\ref{sat_{n,m}}, we have
\begin{equation}
\label{1.6}
E(|S(F_{n,m})|)=
\begin{cases}
\Theta(\varepsilon^{-2})
&\text{if $ \varepsilon < 0$}
\\
\Theta(n^{2/3} )
&\text{if $ \varepsilon = 0$}
\\
2 n\vartheta(\varepsilon)+o(n)
&\text{if $ \varepsilon > 0$}.
\end{cases}
\end{equation}
\end{theorem}

The behavior above, coupled with
the role of the spine in
the proof of Theorem~\ref{sat_{n,m}}, justifies our
identification of
the density of the spine as an order parameter for
the 2-SAT transition. In the language of phase
transitions, Theorem \ref{thm1.2}
implies that the 2-SAT transition
is second-order (or continuous), with
order parameter critical exponent $\beta=1$.
Here, as usual, we say that the order parameter
has critical exponent $\beta$ if
$\lim_{n\to\infty}E(|S(F_{n,m})|)/n = \Theta(\varepsilon^\beta)$
as $\varepsilon \downarrow 0$, see discussion following
Remark~\ref{rem-variance}.

The next theorem states our results for the finite-size
scaling of the spine $S(F_{n,m})$:

\begin{theorem}
\label{thm1.3}
Let $\varepsilon_0$ and $\lambda_0$ be the constants in
Theorem~\ref{sat_{n,m}}.
Suppose $|\lambda_n| \leq \varepsilon_0 n^{1/3}$.  Then
\begin{equation}
E(|S(F_{n,m})|) =
\label{new1.8}
\begin{cases}
\frac{1}{2} \lambda_n^{-2}{n^{2/3}}
 (1 + o(1))
&\quad\text{if $ \lambda_n < -\lambda_0$}
\\
\Theta(n^{2/3} )
&\quad\text{if $ |\lambda_n| \leq \lambda_0$}
\\
4\lambda_n n^{2/3}(1+o(1))
&\quad\text{if $ \lambda_n > \lambda_0$},
\end{cases}
\end{equation}
where the $o(1)$ terms represent errors which go to
zero as $|\lambda_n|\to \infty$ and $\varepsilon=\lambda_n n^{-1/3} \to 0$.
\end{theorem}

\begin{remark} \label{rem-variance}
In the course of proving Theorems
\ref{sat_{n,m}}, \ref{thm1.2} and \ref{thm1.3},
we shall prove bounds on the variance of $|S(F_{n,m})|$
which allow us to generalize the above statements in expectation
to statements in probability.
\end{remark}

Statistical mechanical models with second-order (i.e., continuous)
transitions are often characterized by critical exponents which
describe the behavior of fundamental quantities at or approaching
the critical point.  It turns out (see \cite{BCKS98b} and announcements
in \cite{Cha98} and \cite{CPS99}) that it is possible to read off
some of these exponents from the finite-size scaling form of the order
parameter and the scaling window.  In particular, the scaling of
the order parameter at the critical point allows us to evaluate
the so-called field exponent $\delta$ as
\begin{equation}
\label{exponent-delta}
E(|S(F_{n,m})|) = \Theta\big(n^\frac{\delta}{1+\delta}\big)
\quad\quad \text{if}\quad |\lambda_n| < \lambda_0.
\end{equation}
Similarly (again using \cite{BCKS98b}, \cite{Cha98} and \cite{CPS99}),
the scaling of the window allows us to identify the exponent sum
$2\beta + \gamma$,
according to
\begin{equation}
\label{exponent-sum}
W(n,\delta)=\big(1-\Theta(n^{-1/(2\beta+\gamma)}),
1+\Theta(n^{-1/(2\beta+\gamma)})\big),
\end{equation}
where $\beta$ is the order parameter exponent
described above, i.e.\
\begin{equation}
\label{exponent-beta}
\lim_{n\to\infty}\frac{1}{n}E(|S(F_{n,m})|) = \Theta(\varepsilon^\beta)
\quad\quad \text{as}\quad \varepsilon \downarrow 0,
\end{equation}
and $\gamma$ is the so-called susceptibility exponent.
Comparing equations \eqref{exponent-beta}, \eqref{exponent-sum}, and
\eqref{exponent-delta} to Theorem \ref{thm1.2}, Corollary \ref{cor:window},
and Theorem \ref{thm1.3}, respectively, we get the following.

\begin{corollary}
\label{exponents}
The 2-SAT transition is a second-order (i.e.\ continuous) transition
with critical exponents:
\begin{equation*}
\beta = 1, \quad \gamma = 1, \quad and \quad \delta = 2.
\end{equation*}
\end{corollary}
\noindent
Thus we have proved that the critical exponents of the random
2-SAT problem are identical to those of the random graph.
See \cite{BBCKW00} for a more detailed discussion of the
critical exponents for 2-SAT.

The organization of this paper is as follows.  In Section
\ref{sec:digraph}, we discuss the well-known representation of 2-SAT
formulae as directed graphs, a representation we use
extensively in our proofs.
In that section, we also derive new results on various
representations of
the spine in terms of directed graphs.
While most of the results in Section
\ref{sec:digraph} concern given formulae, not distributions
of formulae, a final result there gives a mapping of a
distribution of certain sets in the graphical representation of
random
2-SAT into the standard random graph model.
In Section \ref{sec:strat},
we state our main technical estimates on the
expectation and variance of the size of the spine, and formulate
an analogue of Theorem~\ref{sat_{n,m}} for the distribution
$F_{n,p}$.  We then outline the strategy of our proof, giving
first our heuristic for the expected size of the spine, and then
showing how this will be used to obtain the size of the scaling
window (Theorem~\ref{sat_{n,m}}).
While the width of the scaling window can be determined from the spine
expectation and variance estimates alone, the rate of approach from
above in Theorem~\ref{sat_{n,m}} requires that a sufficiently large
spine forms with extremely high probability.  In order to prove this,
in Section \ref{sec:sat}, we define structures we call ``hourglasses''
which are basically precursors to the spine, and we state a theorem
giving conditions under which a giant hourglass forms.  The
proof of the hourglass theorem is given in Section \ref{sec:hourglass}.
In Section \ref{sec:sat}, we use the expectation and variance results on
the spine, and the hourglass theorem, to establish the analogue of
Theorem~\ref{sat_{n,m}} for the distribution $F_{n,p}$.  In
Sections \ref{sec:prelim} and \ref{sec:moments},
we develop some machinery from random graph theory
and derive moment bounds, which enable us to prove the
expected size and variance results for the spine in Sections
\ref{sec:spine}
and \ref{sec:var}, respectively.   Appendix A contains the conversion
from our results on $F_{n,p}$ to $F_{n,m}$,
and Appendix B establishes a technical result on the cluster size
distribution in the random graph problem.

\bigskip

\section{The Digraph Representation of 2-SAT}
\label{sec:digraph}

In the digraph representation, each
2-SAT formula corresponds to
a certain directed graph (or digraph) $D_F$.
To motivate the mapping of $F$ into $D_F$, note
that $F$ is satisfiable if and
only if all
clauses in $F$ are satisfiable.  Thus, if $F$
contains a clause $C=x\vee y$, a
satisfying truth assignment with $x$
set to FALSE
requires that $y$ is set to TRUE,
and a satisfying assignment with $y$ set to
FALSE requires that $x$ is set to TRUE.
So the clause $x\vee y$
corresponds to the logical implications
$\overline x=\TRUE\implies y=\TRUE$ and
$\overline y=\TRUE\implies x=\TRUE$.  We shall encode
this fact in the digraph $D_F$ by including
the edges $\overline x\rightarrow y$
and $\overline y\rightarrow x$ in $D_F$ iff
$F$
contains a clause $C=x\vee y$.

To be
precise, given a 2-SAT formula $F$, define
the digraph $D_F$
as the directed graph with vertex set%
\footnote{Note that we deviate from the standard
notation, where $[n]$ stands for the set
$\{1,2,\dots,n\}$.}
\begin{equation}
\label{1.15}
[n]=\{x_1,\dots,x_n,\overline x_1,\dots,\overline x_n\}
\end{equation}
and edge set
\begin{equation}
\label{1.16}
E_F=\{x \rightarrow y\mid (\overline x\vee y)
\text{ is a clause in } F\}.
\end{equation}
Since $(x\vee y)$ and $(y\vee x)$ are considered to
be the same clause, the digraph $D_F$ contains
the edge
$\overline x \rightarrow y$ if and only if it contains
the edge $\overline y \rightarrow x$.
As usual, an oriented path in $D_F$ is
a sequence of vertices
$v_0,v_1, \dots, v_k\in [n]$
and edges $v_i\rightarrow v_{i+1}$
for $i=0,1,\dots, k-1$.
We say that this path is a path from $x$ to $y$ if $v_0=x$ and
$v_k=y$.  We write $x{\underset{D_F}\rightsquigarrow} y$,
or sometimes simply $x\rightsquigarrow y$, if
$D_F$ contains an oriented path from $x$ to $y$.
By convention, we shall say $x{\rightsquigarrow} x$
for all $x$.
Finally, we say that $D_F$ contains a {\em contradictory cycle}
if $x\rightsquigarrow \overline x$ and $\overline x\rightsquigarrow x$
for some $x\in [n]$.

The following lemma connecting the structure of the digraph
$D_F$ with the satisfiability of the formula $F$ is implicit
in all digraph analyses  of 2-SAT,  see e.g.\
 \cite{Goe92}.  For completeness, we shall give an explicit proof here.

\begin{lemma}
\label{lem:sat-digraph}
A  2-SAT formula $F$
is satisfiable if and only if
the  digraph $D_F$ has no contradictory cycle.
\end{lemma}

\begin{proof}
Let us
first assume that $F$ is satisfiable, with satisfying assignment
$\eta_i\in\{0,1\}$, $i=1,\dots,n$.  Consider an edge $\overline x
\rightarrow y$ in the corresponding digraph.  Since
$F(\eta_1,\dots,\eta_n)=1$, the presence of the edge $\overline x
\rightarrow y$ gives the logical implication
$x=\FALSE \implies y=\TRUE$. A
contradictory cycle $x\rightsquigarrow \overline x\rightsquigarrow x$
therefore gives the logical implication
$x=\TRUE \implies x=\FALSE \implies x=\TRUE$,
which is not compatible with any truth assignment for $x$.

We prove the converse by induction on the number $n$ of
variables.  For $n = 1$ there is nothing to prove.  Turning
to the induction step, suppose that the digraph $D_F$ has
no contradictory cycles.  We claim that in this case
$F$ is satisfiable.  To this end, we first recall the
definition of strongly connected components for directed
graphs.  We say that two vertices $x$ and $y$ in a directed
graph are {\it strongly connected}
if $x\rightsquigarrow y\rightsquigarrow x$, i.e.\ if the
directed graph $D_F$ has a cycle containing
$x$ and $y$.  The strongly connected component
of a vertex $x$ is the induced subgraph of $D_F$
containing the set of vertices
\begin{equation}
\label{1.17}
C_S(x)=\{y\mid x\rightsquigarrow y\rightsquigarrow x\}.
\end{equation}
Somewhat loosely, we call
$C_S(x)$ the {\it strong component} of $x$.
Clearly, the strong component partitions the vertex
set $[n]$.  We define a partial order $\leq$
on the set of all strong components by taking
$C_S(x) \leq C_S(y)$ if $x \rightsquigarrow y$,
and so $x^\prime \rightsquigarrow y^\prime$
for all $x^\prime \in C_S(x)$ and
$y^\prime \in C_S(y)$.
 Let $C_S$ be a minimal element in this
partial order, i.e.\ let $C_S$ be
a strong component such that $D_F$ contains
no edge $x \rightarrow y$ with $x \notin C_S$
and $y \in C_S$.  For a set of literals $M$,
let
\begin{equation}
\label{setbar}
\overline M=\{y\mid \overline y\in M\}.
\end{equation}
Since $D_F$ has no contradictory cycle,
$C_S \cap \overline{C_S} = \emptyset$.
Furthermore, since $C_S$ is a minimal element
in our partial order,
$\overline{C_S}$ must be a maximal
element.
If we set all literals in $C_S$ to
FALSE, and so all literals in $\overline{C_S}$
to TRUE, then all clauses in $F$ containing
at least one literal from $C_S \cup \overline{C_S}$
are TRUE.  This process removes all the variables
corresponding to literals in $C_S$ and $\overline{C_S}$
from $[n]$, and all clauses involving these
variables from $F$, leading to a new 2-SAT
formula $F^\prime$.  Since the graph
$D_{F^\prime}$ is a subgraph of $D_F$,
it contains no contradictory cycles either.
Using the inductive hypothesis, we obtain a satisfying
assignment for $F$, which completes the proof
of the converse and hence of the theorem.
\end{proof}

\begin{remark}
If $F$ is a mixture of one and two-clauses,
i.e.\ if it is of the form
$F=H\wedge x_1\wedge\dots\wedge x_k$
where $H$ is a 2-SAT formula and
$x_1,\dots,x_k$ are literals, we define $D_F$
by including the edges $\overline x_i\to  x_i$,
$i=1,\dots,k$, in
addition to the edges in $D_H$.  It is not
hard to see that the above proof applies also
to this situation,  giving again that $F$ is
SAT if and only if $D_F$ contains no contradictory
cycles.
\end{remark}

While the previous lemma says that contradictions
in a formula correspond to cycles in the digraph,
the next lemma says that the spine of a formula
corresponds to ``half-cycles'' in the digraph.
This graphical description of the spine
is central to our analysis.

\begin{lemma}
\label{lem:spine-digraph}
For every 2-SAT formula $F$,
\begin{equation}
\label{1.12}
S(F)=\{x\mid x{\underset{D_F}\rightsquigarrow}\overline x\},
\end{equation}
where $D_F$ is the digraph corresponding to $F$.
\end{lemma}

\begin{proof}
Suppose that $x \rightsquigarrow \overline x$, and let
$x=v_0\to v_1 \to
\cdots \to v_{r-1} \to v_r=\overline x$
be a shortest directed path from $x$ to
$\overline x$ in $D_F$. Then no literal appears twice in the path,
although a literal and its negation may well do so.
Let $\ell$ be the smallest positive integer such that
$v_\ell$ is not strictly distinct from all of $v_0, v_1,
\dots, v_{\ell-1}$, and let $0\le k < \ell$ be such that
$v_\ell=\overline v_k$.
Then
 $H=(\overline v_0 \vee v_1)\wedge
(\overline v_1\vee v_2) \wedge \cdots \wedge
(\overline v_{\ell-1} \vee v_\ell)$ is
a subformula of $F$ that is satisfied by setting each of
$v_0, v_1, \dots , v_{\ell-1}$ to FALSE.
On the other hand, $H\wedge x=H\wedge v_0$ is UNSAT since
in order to satisfy it, we would have to set $v_0$ to TRUE,
then $v_1$ to TRUE, and so on, ending with the requirement
that $v_\ell$ be set TRUE.  However, as $v_k$ is
TRUE, $v_\ell$ is already set FALSE.
This completes the proof that
$x\rightsquigarrow\overline x$
implies that $x\in S(F)$.

Conversely, suppose that $H \subset F$ is SAT and
$H'=H\wedge x= H\wedge (x\vee x)$
is UNSAT. Then $D_{H'}$ has a contradictory cycle
$C=u\rightsquigarrow \overline u \rightsquigarrow u$. Since $D_H$
does not have a contradictory cycle, the cycle $C$ of $H'$
contains the oriented edge $\overline x \to x$, say
$u \rightsquigarrow \overline u \rightsquigarrow \overline x
\rightarrow x \rightsquigarrow u$. But then in $D_H$ we
have $x \rightsquigarrow u \rightsquigarrow \overline u
\rightsquigarrow \overline x$, so $x \rightsquigarrow \overline x$.
Hence if $x\in S(F)$ then
$x{\underset{D_F}\rightsquigarrow} \overline x$.
\end{proof}

Our next lemma gives an alternative representation for the spine
of a 2-SAT formula $F$.  In order to state it, we
introduce the out-graph $D_F^+(x)$ of a vertex $x$ in $D_F$
as the set of vertices and edges that can be reached
from $x$.  $D_F^+(x)$ therefore has the vertex set
\begin{equation}
\label{outset}
L^+(x)=L^+_F(x) = \{y \mid x{\underset{D_F}\rightsquigarrow} y \},
\end{equation}
and contains all edges $y\to z$ in $D_F$ such that
$y\in L^+_F(x)$.
For future reference, we also introduce the in-set
\begin{equation}
\label{inset}
L^-(x)=L^-_F(x) = \{y \mid y{\underset{D_F}\rightsquigarrow} x \}
\end{equation}
and the corresponding in-graph  $D_F^-(x)$.
Note that $x \in L_F^\pm (x)$ since, by our convention,
$x \rightsquigarrow x$ for all $x$.

As we shall see,
the spine of a 2-SAT formula $F$ can equivalently be
described as the set of literals $x$ such that
$L^+_F(x)$ is not strictly distinct, where for
simplicity, we say that a set $M\subset [n]$ is {\em
strictly distinct} (s.d.)
if the literals in $M$ are pairwise strictly distinct.

\begin{lemma}
\label{lem:spine-s.d.}
For every 2-SAT formula $F$
\begin{align}
\label{spine-s.d.}
\{x\mid {\underset{D_F}
{x\rightsquigarrow \overline x}}\}
&=\{x\mid L^+_F(x)\; \text{\rm is not s.d.}\}
\notag
\\
&=
\{x \mid L_F^+(x)\setminus\{x,\overline x\}\;\text{\rm is not s.d.}\}.
\end{align}
\end{lemma}

\begin{proof}
We start with the first equality in (\ref{spine-s.d.}).
If $x \rightsquigarrow \overline x$, then
$\{x,\overline x\}\subset L_F^+(x)$, so $L_F^+(x)$ is not strictly
distinct.  If $L_F^+(x)$ is not strictly
distinct, then $\{y,\overline y\}\subset L_F^+(x)$
for some literal $y\in [n]$, and hence
$x \rightsquigarrow y$
and
$x \rightsquigarrow \overline y$.
But $x \rightsquigarrow \overline y$
implies that $y \rightsquigarrow \overline x$,
which together with $x \rightsquigarrow y$
implies $x \rightsquigarrow \overline x$.

To prove the second equality, we first note
that
the set of literals $x$ for which
$L_F^+(x)\setminus\{x,\overline x\}$ is not strictly
distinct is obviously a subset of
the set of literals $x$ such that
$L_F^+(x)$ is not strictly
distinct.  We are thus left with the
proof that the statement that
$L_F^+(x)$ is not strictly
distinct implies the (apparently stronger)
statement that $L_F^+(x)\setminus\{x,\overline x\}$
is not strictly distinct.  So let us assume
that $L_F^+(x)$ is not strictly
distinct.  By the first equality in
(\ref{spine-s.d.}), this implies
$x \rightsquigarrow \overline x$.
Since the digraph of a 2-SAT formula
does not contain
any direct edges from $x$ to $\overline x$, we
conclude that there must be a literal
$y$ strictly distinct from
$x$ such that
$x \rightsquigarrow y \rightsquigarrow\overline x$.
The latter statement implies that
both $x \rightsquigarrow y$
and $x \rightsquigarrow \overline y$,
so that $L_F^+(x)\setminus\{x,\overline x\}$
is not strictly distinct.
\end{proof}

\begin{remark}
As the above  proof shows,  the first equality in
(\ref{spine-s.d.}) is true for mixed formulas of
1- and 2-SAT clauses as well.  The second
is obviously false for mixed formulas of
1- and 2-SAT clauses, as the simple example of
the formula $F=\overline x$ shows.
\end{remark}

\subsection*{\sc The Trimmed Out-Graph}

We end this section with a construction of
a trimmed version of the out-graph $D_F^+(x)$,
which we denote by $\widetilde D_F^+(x)$
with vertex set denoted by $\widetilde L_F^+(x)$.
The utility of this trimmed graph is that, by projecting it to an
unoriented graph, we shall be able to relate it to the more familiar
random graph.  Given any digraph on a subset of
$\{x_1,\ldots,x_n,\overline x_1, \ldots, \overline x_n\}$ we can
project it to an unoriented graph on a subset of $\{x_1,\ldots,x_n\}$
by dropping negations.  In particular, each literal
$x\in\{x_1,\ldots,x_n,\overline x_1, \ldots, \overline x_n\}$ gets
mapped to its corresponding variable $\langle x\rangle\in
\{x_1,\ldots,x_n\}$, and each clause $x\vee y$ gets mapped to the edge
$\{\langle x\rangle,\langle y\rangle\}$.  We call this the unoriented
projection of the digraph.  Specifically, for $F=F_{n,p}$, we shall
compare the distribution of $\widetilde D_F^+(x)$ for a fixed vertex
$x$ to that of the connected component of a given vertex in the random
graph $G_{n,2p-p^2}$, where, as usual, $G_{n,\widetilde p}$ denotes
the random graph on $\{x_1, \dots, x_n\}$ that is obtained from the
complete graph on $\{x_1, \dots, x_n\}$ by keeping each edge with
probability $\widetilde p$.  We use the symbol $C_{n,\widetilde p}(x)$
to denote the connected component of the vertex $\langle x\rangle$ in
$G_{n,\widetilde p}$.

\smallskip
\noindent {\it Construction of the trimmed out-graph.}

We construct the trimmed out-graph $\widetilde D_F^+(x)$ by doing a
local search in $D_F$ starting from literal $x$, and at the same time
we construct that portion of the random graph $G_{n,2p-p^2}$ which
determines the connected component of vertex $\langle x \rangle$.
Let the ``current
graph'' be that subgraph of $D_F$ which consists of the vertices and
edges that have been examined by the local search.  The ``frontier''
consists of those vertices of the current graph from which further
searching may be done.  Initially the current graph consists of just
the literal $x$, and $x$ is in the frontier.  Eventually the frontier
will be empty, terminating the local search,
at which point $\widetilde D_F^+(x)$
will be defined to be the current graph.  During the search,
certain edges $v\rightarrow w$ will be tested to see if they are in $D_F$,
and the search records whether the results are ``yes'' or ``no''
on the corresponding unoriented edge
$\langle v\rangle\sim\langle w\rangle$.
These test results will later be used to construct the random
graph $G_{n,2p-p^2}$.  Each step in the local search consists of the
substeps listed below.

\begin{enumerate}
\item An arbitrary literal $v$ in the frontier is selected (one
choice is the lexicographically smallest).
\item For each literal $w$ such that neither $w$ nor $\overline w$ is in
the current graph, check if $v\rightarrow w$ is in $D_F$, and record
either ``yes'' or ``no'' on the edge
$\langle v  \rangle\sim \langle w  \rangle$ accordingly.

\item For each literal $w$ for which ``yes'' was recorded, declare
$w$ to be a ``new literal'' --- unless ``yes'' is recorded for both
$w$ and $\overline w$, in which case we declare only one
(say the unnegated one)
of them to be a ``new literal.''
\item Adjoin each new literal $w$ and the edge $v\rightarrow w$ to the
current graph.
\item Adjoin each new literal $w$ to the frontier.  Remove $v$ from the
frontier.
\item Consider each ordered pair of vertices $(w,f)$ such that either
(1) $w$ is new and $f$ is in the frontier but not new, or (2) $w$ and
$f$ are both new, and $w$ is lexicographically smaller.
Test if
$w \rightarrow f$ or $f\rightarrow w$ in $D_F$, and record
either ``yes'' or ``no'' on the edge
$\langle w\rangle\sim\langle f\rangle$ accordingly.
If there is one ``yes,'' adjoin the corresponding edge to the current
graph, if there are two ``yes'''s, adjoin only one of the edges, say
the edge $w\rightarrow f$.
\end{enumerate}

\begin{lemma}
The trimmed out-graph $\widetilde D_F^+(x)$ defined above
has the following properties.
\begin{enumerate}
\item[i)] $\widetilde D_F^+(x)$ is a subgraph of $D_F^+(x)$.
\item[ii)] $\widetilde L_F^+(x)$ is strictly distinct.
\item[iii)] $\widetilde L_F^+(x) = L_F^+(x)$ if and only if
$L_F^+(x)$ is
strictly distinct.
\item[iv)] For $F=F_{n,p}$, the unoriented projection of the digraph
$\widetilde D_F^+(x)$
has the same distribution as
$C_{n,2p-p^2}(x)$.
In particular,
$|\widetilde L_{F_{n,p}}^+(x)|$ and
$|C_{n,2p-p^2}(x)|$ are equidistributed.
\end{enumerate}
\label{trimmed-graph}
\end{lemma}

\begin{proof}
By construction, properties (i) and (ii) are obvious.  Property (iii)
is not much more difficult.  There are certain possible edges leading
out of the vertex set $\widetilde L_F^+(x)$ that were never tested, or
that were tested and present, but then excluded from the trimmed
out-graph $\widetilde D_F^+(x)$ anyway.  But each such edge either led
to a literal already in $\widetilde L_F^+(x)$, or else led to a literal
whose complement was in $\widetilde L_F^+(x)$.  Thus if the literal set
$L_F^+(x)$ were to contain more literals than $\widetilde L_F^+(x)$, then
$L_F^+(x)$ would not be strictly distinct.
On the other hand, if $L_F^+(x)$ and $\widetilde L_F^+(x)$
are identical, then $L_F^+(x)$ is trivially strictly distinct
by property (ii).

Property (iv) is similarly easy.  First, for each literal
$u\in[n]$, we define $[u]=\{u,\overline u\}$.
By induction we shall prove that, at the beginning and
end of each step of the search, the following properties hold:
\begin{enumerate}
\item For
every pair of literals $u$ and $v$ of the current graph, precisely two
edges between $[ u]$ and $[ v]$ have been tested.
\item   For every literal $v$
in the current graph but not in the frontier, and every literal $w$
such that neither $w$ nor $\overline w$ is in the current graph, precisely
two edges between $[ v]$ and $[ w]$ have
been tested, both results being
``no.''

\item   For every literal $v$
in the frontier, and any literal $w$
such that neither $w$ nor $\overline w$ is in the current graph,
none of the edges between $[ v]$ and $[ w]$ have
been tested.

\item If none of $u,\overline u, v, \overline v$ are in the current graph,
then no edges between $[ u]$ and $[ v]$ have
been tested.

\item For any pair of strictly distinct literals $u,v \in [n]$,
either none or precisely one of the four edges between $[u]$ and $[v]$
appears in the current graph.  The latter happens if and only if some
test between $[u]$ and $[v]$ was positive.

\end{enumerate}

Indeed, assume that (1) -- (5) hold at the beginning of a step.
To prove that (1) holds at the end of the step, we first note
that no edge between $[ u]$ and $[ v]$
was tested in the current step if neither $u$ nor $v$ is new.
If $v$ is old and $u$ is new, then either $v$ was the
selected vertex in the frontier, in which case the edges
$v\to u$ and $v\to \overline u$ have been tested in the current
step, or $v$
was not in the frontier, in which case
precisely two edges between
$[ u]$ and $[ v]$
were tested in a previous step (with answer ``no'')
by the inductive assumption (2).  If both $v$ and $u$
are new, then no edge between
$[ u]$ and $[ v]$
was tested in a previous step by the inductive assumption
(3), and precisely two edges
(the edges $u\to v$ and $v\to u$)
between
$[ u]$ and $[ v]$
are tested in the current step.

To prove (2), we note that if $v$ is in the current graph
but not in the frontier, it was in the frontier
in some previous step, and got removed from
the frontier after all edges from $v$
to  vertices $u$, with neither $u$ nor $\overline u$
in the current graph at the time,
were tested. This includes
in particular the vertex $w$ in question, and since
we assume that neither $w$ nor $\overline w$
is in the current graph, it follows that both tests
must have given the result ``no'' at the time.
After that step, $v$ is not in the frontier,
so no edge containing $v$
or $\overline v$ is ever tested again,
implying statement (2).

Statement (3) follows from the observation
that an edge between a vertex $v$
in the current graph and a vertex
$w$ such that neither $w$ nor $\overline w$
is in the current graph is only tested
if $v$ is the selected vertex in the
current step, in which case it is not
in the frontier after this step anymore.

Statement (4) is obvious, since an edge
$f\to w$ is only tested if either $f$ is
in the frontier (and hence in the current
graph before the current step), or
both $f$ and $w$ are new vertices, which means they
are in the current graph after steps (1) -- (6).

To prove (5), we consider three cases.  In the first
case, none of
the vertices $u$, $\overline u$, $v$ and $\overline v$ is in
the current graph, in which case no edge between
$[ u]$ and $[ v]$ appears
in the current graph by the inductive assumption (4).
The second case is
the one in which exactly one of the
four vertices $u$, $\overline u$, $v$ and $\overline v$
is in the current graph.  Without loss of generality,
let us assume that this is the vertex $v$.
Then none of the edges between
$[ u]$ and $[ v]$ appears
in the current graph by (2) and (3).  The third
case is that precisely two of the
four vertices $u$, $\overline u$, $v$ and $\overline v$
are in the current graph, say $u$ and $v$.
Then precisely two of the four edges
between $[ u]$ and
$[ v]$ have been tested by the inductive
assumption (1).  Since the above search procedure always
tests two of the four edges between  $[ u]$ and
$[ v]$ at a given time, and adds one (but not
both) of them precisely when at least one of them tests
positive, we get (4).

We now use the properties (1) -- (5) above to prove statement
(iv) of the lemma.
If we pick the unordered pairs of numbers between
$1$ and $n$ in some
arbitrary order, each time randomly saying ``present'' (with
probability $2p-p^2$) or ``absent'' (with probability $(1-p)^2$),
then even if the order in which we pick the pairs depends on the
previous random choices of present/absent, the result will be the
random graph $G_{n,2p-p^2}$.  This is in effect what the
trimmed local
search does, except that it stops when
the connected component
containing $\langle x \rangle$ has been determined.
Thus unoriented projection of $\widetilde D_F^+(x)$
is just the connected component
containing $\langle x \rangle$ in $G_{n,2p-p^2}$.
\end{proof}

\section{Strategy of the Proof}
\label{sec:strat}
In this section, we shall first state our principal estimates
and results for
the distribution $F_{n,p}$ (to be proved in later sections), and
then give the heuristics for these results.

\subsection{\sc Main Results for the Distribution $F_{n,p}$}

As explained in the last section, the spine of a formula $F$
consists of all literals $x$ for which $x \rightsquigarrow \overline
x$ (see Lemma~\ref{lem:spine-digraph}),
which in turn
is just the set of all literals $x$ such that
$L^+_F(x)$ is strictly distinct
(see Lemma~\ref{lem:spine-s.d.}).
If $F$ is distributed according to the model
$F_{n,p}$, the
expectation and variance of the size of $S(F_{n,p})$ are
therefore
given by the equations
\begin{equation}
\label{2.4}
E(|S(F_{n,p})|)=
\sum_{x\in [n]}
\Pr(x\underset {F_{n,p}}\rightsquigarrow\overline x)
\end{equation}
and
\begin{equation}
\label{2.5}
E(|S(F_{n,p})|^2)-E(|S(F_{n,p})|)^2
=\sum_{x,y\in [n]}
\Pr
\Big(
x\underset {F_{n,p}}\rightsquigarrow\overline x
\text{ and }
y\underset {F_{n,p}}\rightsquigarrow\overline y
\Big)
-
\Pr
\Big(
x\underset {F_{n,p}}\rightsquigarrow\overline x
\Big)
\Pr
\Big(
y\underset {F_{n,p}}\rightsquigarrow\overline y
\Big),
\end{equation}
where $x\underset {F_{n,p}}\rightsquigarrow\overline x$
a shorthand for
$x\underset {D_{F_{n,p}}}\rightsquigarrow\overline x$.

The following two theorems allow us to prove suitable
bounds on the expected size and variance of the spine
of a random 2-SAT formula, and are at the heart of our
proofs.
Before we can proceed,
we unfortunately need a short interlude
on Landau symbols:

In this paper, we shall use Landau's notation
$f=O(g)$ and $f=o(g)$.
As usual, $f=O(g)$ stands for a bound
$|f|\leq c|g|$, where $c$ is a universal constant,
unless otherwise specified.
If we have
a bound of the form $|f|\leq h(g)|g|$,
where $h(g)$ is a function which is bounded above,
though not necessarily uniformly, for
finite $g$, and which is uniformly bounded above
as $g$
goes to zero, we shall use the notation
\renewcommand\O{{\mathbb O}}
$f=\O_0(g)$.
In this notation,
$e^{x^2}-1$ is $\O_0(x^2)$,
but it is not $O(x^2)$.

Our use of the symbol $o(g)$ is slightly stronger
than usual.  Typically, $f=o(g)$ means that
$f/g$ goes to zero as the independent variables
in question tend to their limiting values, but
usually $f=o(g)$ does not require that
$f/g$ is bounded in the whole domain of the
independent variables.  We require both
uniform boundedness and that $f/g$ tends
to zero.
Since it may be ambiguous which independent
variables tend to $\infty$ or $0$ in an expression of
the form $f=o(g)$, we frequently specify the variables
in question.  Thus $f=o_{\lambda,\varepsilon}(g)$ means
that $f/g \rightarrow 0$ as $\lambda \rightarrow \infty$
and $\varepsilon \rightarrow 0$.
For example, in this notation, the $o(1)$ terms
in Theorem~\ref{thm1.3} would be written as
$\o(1)$.

Finally, as mentioned earlier, $f=\Theta(g)$ means
that there exist positive, finite
constants $c_1$ and $c_2$ such that $c_1 \leq f/g \leq c_2$.
Unless noted otherwise, these constants are universal.

\addtolength{\textheight}{1\baselineskip}

\begin{theorem}
\label{spine-exp} There are constants $\lambda_0$ and $\varepsilon_0$,
$0<\lambda_0<\infty$ and $0<\varepsilon_0<1$, such that the
following statements hold for
\begin{equation}
\label{p-epsilon}
p=\frac 1{2n}(1 + \varepsilon)
=\frac 1{2n}(1 + \lambda_n n^{-1/3})
\end{equation}
and $\lambda_0\leq |\lambda_n|\leq\varepsilon_0 n^{1/3}$.

\noindent
i) If $\varepsilon<0$, then
\begin{equation}
\label{left-spine-exp}
\Pr\Big(
x\underset {F_{n,p}}\rightsquigarrow\overline x
\Big)
=  \frac{n^{-1/3}}{4\lambda_n^2}
\big(1 +  \o(1)\big)~.
\end{equation}

\noindent
ii) If $\varepsilon>0$, then
\begin{align}
\label{right-spine-exp}
\Pr\Big(
x\underset {F_{n,p}}\rightsquigarrow\overline x
\Big)
&= \vartheta(\varepsilon)
\big(1+o_{\lambda_n}(1\big)),
\notag
\\
&=2\lambda_n n^{-1/3}
\big(
1+o_{\lambda_n}(1\big)+O(\varepsilon)
\big)
\end{align}
where
$\vartheta(\varepsilon)$ is defined in
\eqref{def-theta}.

\end{theorem}

\begin{theorem}
\label{spine-var} Let
$p$,
$\varepsilon$ and $\lambda_n$ be as in
Theorem~\ref{spine-exp}.  Then the following statements
hold for all
strictly distinct literals $x$ and $y$.

\noindent
i) If $\varepsilon<0$, then
\begin{equation}
\label{left-spine-var}
\Pr
\Big(
x\underset {F_{n,p}}\rightsquigarrow\overline x
\Big)
\Pr
\Big(
y\underset {F_{n,p}}\rightsquigarrow\overline y
\Big)
\leq
\Pr
\Big(
x\underset {F_{n,p}}\rightsquigarrow\overline x
\text{ and }
y\underset {F_{n,p}}\rightsquigarrow\overline y
\Big)
 = O\Big(\frac {n^{-2/3}}{\lambda_n^4}\Big).
\end{equation}

\noindent
ii) If $\varepsilon>0$, then
\begin{equation}
\label{right-spine-var}
0 \leq
\Pr
\Big(
x\underset {F_{n,p}}\rightsquigarrow\overline x
\text{ and }
y\underset {F_{n,p}}\rightsquigarrow\overline y
\Big)
-
\Pr
\Big(
x\underset {F_{n,p}}\rightsquigarrow\overline x
\Big)
\Pr
\Big(
y\underset {F_{n,p}}\rightsquigarrow\overline y
\Big)
  = O\Big(\frac {n^{-2/3}}{\lambda_n}\Big).
\end{equation}
\end{theorem}

\begin{remark}
\label{spine-var-ext}
By monotonicity, the bound
\eqref{left-spine-var}
can be extended to
all $\lambda_n\in [-n^{1/3},-\lambda_0]$.
Indeed, using that the events $x\rightsquigarrow \overline{x}$
and $y\rightsquigarrow \overline{y}$ are monotone events,
we have that
\begin{equation*}
\Pr
\Big(
x\underset {F_{n,p}}\rightsquigarrow\overline x
\text{ and }
y\underset {F_{n,p}}\rightsquigarrow\overline y
\Big)
\leq
\Pr
\Big(
x\underset {F_{n,p_0}}\rightsquigarrow\overline x
\text{ and }
y\underset {F_{n,p_0}}\rightsquigarrow\overline y
\Big)
\end{equation*}
provided $p\leq p_0$.
Setting $p_0=(1-\varepsilon_0)/2n$,
using equation \eqref{left-spine-var} to bound
the right hand side by
$O\big( {n^{-2/3}}/{(n^{1/3})^4}\big)$,
and observing that
${n^{-2/3}}/{(n^{1/3})^4}=O\big( {n^{-2/3}}/{\lambda_n^4}\big)$
provided $\lambda_n\in [-n^{1/3},-\varepsilon_0 n^{1/3}]$,
we obtain that
\begin{equation}
\label{2-point}
\Pr
\Big(
x\underset {F_{n,p}}\rightsquigarrow\overline x
\text{ and }
y\underset {F_{n,p}}\rightsquigarrow\overline y
\Big)
= O\Big(\frac {n^{-2/3}}{\lambda_n^4}\Big)
\quad \text{for all} \,\,
\lambda_n\in [-n^{1/3},-\lambda_0].
\end{equation}

\end{remark}

Given the above two theorems, we shall
prove the following analogue of
Theorem~\ref{sat_{n,m}} for the ensemble $F_{n,p}$.

\begin{theorem}
\label{sat_{n,p}}
There are constants $\lambda_0$ and $\varepsilon_0$,
$0<\lambda_0<\infty$ and $0<\varepsilon_0<1$, such that the
following statements hold for
$
p=\frac 1{2n}(1 + \lambda_n n^{-1/3})
$
and
$\lambda_0\leq |\lambda_n|\leq\varepsilon_0 n^{1/3}$.

\noindent
i) If $\lambda_n<0$, then
\begin{equation}
\label{1.3}
\Pr(F_{n,p} \; \text{\rm is SAT})
=\exp\big(-\Theta\big({|\lambda_n|^{-3}}\big)\big).
\end{equation}

\noindent
ii) If $\lambda_n>0$, then
\begin{equation}
\label{1.4}
\Pr(F_{n,p} \; \text{\rm is SAT})
=\exp\big(-\Theta\big({\lambda_n^3}\big)\big).
\end{equation}
\end{theorem}

For fixed $\lambda_n$,
both (\ref{left-spine-exp}) and
(\ref{right-spine-exp}) are of the form
$\Pr(x\rightsquigarrow\overline x)
=\Theta(n^{-1/3})$.
Together with equation \eqref{2.4},
Theorem~\ref{spine-exp} therefore implies that
the expected size of the spine
scales like $n^{2/3}$, provided
$\lambda_n$ stays bounded as $n\to\infty$.
The heuristics for this scaling with
$n$ will be given in the next subsection, and the actual proof
of the scaling will be given in Sections
\ref{sec:prelim}--\ref{sec:spine}.
Theorem~\ref{spine-var} allows us
to control the deviations of
the random variable $|S(F_{n,p})|$
from its expectation; its proof will be
given in Section \ref{sec:var}.
Together, these two
theorems allow us to prove Theorem~\ref{sat_{n,p}},
which is just the analogue of Theorem~\ref{sat_{n,m}}
in the model $F_{n,p}$.  In the final subsection, we shall describe the
strategy of this proof.  While the
actual proof is easier in the model
$F_{n,p}$, the heuristic argument
is easier in the model
$F_{n,m}$.  Our goal in the last subsection
is therefore
to describe how the scaling
$n^{2/3}$ for the size of the
spine in the model $F_{n,m}$
leads to bounds of the form
\eqref{1.3} and \eqref{1.4}.
The actual proof of Theorem~\ref{sat_{n,p}}
is given in Section \ref{sec:sat}.

\subsection{\sc Heuristics for the Scaling of the Spine}

The proof of Theorem~\ref{spine-exp} (and hence also the
proof of the expected size of the spine, Theorem~\ref{thm1.3})
uses the digraph representation of the last section.  Indeed, by
Lemmas
\ref{lem:spine-digraph},
\ref{lem:spine-s.d.},
\ref{trimmed-graph} (iii) and \ref{trimmed-graph} (iv),
and the fact the probability of the event
$x\rightsquigarrow\overline x$
does not depend on the choice of the literal $x\in [n]$,
we have
\begin{align}
E(|S(F_{n,p})|)
&=
2n\Pr
\left(
x\rightsquigarrow\overline x
\right)
\notag
\\
&=
2n\Pr
\left(
   L^+_{F_{n,p}}(x)
   \neq
   \widetilde L^+_{F_{n,p}}(x)
\right)
\notag
\\
&=
2n\sum_{k=1}^n
\left[
\Pr\left(|\widetilde L^+_{F_{n,p}}(x)|=k\right)
-
\Pr
\left(
   |L^+_{F_{n,p}}(x)|=k,\;  L^+_{F_{n,p}}(x)
   \text{ is s.d}
  \right)
\right]
\notag
\\
&=
2n\sum_{k=1}^n
\left[
\Pr\left(|C_{{n,{2p-p^2}}}(x)|=k\right)
-
\Pr
\left(
   |L^+_{F_{n,p}}(x)|=k,\;  L^+_{F_{n,p}}(x)
   \text{ is s.d}
  \right)
\right].
\label{h.1}
\end{align}
It turns out that for
${2p-p^2}$ near to the random graph threshold
$1/n$,
and
$k\leq \Theta(n^{2/3})$,
the size of the largest component in the random
graph,
the probability that
$L^+_{F_{n,p}}(x)$
is strictly distinct and has size $k$
is well approximated by
$\Pr\left(|C_{{n,{2p-p^2}}}(x)|=k\right)$,
so that the summand in equation (\ref{h.1}) is
approximately zero.  On the other hand,
for ${2p-p^2}$ near $1/n$ and
$k\geq \Theta(n^{2/3})$, only the
sum over
$\Pr\left(|C_{{n,{2p-p^2}}}(x)|=k\right)$
contributes to \eqref{h.1}.
Thus we can approximate
\begin{equation}
\label{h.2}
E(|S(F_{n,p})|)
\approx
2n\Pr
\left(|C_{n,{2p-p^2}}(x)|\geq n^{2/3}\right).
\end{equation}
As the reader might imagine, the above arguments
require a good deal of justification; see Section
\ref{sec:prelim}--\ref{sec:spine}
for  precise bounds.
But for ${2p-p^2}$ near $1/n$, the
probability that $|C_{n,{2p-p^2}}(x)|\geq n^{2/3}$
scales like the probability that
$x$ lies in the largest component
in the random graph, which in turn scales like
$n^{-1/3}$ (see e.g.~ \cite{Bol85}).  This implies that the
expected size of the spine $S(F_{n,p})$
scales like $n^{2/3}$ provided
$p$ is of the form $p=\frac 1{2n}(1\pm\Theta(n^{-1/3}))$.
Observing that the models $F_{n,p}$
and $F_{n,m}$ are equivalent as long as $m$
is near its expected value $\binom{n}{2}p$ (see Appendix A),
we obtain
the scaling of Theorem~\ref{thm1.3}.

\subsection{\sc Heuristics for the Scaling of the Window}

As explained earlier, our goal is to describe how the behavior
$E(|S(F_{n,m})|) = \Theta (n^{2/3})$ leads to the bounds
\eqref{1.3} and \eqref{1.4}.
To this end, consider a process
which builds random formulas
as follows: Given a 2-SAT formula $F_m$, let $F_{m+1}
 =F_m\wedge C$, where  $C=x\vee y$ is chosen
uniformly at random from
the set of all 2-clauses  over $\{x_1, \dots, x_n\}$
that have not yet
been used in $F_m$.  Obviously, the distribution of $F_m$ is
the same as that of $F_{n,m}$.  Furthermore, $F_{m+1}$
is satisfiable if and only if $F_{m}$ is satisfiable
and either $x$ or $y$ does not lie in the spine of $F_m$.
Conditioned on the events that $F_m$ is SAT and that
$S(F_m)$ has size $s$, the probability that $F_{m+1}$
is SAT is therefore equal to
\begin{align}
\Pr\Big(\SAT(
F_{m+1})
\Big|
\SAT(F_{m})\text{ and }|S(F_m)|=s\Big)
&=
1-{\binom{s}{2}}\left({4\binom{n}{2}-m}\right)^{-1}
\notag
\\
&=1-\frac{s(s-1)}{4n(n-1)-2m}.
\label{heur.1}
\end{align}
By the analogue of Theorem~\ref{spine-exp}
for the model $F_{n,m}$
and the monotonicity of the expected size of
the spine, we have that
near the transition point,
$E(|S(F)|)=\Theta(n^{2/3})$.  Neglecting the
difference between statements in expectation
and statements in probability, equation (\ref{heur.1})
therefore implies that the probability that
$F_{m+1}$ is UNSAT, conditioned on $F_m$ being
satisfiable, is $\Theta(n^{4/3}/n^2)=\Theta(n^{-2/3})$.
After $\Theta(n^{2/3})$ steps, a satisfiable formula
therefore becomes UNSAT, giving the a finite-size scaling
window of width $\Theta(n^{2/3})$ in $m$, and hence of
width $\Theta(n^{-1/3})$ in $m/n$.  This is the result
in Corollary \ref{cor:window}.

In order to explain heuristically the error bounds
$\exp(-\Theta(|\lambda_n|^{-3}))$
and $\exp(-\Theta(\lambda_n^{3}))$ in Theorem~\ref{sat_{n,m}},
we proceed as follows.
As a consequence of (\ref{heur.1}),
we have that
\begin{align}
\Pr\Big(\SAT(
F_{m+1})
\mid
\SAT(F_{m})\Big)
=1-
\frac{{E}\Big((|S(F_{m})|^2-|S(F_{m})|)\Big|
\SAT(F_{m})\Big)}
{4n(n-1)-2m},
\label{heur.2}
\end{align}
and hence
\begin{equation}
\Pr(\SAT(F_m))
=
\prod_{k=0}^{m-1}
\bigg(
1-
\frac{{E}\Big((|S(F_{k})|^2-|S(F_{k})|)\Big|
\SAT(F_{k})\Big)}
{4n(n-1)-2k}
\bigg).
\label{heur.2a}
\end{equation}
Neglecting the difference between conditional and
unconditional expectations,
and approximating $|S(F_{k})|^2-|S(F_{k})|$
by  $|S(F_{k})|^2$,
we get
\begin{align}
\Pr(\SAT(F_m))
&\approx
\prod_{k=0}^{m-1}
\bigg(
1-
\frac{1}{4n^2-2k}
{E}\big(|S(F_{k})|^2\big)
\bigg)
\notag
\\
&\approx
\exp
\bigg(-\sum_{k=0}^{m-1}
\frac{1}{4n^2}
{E}\big(|S(F_{k})|^2\big)
\bigg)
.
\label{heur.3}
\end{align}

For $\lambda_n<0$ and
$k=(1+\lambda n^{-1/3})n$, $\lambda\in(-n^{1/3},\lambda_n)$,
we then approximate ${E}\big(|S(F_{k})|^2\big)$
by ${E}\big(|S(F_{n,p})|^2\big)$,
$p=(1+\lambda n^{-1/3})/2n$.  Using Theorems \ref{spine-exp}
and \ref{spine-var}
to estimate this probability, we have
\begin{align}
\Pr(\SAT(F_m))
&\approx
\exp
\bigg(-\frac{1}{4n^2}\int_{k=0}^{m}dk
{E}\big(|S(F_{k})|\big)^2
\bigg)
\notag
\\
&\approx
\exp
\bigg(-\frac{1}{4n^2}\int_{-n^{1/3}}^{\lambda_n}
n^{2/3}d\lambda
\,\Theta\Big(n^{2/3}/\lambda^2\Big)^2
\bigg)
\notag
\\
&=\exp\bigl(\Theta(|\lambda_n|^{-3})\bigr),
\label{heur.3a}
\end{align}
giving the bound (\ref{new1.1}) below the threshold.  In a similar way,
we can integrate the bound (\ref{right-spine-exp}) to obtain a
heuristic derivation of (\ref{new1.1}) above the threshold.

The actual proof of Theorem~\ref{sat_{n,m}}
in the form (\ref{1.3})--(\ref{1.4}) will
be given in the next section, and relies heavily on
Theorems \ref{spine-exp} and \ref{spine-var}, which
in turn are proven in Sections
\ref{sec:prelim}--\ref{sec:var}.  In addition, we shall
need two more technical lemmas, to be proven in
Section \ref{sec:hourglass}.
In  Appendix A, we discuss the
relation between the models $F_{n,m}$ and
$F_{n,p}$.

\bigskip

\section{Probability of Satisfiability} \label{sec:sat}

In this section we prove Theorem~\ref{sat_{n,p}}, which together
with Appendix A establishes Theorem~\ref{sat_{n,m}}.  The lower bounds
depend on the second moment estimates of Theorem~\ref{spine-var},
which is proved in Section~\ref{sec:var}.  The upper bounds depends
on a theorem
showing that with high probability there are many structures, to be
called {\it hourglasses}, that are ``seeds'' for the growth of the
spine.
The hourglass theorem is proven in
Section~\ref{sec:hourglass}, after we develop suitable machinery in
the intervening sections.

To derive the bounds, we shall find it convenient to
view the random formula $F_{n,p}$ as a process, and to
couple the processes
for all possible values of $n$
and $p$.
To do this, for each unordered pair of natural numbers we
pick four random variables uniformly distributed on the interval from
0 to 1, so that the set of all these random variables indexed by 4
copies of $\binom{\mathbb N}{2}$ (sometimes denoted ${\mathbb
N}^{(2)}$) is a set of independent random variables.  As usual
$\binom{X}{2}$ denotes the set of unordered pairs of elements of the
set $X$.  A pair of natural numbers specify two different variables of
the formula, and the four random variables associated with the pair
correspond to the four different clauses that can be made using these
variables, so that each possible clause $C$ has its own independent random
variable $U_C$ distributed uniformly in $[0,1]$.
A clause $C$ appears in the formula on $n$ variables
at probability $p$ precisely when the indices of its two variables are
not larger
than $n$, and $U_C<p$.  We shall call $U_C$ the {\it birthday\/}
of the clause $C$.  The process $\widetilde{F}$ is the collection of
all
these random variables, and it defines a family of formulas by
\begin{equation}
\label{process}
F_{n,p} = \bigwedge_{C:n(C)\leq n,\ U_C<p} C,
\end{equation}
where $n(C)$ denotes the maximum index of the two variables in
clause $C$.
It is easy to see that for each value of $n$ and
$p$, the distribution of the resulting formula $F_{n,p}$ is exactly the
distribution introduced before, which justifies our using the same
notation as before.  By construction, satisfiability of $F_{n,p}$ is
monotone decreasing in $n$ and $p$.

\subsection*{\sc Lower bounds}

To derive the lower bounds of Theorem~\ref{sat_{n,m}},
it is sufficient to consider the coupling
above for a fixed value of $n$, so we suppress the variable $n$.
Given $\widetilde{F}$ (i.e., $F_p$ for each $p$), we define
the reduced formula process
$\widetilde{\Phi} = (\Phi_p)_{p\in[0,1]}$ as follows: $\Phi_{0}$
has no
clauses.  Start with $p=0$, and increase it until $p=1$.
Clauses are added to $F_p$ one at a time; each time we add a clause to
$F_p$, we also add that clause to $\Phi_p$ provided that doing so does
not make $\Phi_p$ unsatisfiable.

Given $\widetilde{\Phi}$ we can define a new formula process
$\widetilde{H}$
as follows: For each clause $C$ appearing in $\widetilde{\Phi}$,
set its
birthday in $\widetilde{H}$ to coincide with its birthday in
$\widetilde{\Phi}$.  For each clause $C$ not appearing in $\widetilde\Phi_{1}$,
there is some smallest value $p_{\min}$ of $p$ for which $\Phi_p\wedge
C$ is not satisfiable.  Pick the birthday of $C$ in $H_p$ uniformly at
random in the interval $[p_{\min},1)$.  Since the $\widetilde{H}$ process
is drawn uniformly at random from the set of $\widetilde{F}$ processes with
reduced formula process $\widetilde{\Phi}$, the $\widetilde{H}$
and $\widetilde{F}$
processes are identically distributed.  Also note that the first time
that $H_p$ differs from $\Phi_p$ is also the first time that $H_p$
becomes unsatisfiable.

Let $U(F)$ denote the number of different clauses $C$ such that
$F\wedge C$ is unsatisfiable.  If $F$ is itself satisfiable, then
$$
U(F) = \binom{|S(F)|}{2},
$$
where as usual $S(F)$ denotes the spine of $F$.

Conditional on the reduced formula process being $\widetilde\Phi$ and $H_p$
being satisfiable (i.e.\ that $H_p=\Phi_p)$, the probability that
$H_{p+\delta}$ is satisfiable is
$$
\left(1-\frac{\delta}{1-p}\right)^{U(\Phi_p)}
=\,\exp\left[-\frac{U(\Phi_p)}{1-p}\delta+O(\delta^2)\right]
$$
provided $\delta$ is small enough.
Multiplying these
probabilities for $p=0,\delta,2\delta,\dots$,
and passing to
the limit $\delta\to 0$,
we find that conditional upon
$\widetilde{\Phi}$, the probability that $H_p$ is satisfiable
is given by
$$ \Pr[\SAT(H_p)|\widetilde{\Phi}]
= \exp\left[-\int_0^p \frac{U(\Phi_s)}{1-s} ds\right].$$
We have

\begin{align*}
\Pr[\SAT(F_p)] &= \Pr[\SAT(H_p)] \\
 &= E_{\widetilde{\Phi}} \big[\Pr[\SAT(H_p)|\widetilde{\Phi}]\big] \\
 &= E_{\widetilde{\Phi}} \left[\exp\left[-\int_0^p
\frac{U(\Phi_s)}{1-s} ds\right]\right] \\
 &\ge \exp E_{\widetilde{\Phi}}\left[-\int_0^p
\frac{U(\Phi_s)}{1-s} ds\right] \\
 &= \exp\left[-\int_0^p \frac{E_{\widetilde{\Phi}}
[U(\Phi_s)]}{1-s} ds\right] \\
 &= \exp\left[-\int_0^p \frac{4\binom{n}{2}
\Pr[\text{$x\rightsquigarrow \overline{x}$ and
$y\rightsquigarrow \overline{y}$ in $\Phi_s$}]}{1-s} ds\right] \\
 &\ge \exp\left[-4\binom{n}{2}\int_0^p
\frac{\Pr[\text{$x\rightsquigarrow \overline{x}$ and
$y\rightsquigarrow \overline{y}$ in $F_s$}]}{1-s} ds\right],
\end{align*}
where $x$ and $y$ are fixed strictly distinct literals.

Next we proceed to estimate the integral.  Since we are principally
interested in the case $p=O(1/n)$, let us assume $p=o_n(1)$ so that
the effect of the denominator of
the integrand is negligible.
Recalling that $p$ and $\lambda_n$ are related by \eqref{p-epsilon},
and
setting
$$
s =s(t)= \frac{1+ t n^{-1/3}}{2n},
$$
we get
\begin{equation}
 \Pr[\SAT(F_p)]
\ge \exp\left[-(1+o_n(1)) n^{2/3}
\int_{-n^{1/3}}^{{\lambda_n}}
\Pr[x\rightsquigarrow \overline{x}
    \text{ and }
   y\rightsquigarrow \overline{y}
   \text{ in }
   F_{s(t)}]
dt\right].
\label{4:1}
\end{equation}

\begin{remark}
It is not hard to derive the analogue of \eqref{4:1}
in the model $F_{n,m}$.  Indeed, starting from
\eqref{heur.2a}, rewriting
$$
E\Big(|S(F_{k})|^2 - |S(F_{k})|\Big| \SAT(F_k)\Big)
=
4\binom n2
\Pr\big(x\rightsquigarrow\overline x
\text{ and }y\rightsquigarrow\overline y
\text{ in } F_k
\big| \SAT(F_k)\big),
$$
where $x$ and $y$ are strictly distinct,
and observing that by the FKG inequality,
$$
\Pr\big(x\rightsquigarrow\overline x
\text{ and }y\rightsquigarrow\overline y
\text{ in } F_k
\big| \SAT(F_k)\big)
\leq
\Pr(x\rightsquigarrow\overline x
\text{ and }y\rightsquigarrow\overline y
\text{ in } F_k),
$$
one gets
$$
\Pr\big(\SAT(F_{n,m})\big)
\geq
\prod_{k=0}^{m-1}
\Big[1-
\big(1+O(m/n^2)\big)
\Pr\big(x\rightsquigarrow\overline x
\text{ and }y\rightsquigarrow\overline y
\text{ in } F_{n,k}\big)
\Big].
$$
\end{remark}

\begin{proof_left}{Proof of the lower bound of Theorem~\ref{sat_{n,p}}}
For $t\in [-n^{1/3}, -\lambda_0]$, the probability in the
integrand in equation (\ref{4:1}) is
$O(n^{-2/3}/t^4)$ by
Theorem~\ref{spine-var}
and Remark \ref{spine-var-ext}.
Integrating, we find
that
\begin{align*}
\Pr[\SAT(F_p)]
\ge \exp\left[\Theta\left(
\frac{1}{{\lambda_n}^3}\right)\right]
= 1
- \Theta\left(
\frac{1}{{|\lambda_n|}^3}\right)
\end{align*}
provided
${\lambda_n}\in [-n^{-1/3},-\lambda_0]$.
\end{proof_left}

\begin{proof_right}{Proof of the lower bound of Theorem~\ref{sat_{n,p}} }
By Theorems~\ref{spine-exp} and \ref{spine-var},
the probability in
the integrand in equation (\ref{4:1}) is
$4n^{-2/3}t^2(1+O(\varepsilon)+o_{\lambda_n}(1))$,
provided $t\in [\lambda_0,\varepsilon_0 n^{1/3}]$.
In the middle region $t\in [-\lambda_0,\lambda_0]$
we upper bound the probability in the integrand
by $\Theta(n^{-2/3})$, and we bound it in the left
region $t\in[-n^{1/3},-\lambda_0]$ by $O(n^{-2/3}/t^4)$ as above.
Integrating,
we find that
\begin{align*}
\Pr[\SAT(F_{p})]
&\ge \exp\left[-\frac{4+o_{\varepsilon,{\lambda_n}}(1) }{ 3}
{\lambda_n}^3\right].\end{align*}
provided $\lambda_n\in[\lambda_0,\varepsilon_0 n^{1/3}]$.
\end{proof_right}

\subsection*{\sc Upper bounds}

To bound from above the probability of satisfiability, we
use
Theorem~\ref{thm:hourglasses} below, which states that with high
probability
there are certain types of structures contained within the directed
graph associated with a formula.

\begin{definition}
\label{hourglass}
An {\em hourglass} is a triple $(v,I_v,O_v)$
where $v$ is a literal, and
$I_v$ and $O_v$
are two disjoint sets of literals not
containing $v$,
such that for each $x\in I_v$,
there is a path $x\rightsquigarrow v$ in ${I_v}\cup\{v\}$,
and for each $y\in O_v$,
there is a path $v\rightsquigarrow y$ in ${O_v\cup\{v\}}$.
Furthermore, we require that $\{v\}\cup I_v\cup O_v$
is strictly distinct.  We call
$v$ the {\em central literal},
$I_v$ the  {\em in-portion},
and $O_v$ the {\em out-portion}
of the hourglass.
\end{definition}

\begin{theorem}
\label{thm:hourglasses}
There are constants $\lambda_0$ and $\varepsilon_0$,
$0<\lambda_0<\infty$ and $0<\varepsilon_0<1$,
such that for
$
p=\frac 1{2n}(1 + \lambda_n n^{-1/3})
$
and
$\lambda_0\leq |\lambda_n|\leq\varepsilon_0 n^{1/3}$,
the following statements hold with probability
at least $1-\exp(-\Theta(|\lambda_n|^3)$.

\noindent i)
If $\lambda_n<0$, then there are at least
$\Theta(|\lambda_n|^3)$ disjoint,
mutually strictly distinct hourglasses with in-portion
and out-portion each of size at least
$n^{2/3}/\lambda_n^2$.

\noindent ii)
If
$\lambda_n>0$, then there is at least one
hourglass with in-portion
and out-portion each of size
$\Theta(\lambda_n n^{2/3})$.
\end{theorem}

The proof of this theorem will be given in Section \ref{sec:hourglass}.
There we shall use the coupling of the trimmed out-graph
$\widetilde G^+_{F_{n,p}}(x)$ to the random graph process
$G_{n,2p-p^2}$ (see Lemma~\ref{trimmed-graph} and its proof) to
explicitly construct $\Theta(|\lambda_n|^3)$ many hourglasses
below threshold.  To prove the theorem above threshold, we shall
show that, when $\lambda$ is increased from its value below the
threshold to its value above the threshold, a constant fraction
of these subcritical hourglasses
will merge into one giant hourglass
of size $\Theta(|\lambda_n|^3)\Theta(n^{2/3}/\lambda_n^2)$.
See Section  \ref{sec:hourglass} for details.

Here, we shall use the hourglasses to derive the upper bounds on
satisfiability both to the left and to the right of the window.

\begin{proof_left}{Proof of upper bound of Theorem~\ref{sat_{n,p}}}
To get the bound on the left, we increase $p$
from
$(1-tn^{-1/3})/2n$
to
$(1-(t/2)n^{-1/3})/2n$,
with $\lambda_0\leq t\leq \varepsilon_0 n^{1/3}$.
For any pair of vertices, whether or
not there was a directed edge between them before, afterwards the
probability of finding such an edge is at least $(t/4)n^{-4/3}$.  For
each hourglass, for each pair of literals $u$ and $v$ in the
out-portion of the hourglass, if the clause $\overline{u}\vee\overline{v}$ appears,
then we claim that the central vertex and the entire in-portion of the
hourglass is afterwards part of the spine of the formula.
Indeed, let $x$ be such a vertex.  Then $x\rightsquigarrow u$
and $x\rightsquigarrow v$ since $u$ and $v$ are in the
out-portion of the hourglass.  But the appearance of the clause
$\overline{u}\vee\overline{v}$ implies that $u\rightarrow \overline v$,
so that we have $x\rightsquigarrow u\rightarrow \overline v$.
Together with $x\rightsquigarrow v$, which is equivalent
to $\overline v \rightsquigarrow \overline x$, we conclude that
$x \rightsquigarrow \overline x$.
 The
probability of the event
that
the clause $\overline{u}\vee\overline{v}$ appears
is at least $\Theta((n^{2/3}/t^2)^2 t
n^{-4/3}) = \Theta(1/t^3)$.  If furthermore a clause appears
that contains two
literals in the in-portion, then the formula is not
satisfiable.  These events are
independent,
so the probability that
they both occur is at least $\Theta(1/t^6)$.  But since with high
probability there are $\Theta(t^3)$ hourglasses, with probability at
least $\Theta(1/t^3)$
the formula becomes unsatisfiable.  Setting $t=2|\lambda_n|$ gives the desired
upper bound on the left.
\end{proof_left}

\begin{proof_right}{Proof of upper bound of Theorem~\ref{sat_{n,p}}}
To get the bound on the right, we start with $p = (1+tn^{-1/3})/2n$
(with $\lambda_0\leq t\leq \varepsilon_0 n^{1/3}$), where the probability
that there is no giant hourglass is at most $\exp(-\Theta(t^3))$, and
then increase it to $(1+2tn^{-1/3})/2n$.
Any clauses of the form
$(\overline{u}\vee\overline{v})$ where $u$ and $v$ are in the
out-portion of the the giant hourglass
beforehand will afterwards appear with probability
at least $tn^{-4/3}$.
This will cause the in-portion of the giant
hourglass to become part of the spine of the formula, except with
probability that can be bounded by
$\left(1-tn^{-4/3}\right)^{\Theta(tn^{2/3})^2}=
\exp(-\Theta(t^3))$.
Furthermore, any
clauses of the form $(u\vee v)$ where $u$ and $v$ are in the
in-portion of the giant hourglass beforehand
will afterwards appear with probability
at least $tn^{-4/3}$.
Therefore the formula will become
unsatisfiable, except with probability that is again exponentially small in
$t^3$.  Setting $t=\lambda_n/2$ completes the proof.
\end{proof_right}

\begin{remark}
Instead of using Theorem~\ref{thm:hourglasses},
one can alternatively use Theorems~\ref{spine-exp}
and \ref{spine-var} to prove that below the window,
the spine has size at least $E(|S(F_{n,p})|)/2$
with probability
$\exp(-O(1))$.
Increasing $p$ as in the above proof,
one obtains an alternative
proof of the fact that, with  probability at
least $\Theta(1/|\lambda_n|^3)$, the formula becomes unsatisfiable
below the threshold.   While a similar argument can
be used to show that above the window, the probability
of satisfiability goes to zero, we cannot use
Theorems \ref{spine-exp}
and \ref{spine-var} alone to prove that it goes to zero exponentially
fast in $\lambda_n^3$.  For this, we need the hourglass theorem.

\end{remark}
\newpage

\section{Machinery from Random Graph Theory} \label{sec:prelim}
\addtolength{\textheight}{-1\baselineskip}

In this section we establish several bounds needed in the
proofs of Theorems \ref{spine-exp} and \ref{spine-var}.
As in earlier sections,
we use the notation $G_{n,p}$ for the (unoriented)
random graph on $\{1,2,\dots,n\}$ with edge
probability $p$.  We also consider $D_{n,p}$,
the random directed graph on
$\{1,2,\dots,n\}$ in which each oriented edge
is chosen independently with probability $p$.
We use the symbol $L_{n,p}^{^+}(x)$
to denote both the set of vertices
$y\in\{1,2,\dots,n\}$
that can be reached from a vertex $x$ in the random
digraph $D_{n,p}$, and the set of vertices
$y\in [n]$ that can be reached from
a vertex $x$
in the digraph corresponding to a random
2-SAT formula $F_{n,p}$.  If the difference
is not clear from the context, we shall
use the notations
$L_{D_{n,p}}^{^+}(x)$
and
$L_{F_{n,p}}^{^+}(x)$
to distinguish the two cases.
We begin this section with a basic lemma
which is implicit in the work of Karp
\cite{Kar90}.

\begin{lemma}\label{basic1}
The probability that in the random digraph
$D_{n,p}$ every vertex can be reached
from a given vertex is precisely the
probability that the random graph $G_{n,p}$
is connected.
\end{lemma}
\begin{proof}
We may assume that the vertex in question is
vertex $1$.
First, we shall inductively define a subtree
$T=T(D_{n,p})$ of $D_{n,p}$, rooted
at $1$, with each edge oriented away from $1$.
To this end, set $X_0=\{1\}$, $Y_0=\emptyset$,
and let $T_0$ be the subtree of
$D_{n,p}$ with the single vertex $1$.  Suppose
that we have defined a pair $(X_i, Y_i)$
of subsets of $\{1,2,\dots,n\}$ with
$Y_i \subset X_i$, and a subtree $T_i$
with $V(T_i)=X_i$. (We think of $Y_i$ as
the set of vertices we have ``exposed'',
i.e.~tested for outgoing edges, and $X_i$
as the set of vertices
we have selected so far.)
If $X_i=Y_i$ then $T_i$ is our tree $T$.
Otherwise, let $x_i$
be the smallest element of $X_i \setminus Y_i$.
Let $\Gamma^+(x_i)$ denote the set of vertices
in $\{1,2,\dots,n\}$
that can be reached by single
edges of $D_{n,p}$ that are oriented outward from $x_i$.  Now
set $X_{i+1}=
X_i\cup \Gamma^+(x_i)$, $Y_{i+1}=Y_i\cup \{x_i\}$,
and take $T_{i+1}$ to be obtained from $T_i$ by
adding to it the vertices in $X_{i+1}\setminus X_i$,
together with all the edges from $x_i$ to
$X_{i+1}\setminus X_i$.  The vertex set of the
subtree $T$ of $D_{n,p}$ constructed in this
way is clearly $L_{n,p}^{^+}(1)$; in particular,
$L_{n,p}^{^+}(1)=\{1,2,\dots,n\}$ iff
$V(T)=\{1,2,\dots,n\}$.  Since the edges
of $T$ are oriented away from $1$, we
may view $T$ as an {\em unoriented\/} tree.

Now, let us construct a subtree $T'=T'(G_{n,p})$
of the random graph $G_{n,p}$ rooted at $1$
by precisely the same algorithm. The lemma
will follow if we show that
\begin{equation}
\Pr(\, T(D_{n,p})=T_0\, )
=\Pr(\, T'(G_{n,p})=T_0\, )
\label{prob1}
\end{equation}
for every tree $T_0$ with vertex
set $\{1,2,\dots,n\}$.

Given $T_0$, we can define the $X_i$'s as above.
We have $T(D_{n,p})=T_0$ if and only
if the random digraph $D_{n,p}$ is such that

1) it contains all the edges of $T_0$ (oriented away from vertex $1$),

2) it contains no edge oriented from $x_i$ to
$\{1,2,\dots,n\}\setminus X_{i+1}$.

\noindent
Similarly, $T'(G_{n,p})=T_0$ if and only if the
random graph $G_{n,p}$ is such that

1) it contains all the edges of $T_0$,

2) it contains no edge from
$x_i$ to $\{1,2,\dots,n\}\setminus X_{i+1}$.

\noindent
Notice that the probability that $D_{n,p}$ contains
a given set $K$ of oriented edges,
and no edge of a second set $K'$ of oriented
edges, is $p^{|K|}(1-p)^{|K'|}$ provided that
$K\cap K'=\emptyset$.  Moreover, this is
equal to the
probability that $G_{n,p}$ contains a set
$\widetilde K$ of unoriented edges, and no edge
of a second set $\widetilde K'$ of unoriented
edges, provided that
$\widetilde K\cap \widetilde K'=\emptyset$,
$|\widetilde K|=|K|$ and $|\widetilde K'|=|K'|$.
Thus
relation (\ref{prob1}) holds, and we are done.
\end{proof}

Returning to the 2-SAT problem $F_{n,p}$,
let $x$ be a fixed literal. The probability that, in a random
2-SAT formula $F_{n,p}$, the set $L_{n,p}^{^+}(x)$ consists of
$k$ strictly distinct literals is trivially independent of
$x$; we shall denote it by $P_{n,p}(k)$:
\begin{equation}
P_{n,p}(k)
=\Pr(\{|L_{n,p}^{^+}(x)|=k\}
   \cap
   \{L_{n,p}^{^+}(x) \text{ is s.d.}\}
).
\end{equation}

\begin{lemma}\label{basic2}
For all $n,k$ and $p$, with $1\le k \le n$ and
$0 < p < 1$, we have
\begin{equation}\label{prob2}
P_{n,p}(k)=2^{k-1}\binom{n-1}{k-1}
(1-p)^{2kn-3k^2/2-k/2}\Pr(\, G_{k,p}
\text{ \rm is connected}\, ).
\end{equation}
\end{lemma}

\begin{proof}
Let $X$ be a set of $k$ strictly distinct literals
with $x\in X$. For $y, z \in X$, the dual of
the implication $y\to z$ involves
no literal in $X$.  Therefore the probability that
$L^{^+}_{F_{n,p}}(x)=X$ is $P_a \, P_b$,
where $P_a$ is the probability that every vertex
of the random digraph $D_{X,p}$ can be reached
from $x$ and $P_b$ is the probability that the
random 2-SAT formula $F_{n,p}$ contains no
implication from the set
$I(X, X^c)=\{y\to z: \, y\in X, z\not \in X\}$.

By Lemma~\ref{basic1}, we have
$P_a=\Pr(G_{k,p} \text{ is connected})$,
so we turn to the task of
calculating the probability $P_b$.

Note that there are $k(2n-k)$ implications
in the set $I(X, X^c)$. However,
a 2-SAT formula $F_{n,p}$ contains
none of the $k$ implications
$y\to \overline y$, $y\in X$.  Also, if
$y\in X$ and $z\in  \overline X$, then $y\to z$ and
$\overline z \to \overline y$ are dual
implications, i.e.\ $F_{n,p}$ contains
$y\to z$ if and only if it contains
$\overline z \to \overline y$. In fact,
both implications $y\to z$ and
$\overline z \to \overline y$ belong to
$I(X, X^c)$ if and only if $y\in X$,
$z\in \overline X$ and $y\not = \overline z$.
Hence $I(X, X^c)$ contains $(k^2-k)/2$ dual
pairs, so that the probability that
$F_{n,p}$ contains no implication from
$I(X, X^c)$ is
\begin{equation*}
P_b=(1-p)^{k(2n-k)-k -(k^2-k)/2}=(1-p)^{2kn-3k^2/2-k/2}.
\end{equation*}
Therefore,
\begin{equation}
\Pr(\, L_{F_{n,p}}^{^+}(x)=X\, )=
\Pr(\, G_{k,p} \, \text{ \rm is connected}\, )
(1-p)^{2kn-3k^2/2-k/2}.
\end{equation}
Since there are $2^{k-1}\binom{n-1}{k-1}$
choices for the set $X$, the lemma is proved.
\end{proof}

In order to transform Lemma~\ref{basic2} into a form suitable
for applications, note that the probability that
$G_{k,p}$ is connected is trivially expressed in terms of
$f(k, m)$, the
number of connected labelled graphs with $k$ vertices and
$m$ edges, and that $f(k,m)$ has a good and fairly
simple approximation when $m-k$ is not too large compared to $k$.
To state this approximation, let us define an array of numbers
$c_{k,\ell}$ by
$f(k, k-1+\ell)=c_{k,\ell} k^{k-2+3\ell/2}$. The somewhat peculiar
choice of the parameter $\ell$ is justified by the fact that
$c_{k,\ell} \neq 0$ if and only if
$0\le \ell \le \binom{k}{2}-k+1$.
Also, if $\ell$ is not too large then
$f(k,k-1+\ell)$ has order about
$k^{k-2+3\ell/2}$.  More precisely, since $f(k, k-1)$ is just the
number of trees with $k$ labelled vertices, by Cayley's theorem
we have $c_{k,0}=1$. Also,
\begin{equation}
\label{ck1}
c_{k,1}=(1+O(k^{-1/2}))(\pi / 8)^{1/2},
\end{equation}
$c_{k, \ell}\le 1$ for all $k\ge 1$ and $\ell \ge 0$,
and for all $2\le \ell \le \binom{k}{2}-k+1$
and some $c<\infty$, we have
\begin{equation}
\label{genbound}
c_{k, \ell}\le (c/\ell)^{\ell/2},
\end{equation}
see~\cite{Bol85}. When $\ell$ is fairly small compared to $k$,
there are
rather detailed estimates for   $c_{k, \ell}$. To be precise,
Wright~(\cite{Wri77},\cite{Wri80})
showed that for $2\le \ell =o(k^{1/3})$ we have
\begin{equation}\label{specbound}
c_{k, \ell}=\gamma (3\pi)^{1/2}
\bigl( \frac{e}{12(\ell-1)}\bigr)^{(\ell-1)/2}
(1+o_\ell(1)),
\end{equation}
where $\gamma =0.159155\ldots$ is the
the limit of a certain bounded increasing
sequence.  Later Meertens proved that $\gamma=1/(2\pi)$
(see \cite{BCM90}).

The probability that $G_{k,p}$ is connected
and has $k-1+\ell$ edges is just
\begin{equation*}
c_{k, \ell}k^{k-2+3\ell/2}p^{k-1+\ell}(1-p)^{\binom{k}{2}-k+1-\ell},
\end{equation*}
so Lemma~\ref{basic2} has the following immediate consequence.

\newcommand{\spk}{\sqrt{\frac{\pi}{8}}\frac{k^{3/2} p}{1-p}
    \left[1+O(k^{-1/2}) + \O_0\Big(\frac{k^{3/2} p}{1-p}\Big)\right]}

\begin{corollary}\label{ppnk/spk}
For all $n$, $k$ and $p$, with $1\le k \le n$ and $0 < p < 1$, we have
\begin{equation}
\label{bprod1}
P_{n,p}(k)=\frac{1}{n}\binom{n}{k} (2pk)^{k-1} (1-p)^{2kn-k^2-2k+1}
          S_p(k),
\end{equation}
where
\begin{equation}\label{bprod2}
S_p(k)=
\sum_{\ell=0}^{\binom{k}{2}-k+1}c_{k, \ell}
\Bigl(\frac{k^{3/2} p}{1-p}\Bigr)^{\ell }.
\end{equation}
If ${k^{3/2} p}/{(1-p)}$ is bounded, then
\begin{equation}\label{spk}
S_p(k)=1+\spk,
\end{equation}
where $\O_0(\cdot)$ is the Landau symbol introduced at the beginning of
Section \ref{sec:strat}.
\end{corollary}

To estimate $P_{n,p}(k)$, we relate it to known bounds on related
events
in random graphs.
To this end, we recall the definition of the
trimmed out-graph $\widetilde D_{F}^+(x)$
in Section \ref{sec:digraph} and its relation to
the random graph
$G_{n,2p-p^2}$
on $n$ vertices with edge probability $2p-p^2$,
 see Lemma~\ref{trimmed-graph}.
Recall that the vertex set of $\widetilde D_{F_{n,p}}^+(x)$
is denoted by $\widetilde L^+_{n,p}(x)$.
We define
\begin{equation}
Q_{n,p}(k)=\Pr(|\widetilde L^+_{n,p}(x)|=k).
\end{equation}
By Lemma~\ref{trimmed-graph} part (iv),
\begin{equation}
\label{def-Qnp}
Q_{n,p}(k) = \Pr\{ |C_{n,2p-p^2} (x)| = k \} ~,
\end{equation}
where $C_{n,2p-p^2} (x)$ is the connected component in $G_{n,2p-p^2}$
containing a fixed vertex $x$.
For all $n,k$ and $p$, with $1\le k \le n$ and
$0 < p < 1$, we have
\begin{align}
Q_{n,p}(k)
&= \binom{n-1}{k-1} (1-2p+p^2)^{k(n-k)}
    \Pr(\, G_{k,2p-p^2} \text{ \rm is connected}\, ) \nonumber \\
&= \frac{1}{n} \binom{n}{k} ((2p-p^2)k)^{k-1}
(1-p)^{2kn-2k^2+k^2-3k+2} S_{2p-p^2}(k).
\label{Qp}
\end{align}

In Section \ref{sec:hourglass}, we shall also need bounds
on the probability $R_{n,p}(k)$ that
$C_{n,2p-p^2}(x)$ is a tree of size $k$,
\begin{equation}
\label{tree-size-dist}
R_{n,p}(k)=\Pr
\left( |C_{n,2p-p^2}(x)|=k \text{ and }
C_{n,2p-p^2}(x)\text{ is a tree}
\right).
\end{equation}
Recalling the derivation of Corollary \ref{ppnk/spk},
we immediately see that $Q_{n,p}(k)$ and $R_{n,p}(k)$
are related by
\begin{equation}
\label{q/r}
Q_{n,p}(k)=R_{n,p}(k) S_{2p-p^2}(k).
\end{equation}

The following lemma will be used to turn well-known
bounds on $Q_{n,p}(k)$ into bounds
on $P_{n,p}(k)$.

\addtolength{\textheight}{1\baselineskip}

\begin{lemma}\label{qpdiff}
For all $0<p<1$,
\begin{equation}
\label{P<Q}
P_{n,p}(k)\leq
 Q_{n,p}(k).
\end{equation}
If $p\leq 1/2$ and $k^{3/2}p$ is bounded, then
\begin{equation}
\label{diff2}
0\leq {Q_{n,p}(k)}-{P_{n,p}(k)}
= \spk {P_{n,p}(k)}~.
\end{equation}
If $p\leq 1/2$ and $k^{3/2}p\geq 1$, then
\begin{equation}\label{P<<Q}
P_{n,p}(k) = O\Big(\ell_0 2^{-\ell_0} Q_{n,p}(k)\Big),
\end{equation}
where
\begin{equation}\label{ellzero}
\ell_0=\ell_0(k)=\min \left\{\frac{k^3 p^2}{12(1-p)^2},
n^{1/5} \right\}.
\end{equation}
\end{lemma}

\begin{proof}
By Lemma~\ref{trimmed-graph},
$L^+_{F_{n,p}}(x)$ is strictly distinct
if and only if $L^+_{F_{n,p}}(x)=\widetilde L^+_{F_{n,p}}(x)$.
As a consequence,
\begin{align}
P_{n,p}(k)
&=
\Pr(L^+_{F_{n,p}}(x)=\widetilde L^+_{F_{n,p}}(x)
\text{ and }|L^+_{F_{n,p}}(x)|=k)
\notag
\\
&=
\Pr(L^+_{F_{n,p}}(x)=\widetilde L^+_{F_{n,p}}(x)
\text{ and }|\widetilde L^+_{F_{n,p}}(x)|=k)
\notag
\\
&\leq
\Pr(|\widetilde L^+_{F_{n,p}}(x)|=k)
=Q_{n,p}(k),
\end{align}
which proves \eqref{P<Q}.
Rewriting (\ref{Qp}) in the form
\begin{align}
 Q_{n,p}(k)
&= \frac{1}{n} \binom{n}{k} (2pk)^{k-1} (1-p/2)^{k-1}
(1-p)^{2kn-k^2-3k+2} S_{2p-p^2}(k)
\notag
\\
&= P_{n,p}(k) (1-p/2)^{k-1} (1-p)^{-k+1} S_{2p-p^2}(k) / S_p(k)
\label{qvsp} \\
&\geq P_{n,p}(k) S_{2p-p^2}(k) / S_p(k) \label{qvsp2}
\end{align}
and observing that $S_p(k)$ is a monotone increasing function of
$p/(1-p)$, and hence of $p$ (see (\ref{bprod2})), we have
$S_{2p-p^2}(k)
\geq S_p(k)$, obtaining an alternative proof of the bound (\ref{P<Q}).

If $p\leq 1/2$ and $k^{3/2}p$ is bounded,
then both $k^{3/2}p/(1-p)$ and $k^{3/2}(2p-p^2)/(1-2p+p^2)$ are
bounded by a constant times $k^{3/2}p$.  By
(\ref{spk}), we therefore have
\begin{equation*}
\frac{S_{2p-p^2}(k)}{S_p(k)} = 1+\spk
\end{equation*}
and hence by (\ref{qvsp})
\begin{equation}
Q_{n,p}(k) = P_{n,p}(k) (1+O(k p)) \left(1+\spk\right).
\label{q/p}
\end{equation}
As a consequence,
\begin{equation*}
\frac{Q_{n,p}(k)}{P_{n,p}(k)}-1 = \spk
\notag
\end{equation*}
which implies the bound (\ref{diff2}).

In order to prove (\ref{P<<Q}), we
decompose $S_p(k)$ as
\begin{align*}
S_p(k)
=\sum_{\ell \ge 0} c_{k, \ell} \Big(\frac{k^{3/2} p}
{1-p}\Big)^{\ell}
&=
\sum_{0 \le \ell < \ell_0} c_{k,\ell} \Big(\frac{k^{3/2} p}
{1-p}\Big)^{\ell}
+\sum_{\ell \ge \ell_0} c_{k,\ell} \Big(\frac{k^{3/2} p}{1-p}
\Big)^{\ell}\\
&=S'_p(k)+ S_p''(k)~,
\end{align*}
where $S'_p(k)$ is the first and $S''_p(k)$ is the second sum above.
If $\ell \le n^{1/5}$ then $\ell =o(k^{1/3})$, so by
(\ref{specbound})
we have
\begin{equation*}
S_p'(k)=O\Big( \sum_{2\le \ell < \ell_0}
\Big(\frac{e^{1/2}k^{3/2} p}{12^{1/2}(\ell-1)^{1/2}(1-p)}
\Big)^{(\ell -1)}
k^{3/2}p\Big).
\end{equation*}
By taking the ratio of successive terms in the series
\begin{align*}
\Big(\frac{e^{1/2} k^{3/2} p}{12^{1/2} (\ell \!-\!1)^
{1/2}(1-p)}\Big)^{\!-(\ell-1)}
\Big(\frac{e^{1/2} k^{3/2} p}{12^{1/2}  \ell    ^
{1/2}(1-p)}\Big)^{\!\ell}\!
&=
  \Big(\frac{e^{1/2} k^{3/2} p}{12^{1/2} (1-p)}\Big)
  \frac{1}{\ell^{1/2}}
  \Big(1 + \frac{1}{\ell \!-\!1}\Big)^{\!-(\ell-1)/2}
  \\
&\geq
  \Big(\frac{e^{1/2} k^{3/2} p}{12^{1/2} (1-p)}\Big)
  \frac{1}{\ell^{1/2}}
  e^{-1/2},
\end{align*}
we see that the summand is an increasing function of
$\ell$ for $2 \le
\ell \le \ell_0+1.$
Since the $(\ell_0+1)$st term is in the sum for $S_p''(k)$,
we then have
$$ S_p'(k) = O(\ell_0 S_p''(k))~.$$

To bound $S_p''(k)$, we note that
\begin{align*}
\frac{2 p}{1-p} = \frac{2 p - p^2}{1-(2p-p^2)} \\
\end{align*}
so that
\begin{align*}
S''_p(k) &= \sum _{\ell \ge \ell_0}
c_{k, \ell} \Big(\frac{k^{3/2}p}{1-p}\Big)^{\ell}\\
&\le  2^{-\ell_0}\sum_{\ell \ge \ell_0} c_{k, \ell}
\Big(\frac{k^{3/2}2p}{1-p}\Big)^{\ell}\\
&=  2^{-\ell_0}\sum_{\ell \ge \ell_0} c_{k, \ell}
\Big(\frac{k^{3/2}(2p-p^2)}{1-(2p-p^2)}\Big)^{\ell}\\
&=  2^{-\ell_0} S''_{2p-p^2}(k)~,
\end{align*}
and therefore
\begin{equation*}
S_p(k)=O\big(\ell_0 2^{-\ell_0} S''_{2p-p^2}(k)\big)=O\big
(\ell_0 2^{-\ell_0} S_{2p-p^2}(k)\big).
\end{equation*}
Combined with (\ref{qvsp2}), this gives
the desired bound (\ref{P<<Q}).
\end{proof}
\medskip

Our later estimates rely heavily on bounds on the {\em expectation}
of the size of the component of a given vertex in a random graph, Lemma
\ref{weakgiant} below.
Although sharper forms of these bounds were already proved by
Bollob\'as~\cite{Bol84}, \L uczak~\cite{Luc90}, and Janson {\it et
al.\/}~\cite{JKLP94}, these previous estimates were proved only for a
restricted range of $\lambda$ (e.g.\ $\lambda\leq n^{1/12}$), whereas
we require the full range (i.e.\ $\lambda_n \le \varepsilon_0n^{1/3}
$).  These estimates turn out to be rather involved; the proof is given
in Appendix~\ref{sec:appendix}.

We remark that the bound we shall
obtain in (\ref{bigp11}) below is closely related to the order of
the giant component:
if
$\lambda_n \to \infty$ and $\lambda_n=o(n^{1/3})$ then, for
$p=(1+\lambda_n n^{-1/3})/n$, with probability tending to $1$, the
random graph $G_{n, p}$ has a unique giant component with
$(2+o(1))\lambda_n n^{2/3}$ vertices.

\begin{lemma}\label{weakgiant}
There are constants $c$,  $\varepsilon_0$
 and $\lambda_0$, $c>0$, $0<\varepsilon_0<1$,
and $0<\lambda_0<\infty$  with
the following property.
Let $\lambda_0 \le \lambda_n \le \varepsilon_0n^{1/3} $ and,
as before, set
$\varepsilon=\lambda_n n^{-1/3}$
and
$p=\big(1 + \lambda_n n^{-1/3}\big)/2n$.
 Then
\begin{equation}\label{bigp11}
\sum_{k\ge \lambda_n  n^{2/3} }
Q_{n,p}(k) =
\vartheta(\varepsilon)(1+O({1/\lambda_n^2}))
\end{equation}
and
\begin{equation}\label{bigp22}
\sum_{k=\lfloor n^{2/3}/\lambda_n\rfloor}^{\lceil\lambda_n n^{2/3}
\rceil}
Q_{n,p}(k)
=  O(e^{-c\lambda_n} n^{-1/3}).
\end{equation}
\end{lemma}

\section{Moment Estimates}
\label{sec:moments}
\addtolength{\textheight}{-1\baselineskip}

In this section we shall bound
the moments of the number of literals $x$
for which
the set $L_{n,p}^{^+}(x)$ consists of
strictly distinct literals.
Recall that $P_{n,p}(k)$ is the probability
that, for a variable $x$,
the set $L_{n,p}^{^+}(x)$ consists of
$k$ strictly distinct literals.
We use Corollary~\ref{ppnk/spk} to get a good estimate of $P_{n,p}(k)$.
Note that
\begin{align}
\label{new1}
\frac{1}{n}\binom{n}{k}(2pk)^{k-1}
 &=\frac{1}{2pnk}\frac{(k/e)^k}{k!} (2npe)^k
  \prod_{i=0}^{k-1}\Big(1-\frac{i}{n}\Big)\notag\\
 &= \frac{\exp[-1/(12 k +\delta_k)]}{2 p n \sqrt{2\pi } k^{3/2}} (2npe)^k
  \prod_{i=0}^{k-1}\Big(1-\frac{i}{n}\Big),
\notag
\end{align}
where $0\leq\delta_k\leq1$ for each $k$.
Recalling that $2np = 1+\varepsilon = 1+{\lambda_n} n^{-1/3}$,
with $|\varepsilon|\leq \varepsilon_0$,
$\varepsilon_0<1$,
 we write
\begin{equation*}
\log (1+\varepsilon)
  = \varepsilon-\frac12\varepsilon^2 + \Theta(\varepsilon^3)
\end{equation*}
and
\begin{equation*}
\log\Big(1-\frac in\Big) =-\frac in-\frac{i^2}{2n^2}-
                           \Theta\Big(\frac{i^3}{n^3}\Big)
\end{equation*}
and obtain
\begin{align*}
\frac{1}{n}\binom{n}{k}(2pk)^{k-1}
 =
 \frac{\exp[-\Theta(1/k)]}{2pn \sqrt{2\pi } k^{3/2}}
       &\exp\Big\{k+ k\varepsilon-\frac{1}{2}k\varepsilon^2 +k
\Theta(\varepsilon^3)\Big\}\\
\times &\exp\Big\{-\frac{k^2}{2n}   + \Theta\Big(\frac{k}{n}\Big)
                  -\frac{k^3}{6n^2} + \Theta\Big(\frac{k^2}{n^2}\Big)
                  -\Theta\Big(\frac{k^4}{n^3}\Big)
            \Big\}.
\end{align*}
Expressing
\begin{align*}
(1-p)^{2kn-k^2-2k+1}
  &= \exp \Big\{-(p+\Theta(p^2))\big(2kn-k^2-2k+1\big)\Big\} \\
  &= \exp \Big\{-k-\varepsilon k + \frac{k^2}{2n}+\frac{
\varepsilon k^2}{2 n}+\frac {k(1+\varepsilon)}{n} - \Theta(p^2 k n)\Big\},
\end{align*}
we therefore get
\begin{align}
P_{n,p}(k) =&
  \frac{1}{2pn\sqrt{2\pi}\ k^{3/2}} \times \notag \\
  &\exp \Big\{
       -\frac{k\varepsilon^2}{2}
       -\frac{k^3}{6 n^2}
       +\frac{\varepsilon k^2}{2 n}
                  + O\Big(\frac{k}{n}\Big)
                  -\Theta\Big(\frac{k^4}{n^3}\Big)
                  +k \Theta(\varepsilon^3)
       -\Theta(1/k)
       \Big\}
  \,S_p(k).
 \label{ppnk-asymp}
\end{align}

\begin{lemma}
\label{lem:pnk}
If $k/n^{2/3}$ is bounded, then
$P_{n,p}(k)$, $Q_{n,p}(k)$ and $R_{n,p}(k)$
can be rewritten in the form
\begin{equation}
\label{pnk-usable}
  \frac{1}{2pn\sqrt{2\pi}\ k^{3/2}}
  \exp \Big\{
       -\frac{k\varepsilon^2}{2}\Big(1-\Theta(\varepsilon)-\Theta\Big(\frac{kn^{-2/3}}
                  {\lambda_n}\Big)\Big)
                  + \O_0\Big(\frac{k^{3/2}}{n}\Big)
                  -\Theta(1/k)
       \Big\}.
\end{equation}
\end{lemma}

\begin{remark}
\label{much_less_than}
It is evident that the multiplicative error terms
$1-\Theta(\varepsilon)-\Theta(kn^{-2/3}/{\lambda_n})$ above can be made arbitrarily
close to $1$ if
$\varepsilon$ is small enough and ${\lambda_n}$ is large enough,
i.e.\ if
$1\ll{\lambda_n}\ll n^{1/3}$, or equivalently $n^{-1/3}\ll
\varepsilon\ll 1$.
Likewise the additive error terms can be made arbitrarily
small when
$1\ll k\ll n^{2/3}$.
We follow the standard convention that ``if $a\ll b$ then
$c=(1+o(1))d$'' means that the $o(1)$ error term can be made
arbitrarily small by taking the ratio $b/a$ large enough, and in our
earlier notation may be re-expressed as ``$c=(1+o_{b/a}(1))d$.''
\end{remark}

\begin{proof_of}{Lemma~\ref{lem:pnk}}
With the assumption that $k/n^{2/3}$ is bounded, we can use
our bound (\ref{spk}) on $S_p(k)$ together with the bounds
$k/n\leq k^{3/2}/n$, $k^4/n^3 = O(k^{3}/n^2)=\O_0(k^{3/2}/n)$, and
$k^3/n^2 = \O_0(k^{3/2}/n)$ to rewrite (\ref{ppnk-asymp}) as
\begin{equation}
\label{pnk-intermediate}
P_{n,p}(k) =
  \frac{1}{2pn\sqrt{2\pi}\ k^{3/2}}
  \exp \Big\{
       -\frac{k\varepsilon^2}{2}(1-\Theta(\varepsilon))
       +\frac{\varepsilon k^2}{2 n}
                  + \O_0\Big(\frac{k^{3/2}}{n}\Big)
                  -\Theta(1/k)
       \Big\}.
\end{equation}
Note that the ratio between the second and first terms in
the exponential
is $k/(n\varepsilon) = k/({\lambda_n} n^{2/3})$,
so we can further
rewrite (\ref{pnk-intermediate}) as (\ref{pnk-usable}).
The estimate \eqref{pnk-usable} for $Q_{n,p}(k)$
follows immediately from
the estimate for $P_{n,p}(k)$ and \eqref{q/p},
and the estimate for $R_{n,p}(k)$ follows from
that for $Q_{n,p}(k)$ and \eqref{q/r}.
\end{proof_of}

\begin{lemma}\label{RC-moments}
There are constants $c$,
$\varepsilon_0$ and $\lambda_0$ with $c>0$,
$0<\varepsilon_0<1$ and $0<\lambda_0<\infty$,
such that the following statements hold
for $a>1/2$ and
$\lambda_0\leq |\lambda|\leq \varepsilon_0 n^{1/3}$.

i) If $\lambda_n<0$, then
\begin{equation}\label{RC-moment1}
\sum_{k \ge n^{2/3}/ |\lambda_n| } k^a Q_{n,p}(k)
=O\left(
      \left[\frac{2}{\varepsilon^2}\right]^{a-1/2} e^{-c|\lambda_n|}
\right),
\end{equation}
where the constant implicit in the Landau symbol $O(\cdot)$ depends on $a$.

ii) For both positive and negative $\lambda_n$, we have
\begin{align}
\label{RC-moment2}
\sum_{k\leq n^{2/3}/|\lambda_n| } k^a\ Q_{n,p}(k)
&= \frac{1+\o(1)}{2pn}\,
  \frac{\Gamma(a-1/2)}{\sqrt{2\pi}}
     \left[\frac{2}{\varepsilon^2}\right]^{a-1/2},
\end{align}
where the $\o(1)$ term depends upon $a$, and for fixed $a$
becomes as small
as we like if $|{\lambda_n}|$ is large enough and
$\varepsilon={\lambda_n}/n^{1/3}$ is small enough.
\end{lemma}
\begin{remark}
Here, as in the rest of this paper, $\Gamma$
denotes the gamma function, which
interpolates factorials.
Recalling that $\Gamma(1/2)=\sqrt{\pi}$,
the lemma immediately implies that below the window the
expected size a component of a given
vertex in $G_{n,2p-p^2}$ can be
estimated by
\begin{equation}
\label{RG-exp}
E(|C_{n,2p-p^2}(x)|)
=\frac 1{{2pn}|\varepsilon|}(1+\o(1)).
\end{equation}
Note also that \eqref{RC-moment1} implies that
below the window
\begin{equation}\label{RC-moment-modified}
\sum_{k \ge n^{2/3}/ |\lambda_n| }  Q_{n,p}(k)
= \sum_{k \ge n^{2/3}/ |\lambda_n| } \frac{k}{k} Q_{n,p}(k)
\leq \frac{\lambda_n}{n^{2/3}}
\sum_{k \ge n^{2/3}/ |\lambda_n| } k Q_{n,p}(k)
=O\left(
      n^{-1/3} e^{-c|\lambda_n|}\right).
\end{equation}

\end{remark}

\begin{proof_of}{Lemma~\ref{RC-moments}}
We start with the proof of ii).  The limit of the summation is
$n^{2/3}/|\lambda| = n^{1/3}/|\varepsilon|$; we first sum the portion up to
$1/\varepsilon$ separately.  For this portion, the sum is
\begin{equation}\label{small-k}
\sum_{k < 1/|\varepsilon|} k^a\ Q_{n,p}(k) = O((1/\varepsilon)^{a-1/2}),
\end{equation} where here and throughout this proof,
constants may depend upon $a$.

Referring to our expression \eqref{pnk-usable}
for $Q_{n,p}(k)$, we see that for the
remainder of the terms in the sum, the additive error terms in the
exponential tend to zero when $k/n^{2/3} \leq 1/|\lambda|$ is small
enough and $k\geq 1/|\varepsilon|$ is large enough, and so they
contribute only a multiplicative error of $1+o(1)$.  The
multiplicative error terms in the exponential can also be made
arbitrarily close to $1$ by taking $|\lambda|$ large enough and
$|\varepsilon|$ small enough.  Therefore we have
\begin{equation}
  \frac{1+\o(1)}{2pn\sqrt{2\pi}} B(2-\delta) \leq
  \sum_{1/|\varepsilon| \leq k \leq n^{1/3}/|\varepsilon|} k^a\ Q_{n,p}(k) \leq
 \frac{1+\o(1)}{2pn\sqrt{2\pi}} B(2+\delta)
 \label{B2+-}
\end{equation}
where $\delta = \o(1)$ and
\begin{align}
B(t) &= \sum_{1/|\varepsilon| \leq k \leq n^{1/3}/|\varepsilon|}
      k^{a-3/2}
      \ \exp\left[-\frac{k\varepsilon^2}{t}\right]. \label{Bt} \\
\intertext{Since the summand in \eqref{Bt} is unimodal,
we may approximate the sum with an integral}
B(t) &= \int_{1/|\varepsilon|}^{n^{1/3}/|\varepsilon|}
      k^{a-3/2} \ \exp\left[-\frac{k\varepsilon^2}{t}\right] dk
 + O\left(\max k^{a-3/2}\ \exp\left[-\frac{k\varepsilon^2}{t}\right]\right)\notag\\
 &= \int_{|\varepsilon|/t}^{n^{1/3}|\varepsilon|/t}
      \left[\frac{t u}{\varepsilon^2}\right]^{a-3/2} \ \exp[-u] \frac{t\ du}{\varepsilon^2}
 + O\left(\left[\frac{1}{\varepsilon^2}\right]^{a-3/2}\right)\notag\\
 &= \left[\frac{t}{\varepsilon^2}\right]^{a-1/2}
    \int_{|\varepsilon|/t}^{|\lambda|/t}
      u^{a-3/2} \ \exp[-u] \ du
 + O\left(\left[\frac{1}{\varepsilon^2}\right]^{a-3/2}\right)\notag\\
 &= \left[\frac{t}{\varepsilon^2}\right]^{a-1/2}
    \left[ \Gamma(a-1/2) +\o(1) \right]
 + O\left(\left[\frac{1}{\varepsilon^2}\right]^{a-3/2}\right)\notag\\
 &= \left[\frac{t}{\varepsilon^2}\right]^{a-1/2}
    \left[ \Gamma(a-1/2) +\o(1) \right]\label{Bt-G}
\end{align}
where the $\o(1)$ depends on $a$.
Combining \eqref{small-k}, \eqref{B2+-} and \eqref{Bt-G},
we get \eqref{RC-moment2}.

To prove i),
we
use the well know fact that the cluster size distribution
in $G_{n,\widetilde p}$, with $\widetilde p=2p-p^2$, is stochastically
dominated by a birth process with binominal offspring
distribution $\Bin(n,\widetilde p)$, which in
turn is stochastically dominated by a Poisson birth process
with  parameter
$
n \log [1/({1-\widetilde p})]=2np(1+O(1/n)).
$
Writing this parameter as
$1+\widetilde\varepsilon$,
we therefore get the estimate

\begin{align}
\label{RC-moments6}
\sum_{k> \frac{n^{2/3}}{|\lambda_n|}}
 k^a Q_{n,p}(k)
&\leq
\frac 1{1+\widetilde\varepsilon}
\sum_{k\geq n^{2/3}/|\lambda_n|}
k^a\,
\frac{k^{k-1}}{k!}
(({1+\widetilde\varepsilon})e^{-({1+\widetilde\varepsilon})})^k
\notag
\\
&\leq
\frac 1{1+\widetilde\varepsilon}
\sum_{k\geq n^{2/3}/|\lambda_n|}^{n}
\frac{k^a}{\sqrt{2\pi k^3}}
\exp\Big({-\frac{\widetilde\varepsilon^2}2k}\Big)
\notag
\\
&= \frac{1+O(1/n)}{2pn\sqrt{2\pi}}
      \sum_{k\geq \frac{n^{2/3}}{|\lambda|}}
      k^{a-3/2}
      \ \exp\left\{-\frac{k\varepsilon^2}{2}(1+o_n(1)) \right\}.
\end{align}
Bounding the sum on the right hand side by
$O\left(
      \left[({ n^{2/3}}){\lambda_n^{-2}}\right]^{a-1/2} e^{-c|\lambda_n|}
\right)$,
we obtain the estimate \eqref{RC-moment1}.
\end{proof_of}

\begin{lemma}\label{moments}
There are constants $c$,
$\varepsilon_0$ and $\lambda_0$ with $c>0$,
$0<\varepsilon_0<1$ and $0<\lambda_0<\infty$,
such that
for $a>1/2$ and
$\lambda_0\leq |\lambda|\leq \varepsilon_0 n^{1/3}$,
we have
\begin{align}
\label{moment}
\sum_{k=1}^{n} k^a\ P_{n,p}(k)
&= \frac{1+\o(1)}{2pn}
  \frac{\Gamma(a-1/2)}{\sqrt{2\pi}}
     \left[\frac{2}{\varepsilon^2}\right]^{a-1/2}
\end{align}
and
\begin{align}
\label{moment1}
\sum_{k\geq n^{2/3}/|\lambda_n| } k^a\ P_{n,p}(k)
&=
\left[\frac{1}{\varepsilon^2}\right]^{a-1/2}
\begin{cases}
O(e^{-\Theta(|\lambda_n|)})&\text{if}\quad\lambda_n<0
\\
O(e^{-\Theta(\lambda_n^{-3/5})})&\text{if}\quad\lambda_n>0,
\end{cases}
\end{align}
where the constants implicit in
both the $\o(1)$ term and the $O$ term depends upon $a$,
and where $\varepsilon = \lambda_n n^{-1/3}$ as usual.
\end{lemma}

\begin{remark}
The particular values of $a$ that we shall need are $a=1$,
$a=3/2$, and $a=2$.
For these values, we have
\begin{equation}
\label{moment-a}
  \sum_k k^a P_{n,p}(k) = \frac{1+\o(1)}{2pn}\times
  \begin{cases} |\varepsilon|^{-1} & a=1\\
                \sqrt{2/\pi}\ \varepsilon^{-2} & a=3/2\\
                |\varepsilon|^{-3} & a=2 ~,
  \end{cases}
\end{equation}
where as before the summation can range over all $k$, or over all
$k\leq n^{2/3}/|\lambda_n|$.
\end{remark}

\begin{remark}
Using the same trick as in \eqref{RC-moment-modified}, we can write
\begin{equation}
\sum_{k \ge n^{2/3}/ |\lambda_n| }  P_{n,p}(k)
= \sum_{k \ge n^{2/3}/ |\lambda_n| } \frac{k}{k} P_{n,p}(k)
= n^{-1/3} \times
\begin{cases}
O(e^{-\Theta(|\lambda_n|)})&\text{if}\quad\lambda_n<0
\\
O(e^{-\Theta(\lambda_n^{-3/5})})&\text{if}\quad\lambda_n>0. \label{big-pnk}
\end{cases}
\end{equation}
\end{remark}

\begin{proof_of}{Lemma~\ref{moments}}
To prove the lemma, we again consider three regions of $k$:
$k\leq n^{2/3}/|\lambda_n|$,
$n^{2/3}/|\lambda_n| \leq k \leq n^{2/3} |\lambda_n| $, and
$k \geq n^{2/3} |\lambda_n| $.
In the first region, we use
\eqref{q/p} to approximate $P_{n,p}(k)$
by $Q_{n,p}(k)$.
Combined with Lemma~\ref{RC-moments}, this gives
\begin{align}
\label{moment.1}
  \sum_{k\leq \frac{n^{2/3}}{|\lambda|}} \! k^a P_{n,p}(k)
   &= \frac{1+\o(1)}{2pn\sqrt{2\pi}}
     \left[\frac{2 n^{2/3}}{{\lambda_n}^2}\right]^{a-1/2}
\Gamma(a-1/2).
\end{align}
For $\lambda_n<0$, we combine the second and the third region.
Using the fact that $P_{n,p}(k)\leq Q_{n,p}(k)$
by \eqref{P<Q}, we then use Lemma~\ref{RC-moments} to
obtains the bound \eqref{moment1} below threshold.

Above threshold, the contribution from the second region is
\begin{align}
  \sum_{n^{2/3}/|\lambda_n |\leq k \leq n^{2/3} |\lambda_n| }
k^a P_{n,p}(k)
  &\leq
  \sum_{n^{2/3}/|\lambda_n |\leq k \leq n^{2/3} |\lambda_n| }
k^a Q_{n,p}(k)
  \tag*{by (\ref{P<Q})\quad}
\\
  &\leq
  (n^{2/3} |\lambda_n|)^a
  \sum_{n^{2/3}/|\lambda_n |\leq k \leq n^{2/3} |\lambda_n| }
Q_{n,p}(k)
\notag
\\
  &\leq
  (n^{2/3} |\lambda_n|)^a
  \exp[-\Theta(|{\lambda_n}|)] n^{-1/3}
  \tag*{by (\ref{bigp22})\quad}
\\
  &=
  \left[\frac{2 n^{2/3}}{\lambda_n^2} \right]^{a-1/2} \times
  \exp[-\Theta(|{\lambda_n}|)] ~.
\label{moment2}
\end{align}

In the third region we use (\ref{P<<Q}) to write
\begin{equation}\label{high-moment}
  \sum_{k \geq n^{2/3} |\lambda_n| } k^a P_{n,p}(k)
  =
  O\bigg( \sum_{k \geq n^{2/3} |\lambda_n| }
  k^a\ \ell_0(k)\ 2^{-\ell_0(k) }\ Q_{n,p}(k)\bigg) ~.
\end{equation}
Consider the prefactor, i.e.\ the summand ignoring the
factor of $Q_{n,p}(k)$.
Provided $\ell_0(k+1)<n^{1/5}$, recalling the definition of
$\ell_0$ in (\ref{ellzero}), we can differentiate the logarithm
of the
prefactor to get
$$ \frac{d}{dk} \log\left[k^{a+3} p^2/(12(1-p)^2) 2^
{-k^3 p^2/(12(1-p)^2)} \right]
=
  \frac{a+3}{k} - (\log 2) 3k^2 \frac{p^2}{12(1-p)^2}
$$
which will be nonpositive provided
$$  \frac{4(a+3)}{\log 2} \leq k^3 p^2 = \frac{k^3}{4n^2}
(1+\varepsilon)^2,$$
which will hold if $|{\lambda_n}|$ is large enough,
where ``large enough''
depends upon $a$.
Provided that these conditions are met, the prefactor
takes on its
largest value in the first term, $k=n^{2/3}|{\lambda_n}|$.
Eventually, if
$k$ gets large enough, $\ell_0(k)$ may stop increasing,
and take on
the value $n^{1/5}$.  If $a\geq0$ then the prefactor will
then increase up to
$n^{a+1/5} 2^{-n^{1/5}}$.  Thus the maximum value of the
prefactor is
$$ \max\left\{(n^{2/3}{\lambda_n})^a \frac{{\lambda_n}^3}{(48+\o(1))}
2^{-{\lambda_n}^3/(48+\o(1))},\ \
               n^{a+1/5} 2^{-n^{1/5}}\right\}. $$
Since by assumption $\varepsilon={\lambda_n} n^{-1/3} \leq 1$, we have
$n\geq \Theta({\lambda_n}^3)$, so the second term in the $\max$ is
$n^a \exp[-\Theta(n^{1/5})] \exp[-\Theta({\lambda_n}^{3/5})]$.
Thus we can write the maximum prefactor as
$$ \max\left\{n^{(2/3)a} \exp[-\Theta({\lambda_n}^3)],\ \
               n^a \exp[-\Theta(n^{1/5})] \exp[-\Theta
({\lambda_n}^{3/5})]
       \right\}
       \leq n^{(2/3)a} \exp[-\Theta({\lambda_n}^{3/5})]. $$
Upon substituting the above expression into (\ref{high-moment}), we find
\begin{align}
  \sum_{k \geq n^{2/3} |\lambda_n| } k^a\ P_{n,p}(k)
 &\leq
  \sum_{k \geq n^{2/3} |\lambda_n| } \
   n^{(2/3)a}
  \exp[-\Theta({\lambda_n}^{3/5})]  Q_{n,p}(k)
\notag
\\
 &\leq
  n^{(2/3)a -1/3}
  \exp[-\Theta({\lambda_n}^{3/5})]
\notag
\\
 &=
  \left[\frac{2 n^{2/3}}{{\lambda_n}^2} \right]^{a-1/2}
  \exp[-\Theta({\lambda_n}^{3/5})]
\label{moment3}
\end{align}
where we used  (\ref{bigp11})
 to get the second line.

The bounds \eqref{moment2}
and \eqref{moment3} imply \eqref{moment1} above threshold,
and \eqref{moment1} and \eqref{moment.1}
imply the estimate \eqref{moment}.
\end{proof_of}

\section{Expected Size of Spine} \label{sec:spine}

After the preparation of the last three sections, we are ready to prove
the bounds (\ref{left-spine-exp}) and
(\ref{right-spine-exp}) in
Theorems~\ref{spine-exp}.

\Proof{Proof of Theorem~\ref{spine-exp}.}
\nobreak
 We begin with a derivation which
holds for $\varepsilon$ of both signs, and later distinguish $\varepsilon < 0$
from $\varepsilon > 0$.
Since $\sum_{k\ge 1} Q_{n,p}(k)=1$, we have
\begin{equation*}
\Pr
\Big(
{x\underset {F_{n,p}}\rightsquigarrow \overline x}
\Big)
=  1- \sum_{k\geq 1} P_{n,p}(k)
= \sum_{k} \big(Q_{n,p}(k)- P_{n,p}(k)\big) .
\end{equation*}
By (\ref{diff2}),
\begin{align*}
 \sum_{k \leq \frac{n^{2/3}}{|\lambda|}} \!\big(Q_{n,p}(k)-
P_{n,p}(k)\big)
&= \sqrt{\frac{\pi}{8}}
   \sum_{k \leq \frac{n^{2/3}}{|\lambda|}}
     \frac{p}{1-p}
     k^{3/2}  \left[1+O(k^{-1/2}) + O\Big(\frac{k^{3/2} p}
{1-p}\Big)\right]
     P_{n,p}(k) \\
&= \sqrt{\frac{\pi}{8}}
   (1+o_{\lambda_n}(1))
    \sum_{k \leq \frac{n^{2/3}}{|\lambda|}}
        p \left[k^{3/2} +O(k)\right] P_{n,p}(k) \\
&= \sqrt{\frac{\pi}{8}}
   (1+o_{\lambda_n}(1))
   \left[(1+\o(1))\frac{p}{2pn}\frac{\sqrt{2/\pi}}{\varepsilon^2}
         +O(p/\varepsilon)
   \right] \tag*{by (\ref{moment-a})} \\
&= \frac{1+\o(1)}{4 n \varepsilon^2}= \frac{1+\o(1)}{4 n^{1/3}
{\lambda_n}^2}.
\end{align*}

So far everything we have done in this section holds for both the
subcritical region and the supercritical region.
We now consider these regions separately.

\Proof{Subcritical regime:}

We use (\ref{RC-moment-modified}) to sum over the remaining values of $k$:
\begin{equation*}
0\leq
\sum_{k > n^{2/3} / |\lambda_n| }
\big(Q_{n,p}(k)- P_{n,p}(k)\big) \leq
\sum_{k > n^{2/3} / |\lambda_n| }
Q_{n,p}(k) =O(e^{-c|\lambda_n|}n^{-1/3})~.
\end{equation*}
Putting these two ranges together we get
\begin{align*}
 \Pr
\Big(
{x\underset {F_{n,p}}\rightsquigarrow\overline x}
\Big)
&=
 \sum_{k}
     \big(Q_{n,p}(k)- P_{n,p}(k)\big) \\
&=
  \sum_{k \leq n^{2/3} / |\lambda_n| }
     \big(Q_{n,p}(k)- P_{n,p}(k)\big)
   +
  \sum_{ n^{2/3} / |\lambda_n| < k}
     \big(Q_{n,p}(k)- P_{n,p}(k)\big) \\
&=
 \frac{1+\o(1)}{4 n^{1/3} {\lambda_n}^2}
 +
 O\left(\frac{e^{-c|\lambda_n|}}{n^{1/3}}\right)
\\
&=
 \frac{1+\o(1)}{4 n^{1/3} {\lambda_n}^2},
\end{align*}
which completes the proof of \eqref{left-spine-exp}.
\lqed

\Proof{Supercritical regime:}

If $k$ is ``mid-size,'' we can apply (\ref{P<Q}) and
(\ref{bigp22}):
\begin{equation*}
0\leq
\sum_{ n^{2/3} / \lambda_n < k \leq \lambda_n n^{2/3}}
     \big(Q_{n,p}(k)- P_{n,p}(k)\big)
  \leq
\sum_{ n^{2/3} / \lambda_n < k \leq \lambda_n n^{2/3}}
     Q_{n,p}(k)
   =O(e^{-c\lambda_n}/n^{1/3})~.
\end{equation*}

If $k$ is ``large,'' we proceed as in the proof of
\eqref{moment3} to write
\begin{align*}
\sum_{ k>\lambda_n n^{2/3}}\!
     \big(Q_{n,p}(k)- P_{n,p}(k)\big)
 &=
   \big(1-O(\exp[-\Theta({\lambda_n}^{3/5})])\big) \sum_{ k>\lambda_n n^{2/3}}
     Q_{n,p}(k)
 \\
 &= \vartheta(\varepsilon) (1+o_{{\lambda_n}}(1))
  \tag*{\text{by (\ref{bigp11})}}\\
 &= 2 {\lambda_n} n^{-1/3} (1+o_{{\lambda_n}}(1)+O(\varepsilon))~.
\end{align*}

Putting these three ranges of $k$ together, we find
\begin{align*}
\Pr
\Big(
{x\underset {F_{n,p}}\rightsquigarrow \overline x}
\Big)
&=
\sum_{k}\big(Q_{n,p}(k)- P_{n,p}(k)\big) \\
&=
 \left[
  \sum_{k \leq n^{2/3} / \lambda_n }
   +
  \sum_{ n^{2/3} / \lambda_n < k \leq \lambda_n n^{2/3}}
   +
  \sum_{ k>\lambda_n n^{2/3}}
 \right]
     \big(Q_{n,p}(k)- P_{n,p}(k)\big) \\
&=
 \frac{1+\o(1)}{4 n^{1/3} {\lambda_n}^2}
 +
 O\left(\frac{e^{-c\lambda_n}}{n^{1/3}}\right)
+
\vartheta(\varepsilon) (1+o_{{\lambda_n}}(1))
\\
&=
\vartheta(\varepsilon) (1+o_{{\lambda_n}}(1))
\\
&=
 2 {\lambda_n} n^{-1/3} (1+o_{{\lambda_n}}(1)+O(\varepsilon))~.
\end{align*}
which yields \eqref{right-spine-exp} and completes
the proof of Theorem~\ref{spine-exp}.
\rqed

\section{Variance of the Size of the Spine} \label{sec:var}

\newcommand{\A}[2]{#2\underset{#1}\rightsquigarrow\overline #2}
\newcommand{\NA}[2]{#2\underset{#1}{\not\rightsquigarrow}\overline #2}
\newcommand{\B}{\widehat{P}_{n,p}(k;y)}

In this section we shall prove the bounds (\ref{left-spine-var})
and
(\ref{right-spine-var}) on the variance of the size of the spine
given
in Theorem~\ref{spine-var}.  First, we note that the lower bound is
immediate from the Harris-Kleitman correlation inequality (see
\cite{Har60}, \cite{Kle66}), which was later generalized to the FKG
inequality~\cite{FKG71}.
We use $\A{n}{x}$ as a shorthand for $\A{F_{n,p}}{x}$, and if $M$ is a
set of literals, $\A{M}{x}$ means that there is a path from $x$ to
$\overline x$ using only literals in $M$.

\begin{lemma}
\label{iteration}
For strictly distinct literals $x,y$, we have
\begin{equation*}
  \Pr\big( \A{n}{x} , \A{n}{y}\big)
 -\Pr\big( \A{n}{x}\big)^2
  \leq
  \sum_{k\geq 1} P_{n,p}(k)
          \Big( \Pr\big(\A{n}{x}\big)
              - \frac{n-k}{n-1 }\Pr\big(\A{n-k}{x}\big)\Big)~.
\end{equation*}
\end{lemma}
\begin{proof}
\newcommand{\ppynk}{\B}
Let $\ppynk$ denote the event that
$L^{+} (y)$ is strictly distinct, and $|L^{+} (y)|=k$.
So in particular $\Pr\big(\ppynk\big) = P_{n,p}(k)$.
Using the resolution of the identity
$$ 1 = \mathbb I_{y\rightsquigarrow\overline y} + \sum_{k=1}^{\infty}
\mathbb I_{\ppynk},$$
where $\mathbb I_A$ denotes the indicator function of the
event $A$, we can decompose $\Pr\big(\A{n}{x}\big)$ in two different
ways to obtain
\begin{align}
\Pr\big( \A{n}{x} , \A{n}{y}\big)
&  +\!\sum_k \Pr\big(\A{n}{x} , \ppynk \big)
= \Pr\big( \A{n}{x}\big) \Big(\Pr\big(\A{n}{y}\big)+\!\sum_k
\Pr\big(\ppynk\big)\Big),
   \notag\\
\intertext{so that}
\Pr\big( \A{n}{x} , \A{n}{y}\big)
&- \Pr\big( \A{n}{x}\big) \Pr\big(\A{n}{y}\big)
  =\!\sum_k
          \Big( \Pr\big(\A{n}{x}\big)
              -\Pr\big(\A{n}{x} | \ppynk\big)\Big)
          P_{n,p}(k).
   \label{and-decomp}
\end{align}
To estimate the probability on the right, consider the event
$\ppynk$.  With probability
$(k-1)/(n-1)$ either $x$ or $\overline{x}$ is in $L^+(y)$.
If $x\in L^+(y)$
then $x\not\rightsquigarrow\overline{x}$.  If $\overline{x}\in L^+(y)$, the
situation is more complicated, so we shall bound the probability
below by $0$.  If $L^+(y)$ contains neither $x$ nor $\overline{x}$, then any
path from $x$ to $\overline x$ avoids literals in $L^+(y)$.
In this case we may as well
explore the out-graph $L^+(x)$ restricted to avoid the variables in
$L^+(y)$; with probability $\Pr\big(\A{n-k}{x}\big)$ the restricted
out-graph will contain $\overline{x}$.  Thus we conclude
\begin{align*}
 \Pr\big(\A{n}{x} | \ppynk\big) &=
   \frac{n-k}{n-1} \Pr\big(\A{n-k}{x}\big) +
   \frac12 \frac{k-1}{n-1} \Pr\big(\A{n}{x} | \overline x\in L^+(y), \ppynk\big) \\
 &\geq  \frac{n-k}{n-1} \Pr\big(\A{n-k}{x}\big).
\end{align*}
Substituting this into (\ref{and-decomp}) proves the lemma.
\end{proof}

Next we relate the probabilities of the events
$\A{n-k}{x}$ and $\A{n}{x}$.

\begin{lemma}\label{telescope}
If ${\lambda_n} = \varepsilon n^{1/3}$ is large enough, and
$\varepsilon = {\lambda_n} n^{-1/3}$ is small enough, then
$$ \Pr\big( \A{n+1}{x}\big) - \Pr\big( \A{n}{x}\big) \leq
\begin{cases}
\displaystyle\frac{1+\o(1)}{2n|\lambda_n|^3}
  & \text{ if } \lambda_n<0
\\
\\
 \displaystyle\frac{2+\o(1)}{n}
  & \text{ if } \lambda_n>0~.
\end{cases}
$$
\end{lemma}

\begin{proof}
Suppose $\NA{n}{x}$.
Let $X$ denote $L^+_{n,p}(x)$, which then must
be strictly distinct.
We shall consider four cases depending on whether or not
$X\rightarrow x_{n+1}$ and whether or not $X \rightarrow
\overline{x}_{n+1}$:

\noindent{\it Case 1.}
$X\not\rightarrow x_{n+1}$ and $X\not\rightarrow
\overline{x}_{n+1}$.
Then $\NA{n+1}{x}$.

\noindent{\it Case 2.}
$X\rightarrow x_{n+1}$ and $X \rightarrow \overline{x}_{n+1}$.
Then $\A{n+1}{x}$.

\noindent {\it Case 3.}
 $X\rightarrow x_{n+1}$ and $X \not\rightarrow \overline{x}_{n+1}$.
It is clear that if $\A{[n+1]\setminus[X]}{x_{n+1}}$,
where $[X]$ denotes $X\cup\overline{X}$, then $\A{n+1}{x}$.
Suppose that conversely $\A{n+1}{x}$.
Since $\NA{n}{x}$, either $x_{n+1}$
or $\overline{x}_{n+1}$ is in the path from $x$ to $\overline{x}$ in $[n+1]$.
The first such occurence must be $x_{n+1}$, so we have
$x_{n+1}\underset{n+1}\rightsquigarrow\overline{x}$, and its contrapositive
$x\underset{n+1}\rightsquigarrow\overline{x}_{n+1}$.  Since $x\not\rightsquigarrow
\overline{x}_{n+1}$
within $[n]\cup\{\overline{x}_{n+1}\}$, it must be that $x_{n+1}$
occurs in the
path from $x$ to $\overline{x}_{n+1}$.  In particular $\A{n+1}{x_{n+1}}$.
We may assume without loss of generality that $x_{n+1}$ and
$\overline{x}_{n+1}$ occur
only at the endpoints of this path.  If a literal in $X$
occurred in the path
from $x_{n+1}$ to $\overline{x}_{n+1}$, consider the last such one.
The next
literal cannot be $\overline{x}_{n+1}$ by assumption, nor can it be
$x_{n+1}$, since
this occurs only at the beginning of the path.  So the next
literal would
have to lie inside $X$, a contradiction.  If a literal in
$\overline{X}$
occured in the path from $x_{n+1}$ to $\overline{x}_{n+1}$, we
could take the contrapositive
and similarly derive a contradiction.
Thus $x_{n+1} \rightsquigarrow\overline{x}_{n+1}$ using only
literals in $[n+1]\setminus [X]$,
i.e.\ $\A{[n+1]\setminus[X]}{x_{n+1}}$.  Thus conditional
upon $L^{+}_{n}(x)=X$, we can write
\newline
\mbox{\vtop{
\hbox{\hspace{0.45in}$X\rightarrow x_{n+1}$,
       $X\not\rightarrow \overline{x}_{n+1}$,
       $\A{n+1}{x}$
    \ iff \
      $X\rightarrow x_{n+1},\;
      X\not\rightarrow \overline{x}_{n+1},\;
      \A{[n+1]\setminus[X]}{x_{n+1}}$}
\hbox{\hspace{0.3in}$\Pr[\text{$X\rightarrow x_{n+1}$,
           $X\not\rightarrow\overline{x}_{n+1}$,
           $\A{n+1}{x}$\ }|L^{+}_{n}(x)=X] = $ \hfil}
 \hbox{\hspace{2.6in}  $
  \Pr[X\rightarrow x_{n+1}]
  \Pr[X\not\rightarrow \overline{x}_{n+1}]
  \Pr[\A{[n+1]\setminus[X]}{x_{n+1}}]$}
}}
\newline
since the events on the right and the event
$L^{+}_{n}(x)=X$ are determined by pairwise
disjoint sets of variables.

\noindent {\it Case 4.}
 $X\not\rightarrow x_{n+1}$ and $X \rightarrow \overline{x}_{n+1}$.
By symmetry, cases 3 and 4 have the same probability.

\smallskip
Putting these four cases together we see
\begin{align}
  \Pr\big(\A{n+1}{x} |L^+_{n,p}(x) = X \big) =&
   \big[1-(1-p)^{|X|}\big]^{2} +
  2\big[1-(1-p)^{|X|}\big](1-p)^{|X|}
    \Pr [\A{n+1-|X|}{x}] \nonumber\\
 \leq&
   p^{2}|X|^{2} + 2 p|X| \Pr [\A{n+1-|X|}{x}]~. \label{Amx|L+}
\end{align}
As a consequence,
\begin{align}
\Pr\big( \A{n+1}{x}\big) - \Pr\big( \A{n}{x}\big)
&= \Pr\big( \A{n+1}{x} \setminus \A{n}{x}\big)
\notag
\\
&= \sum_{\substack{X \subset [n]
\notag
\\ \text{$X$ strictly distinct}}}
\Pr\big( L^+_{n,p} (x) = X \big)
\Pr\big( \A{n+1}{x}   |L^+_{n,p} (x) = X \big)
\notag
\\
&\leq \sum_{k} P_{n,p}(k) ( p^2 k^2 + 2 p k
\Pr [\A{n+1-k}{x}])
\notag
\\
&\leq p^{2}\sum_{k} P_{n,p}(k) k^{2} + 2 p \Pr [\A{n}{x}]
\sum_{k} P_{n,p}(k) k.
\label{xxx}
\end{align}
Hence by \eqref{moment-a}, we have
\begin{align*}
\Pr\big( \A{n+1}{x}\big) - \Pr\big( \A{n}{x}\big)
&\leq p^2\frac{1+o(1)}{2pn}\frac{1}{|\varepsilon|^3} +
  2 p \Pr [\A{n}{x}] \frac{1+o(1)}{2pn}\frac{1}{|\varepsilon|}
\\
&= \frac{1+o(1)}{4 n^2 |\varepsilon|^3} +
\frac{1+o(1)}{n|\varepsilon|} \times \begin{cases}
   \displaystyle\frac{1+o(1)}{4 n \varepsilon^2} &
\text{\shortstack{if $\lambda_n<0$
 by ~(\ref{left-spine-exp})
}} \\
    (2+o(1))\varepsilon & \text{\shortstack{if $\lambda_n>0$
 by ~(\ref{right-spine-exp})
}}
     \end{cases}
\\
&= [1+o(1)]\times \begin{cases}
 1/[2 n^2 |\varepsilon|^3] & \text{if $\lambda_n<0$,}\\
 2/n & \text{if $\lambda_n>0$.}
\end{cases} \\
&= [1+o(1)]\times \begin{cases}
 1/[2 n |{\lambda_n}|^3] & \text{if $\lambda_n<0$,}\\
 2/n & \text{if $\lambda_n>0$.}
\end{cases}
\end{align*}
In the above, all of the $o(1)$ terms are $\o(1)$.
\end{proof}

\newcommand{\leftright}[2]{\left\{\genfrac{}{}{0pt}{0}{\displaystyle #1}{\displaystyle #2}\right\}}

\begin{corollary}\label{telescoped}
Provided $k \leq n^{2/3}/|{\lambda_n}|$, we have
$$
Pr\big(\A{n}{x}\big) - \Pr\big(\A{n-k}{x}\big) \leq
\begin{cases}
\displaystyle{(1+\o(1))\frac{k}{2n|\lambda_n|^3}}
  & \text{ if } \lambda_n<0
\\
\\
 \displaystyle{(1+\o(1))\frac{2k}{n}}
  & \text{ if } \lambda_n>0~.
\end{cases}
$$
\end{corollary}

\begin{proof}
We seek
$$\Pr\big(\A{n}{x}\big) - \Pr\big(\A{n-k}{x}\big) = \sum_{m=n-k}^{n-1}
\left[\Pr\big(\A{m+1}{x}\big)-\Pr\big(\A{m}{x}\big)\right],$$
and need only show that each summand is well-approximated by
$\Pr\big(\A{n+1}{x}\big)-\Pr\big(\A{n}{x}\big)$.
Define ${\lambda_n}'$ by
\begin{align*}
p\equiv\frac{1+{\lambda_n} n^{-1/3}}{2n} &=
\frac{1+{\lambda_n}' m^{-1/3}}{2m}, \\
\intertext{so that}
\frac{m}{n}\left[1+{\lambda_n} n^{-1/3}\right]-1 &=
{\lambda_n}' m^{-1/3} \\
\intertext{and thus}
\frac{m-n}{n} m^{1/3}+{\lambda_n}
\left(\frac{m}{n}\right)^{4/3} &= {\lambda_n}' ~.
\end{align*}
Since $m=(1+\o(1))n$, and
$|n-m|\leq k \ll n^{2/3}$, it follows that
${\lambda_n}' = (1+\o(1)){\lambda_n}$.
Thus when we apply Lemma~\ref{telescope}
with $m$ and ${\lambda_n}'$ rather than $n$ and ${\lambda_n}$,
we obtain an
answer that differs by a factor of $1+\o(1)$, as desired.
\end{proof}

\begin{proof_of}{Theorem~\ref{spine-var}}
We can now complete our estimate of the covariances.
For strictly
distinct literals $x$ and $y$
we have
\begin{align*}
\operatorname{Cov}&(\A{n}{x},\A{n}{y})
 =
  \Pr\big( \A{n}{x} , \A{n}{y}\big)
 -\Pr\big( \A{n}{x}\big) \Pr\big( \A{n}{y}\big) \\
 &\leq
  \sum_{k\geq 1} P_{n,p}(k)
          \Big( \Pr\big(\A{n}{x}\big)
              - \frac{n-k}{n-1 }\Pr\big(\A{n-k}{x}\big)\Big)
\tag*{by Lemma~\ref{iteration} }
\\
&=
  \sum_{k\geq 1} P_{n,p}(k) \frac{k-1}{n-1}
          \Pr\big(\A{n}{x}\big)
  +
  \sum_{k\geq 1} P_{n,p}(k)
          \frac{n-k}{n-1}
          \Big( \Pr\big(\A{n}{x}\big)
              - \Pr\big(\A{n-k}{x}\big)\Big)
\\
 &\leq
  \Pr\big(\A{n}{x}\big)
  \frac{1}{n}
  \sum_{k\geq 1} P_{n,p}(k) k
  +
  \sum_{k\geq 1}
          P_{n,p}(k)
          \Big( \Pr\big(\A{n}{x}\big)
              - \Pr\big(\A{n-k}{x}\big)\Big) ~.
\end{align*}
As a consequence,
\begin{align*}
\operatorname{Cov}&(\A{n}{x},\A{n}{y})
 \leq
\\
 &\leq
  \Pr\big(\A{n}{x}\big)
  \frac{1}{n}
  \sum_{k\geq 1} P_{n,p}(k) k
  +
  \sum_{k\geq 1}
          P_{n,p}(k)
          \Big( \Pr\big(\A{n}{x}\big)
              - \Pr\big(\A{n-k}{x}\big)\Big) \\
 &=
  \Pr\big(\A{n}{x}\big)
  \frac{1+o(1)}{n\varepsilon}
  +\left[
  \sum_{k\leq n^{2/3}/|{\lambda_n}|}+
  \sum_{k> n^{2/3}/|{\lambda_n}|}\right]
          P_{n,p}(k)
          \Big( \Pr\big(\A{n}{x}\big)
              - \Pr\big(\A{n-k}{x}\big)\Big) \\
 &\leq
  \Pr\big(\A{n}{x}\big)
  \frac{1+o(1)}{n\varepsilon}
  +
  \sum_{k\leq n^{2/3}/|{\lambda_n}|}\!\!\!
          P_{n,p}(k) \leftright{\frac{(1+o(1))k}
{2n|{\lambda_n}|^3}}{\frac{(2+o(1))k}{n}}
  +
  \sum_{k> n^{2/3}/|{\lambda_n}|}\!\!\!  P_{n,p}(k) \Pr\big(\A{n}{x}\big)
  \\
 &=
  \leftright{\frac{1+o(1)}{4\varepsilon^2 n}}{(2+o(1))\varepsilon}
  \frac{1+o(1)}{n\varepsilon}
  +
  \leftright{\frac{1+o(1)}{2n|{\lambda_n}|^3\varepsilon}}
{\frac{2+o(1)}{n\varepsilon}}
  +
  \leftright{\exp[-\Theta({\lambda_n})] n^{-1/3}
             \times
             \frac{1+o(1)}{4\varepsilon^2 n}}
            {\exp[-\Theta({\lambda_n}^{3/5})] n^{-1/3}
            \times
            (2+o(1))\varepsilon}
  \\
 &=
  \leftright{\frac{1+o(1)}{2n|{\lambda_n}|^3\varepsilon}}
{\frac{2+o(1)}{n\varepsilon}}
  =
  \leftright{\frac{1+o(1)}{2n^{2/3}{\lambda_n}^4}}{\frac{2+o(1)}
{n^{2/3}{\lambda_n}}}
\end{align*}
where the upper entry applies to the subcritical regime, and the
lower entry applies to the supercritical regime.
In the above, the $o(1)$ terms are all $\o(1)$.
This gives the variance bound above threshold.
The bound on the second moment below threshold
follows by combining the
above variance bound with Theorem~\ref{spine-exp}.
\end{proof_of}

\section{The Existence of Hourglasses} \label{sec:hourglass}

In this section we prove Theorem~\ref{thm:hourglasses}, which was
central to the proof of the upper bounds in Theorem~\ref{sat_{n,p}}.
\medskip

\begin{proof_left}{Proof of  Theorem~\ref{thm:hourglasses} (i)
on the
existence of many hourglasses}

Let $t$ be a large positive number, but still small compared to
$n^{1/3}$, and let $p=(1-t n^{-1/3})/2n$.
We shall construct a process for growing
hourglasses one by one.
However, as we deplete variables, the distribution changes,
so that the hourglasses so constructed are not identically
distributed.  We therefore use the following two variations
of this naive procedure:  In the first process, instead
of drawing from $n$ variables, we draw only from $n'=n-tn^{2/3}$
variables, so that we have a ``buffer'' of $tn^{2/3}$,
and replenish the variables as necessary.  However, even
this process can lead to trouble in the unlikely event
that we have a very long run which uses up our entire
reserve of variables.  So we construct another process
which aborts the growth of an hourglass when
it would use up too many variables.

To be explicit, consider the trimmed
out-graphs and trimmed in-graphs of various literals, where the
in- and out-graphs
are restricted to sets of $n' = n-t n^{2/3}$ variables.
Recall from Lemma~\ref{trimmed-graph} that the unoriented projection
of the trimmed out-graph $\widetilde D_{F_{n',p}}^+(x)$ of a
literal $x$ is identically distributed to the connected component
$C_{n',2p-p^2}(x)$ in $G_{n',2p-p^2}$.
We shall follow the
convention that no matter how many trimmed out-graphs or trimmed
in-graphs we have explored in the past, when exploring another
in-graph or out-graph, we shall always restrict our attention to
variables that are in none of the in-graphs or out-graphs found so
far (except, possibly, for the root), and we shall add
variables to ensure that there are $n' = n-t
n^{2/3}$ variables for the tree to grow within.  (Recall that for the
upper bounds we assume that we have a variable for each natural
number.) In this way the sizes and structures of all the trees will be
independent and identically distributed.  As we alluded to above,
later we shall consider a
variation on this process.

Since there are somewhat fewer than $n$ variables in which we explore
the out-graphs and in-graphs, this decreases the average out-degree of
each literal, and has the effect of shifting the formula
further into the
subcritical regime.  Specifically, if we define $t'$ by
$$ \frac{1-t n^{-1/3}}{2n} = p = \frac{1-t' (n')^{-1/3}}{2n'}, $$
we see that $t' = (2+o(1)) t$.

Pick a literal $u$, and look at its trimmed out-graph $T$ within
$n'$ unused variables.
Recall the definition of $R_{n,p}(k)$, the probability that
$C_{n,2p-p^2}(u)$ is a tree of size $k$.
By using Lemma~\ref{lem:pnk} for $R_{n,p}(k)$,
we see that, for some $c$, there is
a $(c+o(1))(t/n^{1/3})$ chance that $T$ is a tree of size
between $2 n^{2/3}/t^2$ and $4 n^{2/3}/t^2$.
Here, as explained in the first sentence of the proof, we are
assuming $1 \ll t \ll n^{1/3}$, so by Remark \ref{much_less_than},
in the remainder of this proof all $o(1)$ terms without subscripts
are of the form $o_{t/n^{1/3},t}(1)$.
By comparing with random
graphs, if $T$ is a tree, then it is uniformly distributed amongst the
spanning trees on $|T|$ vertices.
Using the structural properties of random spanning trees (see e.g.\
\cite{Aldous90}),
if we pick a random vertex $w\neq u$ in  $T$, with
probability at least $(1-o_{|T|}(1))e^{-1}$ the path connecting $u$ to
$w$ has length at least $\sqrt{2|T|}$.  Let $v$ be the middle vertex in
the path from $u$ to $w$ (in case of tie, we choose the vertex
closer to $w$).  Either the majority of the remaining vertices are connected
to the first part of the path or to the second part of the path.  Since the
spanning tree is uniformly random, by symmetry,
with probability at least $1/2$, at least half of the remaining
vertices of the tree will be connected to $v$ via the path from $v$ to
$w$ --- these vertices will be in the out-graph of $v$.  In the event
that
the path from $u$ to $v$ has length at least
$({|T|/2})^{1/2} \geq n^{1/3}/t$, and
$v$ has at least $\frac12|T| \geq  n^{2/3}/t^2$
descendents in the trimmed
out-graph of $u$, say that vertex $v$ is ``promising,'' and that the
path from $u$ to $v$ is the tail.  We thus have
shown
$$\Pr\left[\text{a random vertex
$v$ in $\widetilde D^+_{F_{n',p}}(u)$ is promising}\right]
 = (1+o(1))\frac{c}{2 e} \frac{t}{n^{1/3}}.$$

In the event that
 vertex $v$ is promising,
we proceed to explore the trimmed in-graphs of the first $n^{1/3}/t$
vertices on $v$'s tail.  Again by Lemma~\ref{lem:pnk}, each
individual in-graph will have size at
least $n^{2/3}/t^2$
with probability at least $(c+o(1))(t/n^{1/3})$.
The probability that none of them is so large is at most
$(1+o(1))e^{-c}$.  If any of them is so large,
we shall call $v$ a
``central variable,'' and $v$ together with its explored out-graph and
in-graph both of size at least $n^{2/3}/t^2$ an ``hourglass''
(see Definition~\ref{hourglass}).
Each time that we pick a literal and look for an hourglass as
described above, we find one with probability at least
$$
(1+o(1))\frac{c(1-e^{-c})}{2 e} \frac{t}{n^{1/3}}.
$$

Next let us compute how many variables we expect to use up while
exploring the trimmed out-graph and trimmed in-graphs.
At this point we recall from
Lemma~\ref{RC-moments} and equation \eqref{RG-exp}
that if $G$ is either the trimmed in-graph or
trimmed out-graph of a vertex,
$$
E[|G|] =
(1+o(1))\frac{(n')^{1/3}}{t'} = (1+o(1))\frac{n^{1/3}}{2t}.
$$
We always explore one out-graph, and with probability
$(1+o(1))[c/(2 e)] t/n^{1/3}$ we explore $\lceil n^{1/3}/t\rceil$
in-graphs.  Thus the expected number of variables used up is
$$
(1+o(1))\left[1 + \frac{c}{2 e}\right]\frac{n^{1/3}}{2t}.
$$

If we look for an hourglass for
$$ \frac{4 e}{c(1-e^{-c})} n^{1/3}/t
$$
times, then the
probability that we fail to find an hourglass is $(1+o(1)) e^{-2}$,
and the expected number of variables that we use up is
$$
(1+o(1)) \left[\frac{2 e}{c} +1\right] \frac{1}{1-e^{-c}}
\frac{n^{2/3}}{t^2}.
$$
The probability that we use up more than $3(1+o(1))$
times the expected number of variables is at most
$1/3+o(1)$.  Therefore, with
probability at least $1-e^{-2} - 1/3 -o(1) \ge 1/2$
(for large enough
$t$ and small enough $t n^{-1/3}$)
we both find an hourglass, and do not use up
more than $b n^{2/3}/t^2$ variables,
where $b = 3 [1+2e/c]/(1-e^{-c})$.

Now consider the following modification of the above procedure: as
above we use the local search procedure in Section 2 to explore the
trimmed out-graphs and trimmed in-graphs, hoping to find an hourglass,
but as soon as we use up
$b n^{2/3}/t^2$ variables, we abort and stop looking for the
hourglass.  Then we can repeat this procedure $t^3/b$ times, be
guaranteed to use no variables other than the first $n$ of them, and
find a number of disjoint hourglasses that stochastically dominates
the binomial distribution $\Bin(t^3/b,1/2)$.
By Chernoff's inequality \cite{Che52} (see also \cite{McDiarmid89}), the
probability that we find fewer than half as many hourglasses as we
expect will be no larger than $\exp[-t^3/(8b)]$.
\end{proof_left}

\begin{proof_right}{Proof of  Theorem~\ref{thm:hourglasses} (ii)
on the existence of a giant hourglass}
Let $t$ be a large positive number, but still small compared to
$n^{1/3}$.  When $p=(1-t n^{-1/3})/2n$, as we have just seen,
there will be $\Theta(t^3)$ hourglasses with in- and out-portion
of size at least $n^{2/3}/t^2$,
except with probability $\exp(-\Theta(t^3))$.
We now increase $p$ by
a suitably large constant times $t n^{-4/3}$,
say $Mt n^{-4/3}$.
For any two hourglasses, the probability of an edge
from the out-portion of the first hourglass to the
in-portion of the
second hourglass is therefore at least
$M/t^3$.  Then the central
variable of the first hourglass implies the central variable of the
second hourglass.  Conceptually we can think of the directed graph
whose nodes are the hourglasses, and place a directed edge from one
node to another whenever the hourglasses connect up like this.  The
edges of this graph occur independently of one another, and the
average
out-degree of the graph is $\Theta(Mt^3/t^3)=\Theta(M)$.  By
choosing $M$ large enough, we can make the
average out-degree to be any convenient constant that we like.  In
particular, if the average out-degree is a constant larger than 1,
then we might expect the connections to percolate, so that there is
some node $v$ that can reach a constant fraction of the other nodes
through edges of this graph, and is reachable by a constant fraction
of the other nodes.  Provided this happens, each literal in the
out-portions of the nodes reachable by $v$ is implied by the central
variable of node $v$, and each literal in the
in-portions of the nodes
that can reach $v$ will imply the central variable of $v$, thereby
giving the desired giant hourglass with in-portion and
out-portion each
of
size $\Theta(t n^{2/3})$.

It remains to be shown
(in Lemma~\ref{lem:super-high-perc}) that we get the
requisite percolation except with probability that is exponentially
small in the number of nodes of the graph.
\end{proof_right}

The following lemma is related to one proved by Karp in \cite{Kar90}
which showed that with high probability there is a giant component
of size $\Theta(N)$ in supercritical directed percolation.  For our
purposes ``with high probability'' is not sufficient; we need the
exceptional events to be exponentially rare.

\begin{lemma}
\label{lem:super-high-perc}
In a random directed graph on $N$ vertices, in which each directed
edge occurs independently with probability $6/N$, then except with
probability that is exponentially small in $\Theta(N)$, there is a
vertex $v$ with out-graph of size $\Theta(N)$ and in-graph of size
$\Theta(N)$.
\end{lemma}

\begin{proof}
For convenience, let $N' = \lceil N/3\rceil$,
so that there are at least $3 N' - 2$ vertices,
and the probability of each directed edge
is at least $2/N'$.  (We can throw out some of the edges
to make the probability exactly $2/N'$.)
To search for a node $v$ with a large in-graph and out-graph,
we consider candidate vertices one at a time, and explore the in-graph
and out-graph of the candidate, restricting the explorations of the
in-graph and out-graph to disjoint sets of $N'-1$ vertices, none of
which have yet been explored in the course of examining a previous
failed candidate.  In this way we ensure that the sizes of the
in-graph and out-graph are independent, and both distributed in the
same manner as the size of the component containing a particular
vertex in the random graph $G_{N',2/N'}$.  We do the explorations in a
parallel interleaved fashion, so that if either the in-graph or
out-graph is found to be too small, then exploration of the other is
immediately halted.  We shall show that
for large enough $N'$, except with probability
exponentially small in $N'$, after looking at $N'/50$ candidates,
the failed candidates have not wasted more than $N'$ variables,
and we find a successful candidate with in-graph and
out-graph each of size at least $N'/10$.

We first claim that
for any real $\gamma$ and integers $k$ and $N'$ such that
$0\leq\gamma\leq N'$ and
$1\leq k\leq N'$, there is an $s$ between $1/2$ and $1$ so that
\begin{equation}
\label{claim}
\Pr\big(|C(x)|=k\text{ in }G_{{N'},\gamma/{N'}}\big) =
  k^{k-2}\binom{{N'-1}}{k-1}\left(\frac{\gamma}{N'}\right)^{k-1}
  \left(1-\frac{\gamma}{N'}\right)^{k {N'} - s k^2}.
\end{equation}

Indeed,
the probability that the component $C(x)$ containing the
vertex $x$ has size $k$
can be bounded below by the probability that $C(x)$ is a tree
of size $k$.  To bound $\Pr\big(|C(x)|=k\big)$  from above,
we note that the connectedness of $C(x)$ implies that
$C(x)$ contains a tree of size $k$.  Summing over all
possibilities for this spanning tree, we therefore get
$$
 \sum_{\substack{\text{trees $T\ni x$}\\|T|=k}}
\Pr\big(T = C(x)\big)
\leq
\Pr\big(|C(x)|=k\big)
\leq
 \sum_{\substack{\text{trees $T\ni x$}\\|T|=k}}
\Pr\big(T\subseteq C(x), |T|=|C(x)|\big).$$
This give a lower bound of
$$
 k^{k-2}\binom{{N'-1}}{k-1}\left(\frac{\gamma}{N'}\right)^{k-1}
  \left(1-\frac{\gamma}{N'}\right)^{k (N'-k)
                                 + \left(\binom k2-(k-1)\right)}
$$
and an  upper bound of
$$
 k^{k-2}\binom{{N'-1}}{k-1}\left(\frac{\gamma}{N'}\right)^{k-1}
  \left(1-\frac{\gamma}{N'}\right)^{k (N'-k)}~,
$$
which establishes the claim~\eqref{claim}.

Next define $X$ by
$$ X = \begin{cases} |C(x)| & \text{if $|C(x)| < N'/10$} \\ 0 &
  \text{if $|C(x)| \geq N'/10$}\end{cases}.$$
Using the bound \eqref{claim} on $\Pr\big(|C(x)|=k\big)$ we see that
\begin{align}
\Pr\big(|C(x)|=k\big) &\leq
  k^{k-2}\binom{{N'-1}}{k-1}\left(\frac{\gamma}{N'}\right)^{k-1}
  \left(1-\frac{\gamma}{N'}\right)^{k {N'} - k^2}\notag\\
&\leq \frac{k^{k-1}}{k!}  \gamma^{k-1}
      e^{-k\gamma + \gamma k^2/{N'}}. \label{pcx=k}\\
\intertext{For $k\leq N'/10$, we get}
\Pr\big(|C(x)|=k\big)
&\leq \frac{k^{k-1}e^{-k}}{k!}  \frac{1}{\gamma}
      e^{k(1-\gamma + \log \gamma + \gamma/10)}. \notag\\
\intertext{Assuming $\gamma = 2$, this gives}
\Pr\big(|C(x)|=k\big)
 &\leq \frac{k^{k-1}e^{-k}}{k!}  \frac{1}{2}
      e^{-k/10} \label{pcx=k,restricted}\\
\intertext{so that}
\sum_{k=1}^{N'/10} \Pr\big(|C(x)|=k\big) e^{k/10}
 &\leq \sum_{k=1}^{N'/10} \frac{k^{k-1}e^{-k}}{k!}  \frac{1}{2}
  \leq \frac12 \sum_{k=1}^{\infty} \frac{k^{k-1}e^{-k}}{k!} =
  \frac12. \notag\\
\intertext{Consequently}
E[e^{X/10}] &\leq \frac12 + \Pr\big(|C(x)|\geq N'/10\big) \times e^{0/10}
 \leq 3/2.
\notag
\end{align}

Let $Y_i$ be the number of variables lost on the $i$th candidate if it
is a failure; $Y_i=0$ if the $i$th candidate is successful:
$$ Y_i = \begin{cases} 0 & \text{$i$th candidate successful} \\
                       2k-1 & \text{$i$th candidate failed because in-graph had
                                    size $k<N'/10$} \\
                       2k   & \text{$i$th candidate failed because out-graph had
                                    size $k<N'/10$.}
         \end{cases}$$
Thus $\Pr(Y=j) \leq \Pr\big(X=\lceil j/2\rceil\big)$, and hence
$$ E[e^{Y_i/20}] = \sum_j \Pr(Y_i=j) e^{j/20} \leq
   2 \sum_{\text{$j$ even}} \Pr(X=j/2) e^{j/20} = 2 E[e^{X/10}] \leq 3.$$
Letting $$S = \sum_{i=1}^{\beta N'}
Y_i$$ be the total number of variables lost on failed candidates, we see
\begin{align*}
\Pr(S > N')
&=
  \Pr\big(e^{S/20} > e^{N'/20}\big) \\
&\leq
  \frac{E[e^{S/20}]}{e^{N'/20}} \\
&=
  \frac{E[e^{Y/20}]^{\beta N'}}{e^{N'/20}} \\
&\leq
  \frac{3^{\beta N'}}{e^{N'/20}},
\end{align*}
which is exponentially small in $N'$ for $\beta=1/50$.

Next consider $\Pr\big(|C(x)| \geq N'/10\big)$.
By \eqref{pcx=k} and \eqref{pcx=k,restricted} we have
\begin{align*}
\Pr\!\left(|C(x)| \geq \frac{N'}{10}\right) &=
 1 - \sum_{k=1}^{10\log N'} \Pr\big(|C(x)|=k\big)
   - \sum_{k=10\log N'}^{N'/10} \Pr\big(|C(x)|=k\big) \\
&\geq
 1 - \exp[200\log^2 N' / N'] \!\sum_{k=1}^{10\log N'}\!\!
                               \frac{k^{k-1}}{k!} \gamma^{k-1} e^{-k\gamma}
   - \!\!\!\sum_{k=10\log N'}^{N'/10} \!\!\! \frac{k^{k-1}e^{-k}}{k!}  \frac{1}{2}
   e^{-k/10} \\
&>
 1 - \exp[200\log^2 N' / N'] \sum_{k=1}^{\infty}
                               \frac{k^{k-1}}{k!} \gamma^{k-1} e^{-k\gamma}
   - \frac1{2 N'} \\
&=
 1 - \exp[200\log^2 N' / N'] (1-\vartheta(\gamma-1))
   - \frac1{2 N'} \\
&= \vartheta(\gamma-1) - O(\log^2 N' / N').
\end{align*}
Thus with probability $(1-o(1))\vartheta(\gamma-1)^2 > 0.63$ (for
large enough $N'$), both the
in-graph and out-graph have size at least $N'/10$.

Therefore, if we try $N'/50$ candidates, the probability that we lose too
many variables on failed candidates, or have enough variables but
still fail to find a vertex with in-graph and out-graph of size at
least $N'/10$, is bounded by
$$ \frac{3^{N'/50}}{e^{N'/20}} +  0.37 ^{N'/50},$$
which establishes the lemma.
\end{proof}

\appendix

\section{Relation between $F_{n,p}$ and $F_{n,m}$}

While the literature focuses on $F_{n,m}$, where a given
number of clauses are specified, most of our theorems and
proofs are done for $F_{n,p}$, where each clause has some
independent chance of appearing in the formula, and the total number
of clauses is random.
In this appendix, we discuss the
relation between the models $F_{n,m}$ and
$F_{n,p}$.
With respect to monotone properties, these two models
are practically interchangeable, provided $m$ is about
$4\binom{n}{2}p$, the expected number of clauses in $F_{n,p}$.
Write $N$ for the number of 2-clauses on
$x_1, \dots , x_n$, so that $N=2n(n-1)$, and let $M_{N,p}$ be
a binomial random variable with parameters $N$ and $p$. Then
we have
\begin{align*}
\Pr(\SAT(F_{n,p})) &= \sum_{m=0}^N \Pr(M_{N,p}=m) \Pr(\SAT(F_{n,m})) ~, \\
E(S(F_{n,p})) &= \sum_{m=0}^N \Pr(M_{N,p}=m) E(S(F_{n,m})) ~,\\
\intertext{and}
E(S^2(F_{n,p})) &= \sum_{m=0}^N \Pr(M_{N,p}=m) E(S^2(F_{n,m})).
\end{align*}

Since $\Pr (\SAT(F_{n,m}))$ is a monotone
decreasing function of $m$, for every $0 < m < N$ we have
\begin{equation*}\label{satbound}
\Pr (\SAT(F_{n,m})) -\Pr (M_{N,p} > m)
 \le \Pr (\SAT(F_{n,p}))
 \le \Pr (\SAT(F_{n,m})) +
                \Pr (M_{N,p} < m).
\end{equation*}
Similarly, we have
\begin{align*}\label{unsatbound1}
\Pr(\UNSAT(F_{n,m})) -\Pr(M_{N,p} < m)&
\\
\le \Pr(\UNSAT&(F_{n,p}))
\\
&\le \Pr(\UNSAT(F_{n,m})) +
     \Pr(M_{N,p} > \! m)~,
\end{align*}
\begin{equation*}\label{unsatbound2}
     E(S(F_{n,m})) - 2n\Pr(M_{N,p} < m)
 \le E(S(F_{n,p}))
 \le E(S(F_{n,m})) + 2n\Pr(M_{N,p} > \! m)~,
\end{equation*}
and
\begin{equation*}\label{unsatbound3}
     E(S^2(F_{n,m})) - 4n^2\Pr(M_{N,p} < m)
 \le E(S^2(F_{n,p}))
 \le E(S^2(F_{n,m})) + 4n^2\Pr(M_{N,p} > \! m)~.
\end{equation*}

In bounding the probability in the tail of the binomial
distribution, we shall make use of the following
Chernoff type inequality
(see e.g.\ \cite{McDiarmid89}):
\begin{equation}\label{tailprob}
\Pr ( |M_{N,p}-pN|
   \ge \rho pN )
       \le  e^{-\rho^2pN/3},
\end{equation}
provided $0<p\le 1/2$.

To bound the probability of unsatisfiability in the subcritical
regime, the expected size of the spine, and the second moment of
the size of the spine, we set
$$ p = \frac{1+\lambda n^{-1/3}}{2n}$$
and
$$ m = (1+\lambda' n^{-1/3}) n,$$
with $$\lambda' = \lambda \pm n^{-1/12}.$$
Then the probability that $M_{N,p}$ is too large or too small
(compared with $m$) is $$\exp(-\Theta(n^{1/6})) = o(1/n^2).$$
From this we
see that our bounds for $F_{n,p}$ imply the desired bounds for
$F_{n,m}$ in the subcritical regime.

To convert the bounds on the probability of satisfiability on the
right from $F_{n,p}$ to $F_{n,m}$, because this probability is so
small, in order for it to dwarf $\Pr(M_{N,p}<m)$ and $\Pr(M_{N,p}>m)$,
we need to have a larger gap between $\lambda$ and $\lambda'$.
We set $$\lambda' = \lambda \pm \frac{\lambda^{6/5}}{n^{1/15}}.$$
Then the probability of a large deviation is at most
$$\exp[-\Theta(\lambda^{12/5}n^{1/5})] =o\big(
\exp[-\Theta(\lambda^3)]\big),$$
if $\lambda$ is small enough compared to $n^{1/3}$.
Furthermore, $\lambda'/\lambda$ is arbitrarily close to $1$
provided $\lambda$ is
sufficiently small compared to $n^{1/3}$.  Thus our bounds for
satisfiability of $F_{n,p}$ in the supercritical regime carry over
to corresponding bounds for $F_{n,m}$.

\section{Proof of Lemma~\ref{weakgiant}}
\label{sec:appendix}

In this appendix, we prove Lemma~\ref{weakgiant}.
We start with the proof of statement (i).  To
this end, we again use
that the cluster size distribution
in $G_{n,\widetilde p}$, $\widetilde p=2p-p^2$, is stochastically
dominated by a Poisson birth process
with  parameter
\begin{equation}
\label{b.1}
\hat\kappa=n \log 1/({1-\widetilde p})=2n\log 1/(1-p) = 2np(1+O(1/n)).
\end{equation}
Writing this parameter as
$\hat\kappa=1+\hat\varepsilon$,
and observing that
the probability that a Poisson birth tree with parameter
$\hat\kappa$ has size $k$ is
$(1/\hat\kappa)
(k^{k-1}/k!)
({\hat\kappa}e^{-{\hat\kappa}})^k$,
we therefore get
\begin{align}
\sum_{k\ge   n^{2/3}/\lambda_n }
Q_{n,p}(k)
&\leq
\vartheta(\hat\varepsilon)
+\frac 1{\hat\kappa}\sum_{k\ge   n^{2/3}/\lambda_n }
\frac{k^{k-1}}{k!}
({\hat\kappa}e^{-{\hat\kappa}})^k
\notag
\\
&\leq
\vartheta(\hat\varepsilon)
+\frac 1{\hat\kappa}\sum_{k\ge   n^{2/3}/\lambda_n }
\frac{1}{\sqrt{2\pi k^3}}
({\hat\kappa}e^{1-{\hat\kappa}})^k~
\label{b.3}
\end{align}
where
\begin{equation}
\label{b.2}
\vartheta(\hat\varepsilon)
=
1-\frac 1{\hat\kappa}\sum_{k=0}^\infty
\frac{k^{k-1}}{k!}
({\hat\kappa}e^{-{\hat\kappa}})^k
\end{equation}
is the survival
probability in a Poisson birth process
with parameter
$\hat\kappa$.
If $\varepsilon\leq\varepsilon_0$, then
\begin{equation}
\label{b.4}
{\hat\kappa}e^{1-{\hat\kappa}}
\leq e^{-c\varepsilon^2}
\end{equation}
for some constant $c=c(\varepsilon_0)$, so that
\begin{align}
\sum_{k\ge  n^{2/3}/\lambda_n  }
Q_{n,p}(k)
\leq
\vartheta(\hat\varepsilon)
+O\Big(\frac 1{\sqrt{\lambda_n  n^{2/3}}}
e^{-c\lambda_n}
\Big)
\leq
\vartheta(\hat\varepsilon)
(1+O(e^{-c\lambda_n}))
\label{b.5}
\end{align}
where we have used that
$\vartheta(\hat\varepsilon)=
\Theta(\hat\varepsilon)
=\Theta(\lambda_n n^{-1/3})$
in the last step.
Since
$\hat\varepsilon=\varepsilon+O(n^{-1})
=\varepsilon(1+O(\lambda_n^{-3}))$,
we conclude that
\begin{align}
\sum_{k\ge  n^{2/3}/\lambda_n  }
Q_{n,p}(k)
&\leq
\vartheta(\varepsilon)
(1+O(\lambda_n^{-3})) \label{b.5a}\\
&=\Theta(\lambda_n n^{-1/3})~. \label{b.5b}
\end{align}

To prove a lower bound,
we show that
\begin{align}
\sum_{k<  n^{2/3}/\lambda_n  }
Q_{n,p}(k)
\leq
1-\vartheta(\varepsilon)
(1+O(\lambda_n^{-2})).
\label{b.6}
\end{align}
Together with
(ii), which will be proved below,
this gives
the desired bound.
To prove \eqref{b.6}, we will first show
that for $k\leq n^{2/3}/\lambda_n$,
\begin{align}
Q_{n,p}(k)
\leq
\frac 1{1+\varepsilon'}
\frac{k^{k-1}}{k!}
((1+\varepsilon')e^{-(1+\varepsilon')})^k
\Big[1+O\Big(\frac{k^{3/2}}{n}\Big)\Big],
\label{b.11}
\end{align}
where
$\varepsilon'$ is defined as the positive solution of
\begin{equation}
(1+\varepsilon')
e^{-\varepsilon'}
=(1+\varepsilon)e^{-\varepsilon}
e^{\varepsilon^2\lambda^{-2}/2}
.
\label{b.9}
\end{equation}
Indeed, using \eqref{Qp}, we bound
\begin{align}
Q_{n,p}(k)
&\leq \frac{1}{n} \binom{n}{k} (2pk)^{k-1}
e^{-p({2nk-k^2-3k+2})} S_{2p-p^2}(k)
\notag
\\
&=\frac{(2pnk)^{k-1}}{k!}
\Big[\prod_{i=0}^{k-1}\Big(1-\frac in\Big)\Big]
e^{-p({2nk-k^2-3k+2})} S_{2p-p^2}(k)
\notag
\\
&\leq \frac{(2pnk)^{k-1}}{k!}e^{-k(k-1)/2n}
e^{-p({2nk-k^2-3k+2})} S_{2p-p^2}(k)
.
\label{b.7}
\end{align}
For $k\leq n^{2/3}/\lambda_n$,
we have $S_{2p-p^2}(k)=1+O(k^{3/2}/n)$.
Combined with the observation that
$pk^2-k^2/2n=\varepsilon k^2/2n
\leq k\varepsilon^2/2\lambda^2
$
we therefore get
\begin{align}
Q_{n,p}(k)
&\leq \frac{(2pnk)^{k-1}}{k!}
e^{-(2pn-\varepsilon^2\lambda^{-2}/2)k}
\Big[1+O\Big(\frac{k^{3/2}}{n}\Big)\Big]
\notag
\\
&=
\frac 1{1+\varepsilon}
\frac{k^{k-1}}{k!}
\left((1+\varepsilon)
e^{-(1+\varepsilon-\varepsilon^2\lambda^{-2}/2)}
\right)^k
\Big[1+O\Big(\frac{k^{3/2}}{n}\Big)\Big]
.
\label{b.8}
\end{align}
Using the definition \eqref{b.9} of
$\varepsilon'$ and observing
that
 $\varepsilon'\leq \varepsilon$,
we get \eqref{b.11}.

As a consequence of  \eqref{b.11},
we now have
\begin{align}
\sum_{k<  n^{2/3}/\lambda_n  }
Q_{n,p}(k)
&\le
\sum_{k<   n^{2/3}/\lambda_n }
\frac 1{1+\varepsilon'}
\frac{k^{k-1}}{k!}
((1+\varepsilon')e^{-(1+\varepsilon')})^k
\Big[1+O\Big(\frac{k^{3/2}}{n}\Big)\Big],
\\
&\leq
1-\vartheta(\varepsilon')
+
O\Big(\sum_{k<   n^{2/3}/\lambda_n }
\frac{k^{3/2}}{n}
\frac 1{1+\varepsilon'}
\frac{k^{k-1}}{k!}
((1+\varepsilon')e^{-(1+\varepsilon')})^k
\Big)
.
\label{b.12}
\end{align}
Observing that $\varepsilon\leq\varepsilon_0$
implies $\varepsilon'\leq\varepsilon_0$,
which implies a bound of the form
\eqref{b.4} for $\kappa'=1+\varepsilon'$,
we now bound the sum over $k$ as
follows:
\begin{align}
\sum_{k< n^{2/3}/\lambda_n}
\frac{k^{3/2}}{n}
\frac 1{\kappa'}
&\frac{k^{k-1}}{k!}
(\kappa'e^{-\kappa'})^k
\leq
\sum_{k<n^{2/3}/\lambda_n}
\frac{k^{3/2}}{n}
\frac{1}{\sqrt{2\pi k^3}}
(\kappa'e^{1-\kappa'})^k
\notag
\\
&=
O\Big(\frac{1}{n}\Big)
\sum_{k< n^{2/3}/\lambda_n}
e^{-c\varepsilon^2 k}
=
O\Big(\frac{1}{n}\Big)\frac 1{1-e^{-c\varepsilon^2}}
\notag
\\
&=O\Big(\frac 1{n\varepsilon^2}\Big)
=\vartheta(\varepsilon)O(\lambda_n^{-3}).
\label{b.13}
\end{align}
Together with the observation that
$\varepsilon'=\varepsilon(1-\Theta(\lambda_n^{-2}))$
for $\varepsilon$ small enough,
the bounds \eqref{b.12} and \eqref{b.13} imply
\eqref{b.6}.

\newcommand{\p}{\widetilde{p}}
Next we prove statement (ii).  To this end,
we use a refinement (due to Alon and Spencer \cite{alon-spencer})
of the viewpoint employed by Karp \cite{Kar90} and used here in the
proof of Lemma~\ref{basic1}.
Define $N_0=n-1$, and for positive $t$,
$N_t=\Bin(N_{t-1},1-\p)$.  Then let $Y_t = n-t-N_t$, and define $T$ to
be the least $t$ such that $Y_t=0$.  This random variable $T$ has the
same distribution as the size of the connected component containing a
given vertex in $G_{n,\p}$ \cite{alon-spencer}.

Condition the $N_t$ process to be small enough that $Y_t$ is positive
whenever $t<n^{2/3}/\lambda$ (i.e.\ $T\geq n^{2/3}/\lambda$).  How
does this affect the distribution of $N_t$ for larger values of $t$?
We can think of the $N_t$'s as being determined by a collection of
i.i.d.\ 0-1 random variables with probability $1-\p$.  Each $N_t$ is
monotone increasing in these variables.  By FKG, the above conditioning can
only decrease the distribution of $N_t$ (increase the distribution of
$Y_t$) for any given value of $t$, and in particular makes it less
likely that $Y_t\leq 0$ for some $t$ in a given range.  Thus we have
\begin{align}
 \Pr\Big(n^{2/3}/\lambda \leq T \leq n^{2/3} \lambda\Big)\!
   &= \Pr\Big( n^{2/3}/\lambda \leq T \Big)
      \Pr\Big(\exists t\!: Y_t \leq 0,\ n^{2/3}/\lambda \leq t
\leq \lambda n^{2/3} | n^{2/3}/\lambda \leq T\Big) \notag\\
  &= O\Big(\frac{\lambda}{n^{1/3}}\Big)
      \Pr\Big(\exists t\!: Y_t \leq 0,\ n^{2/3}/\lambda \leq t
\leq \lambda n^{2/3} \Big) \tag*{by \eqref{b.5b}\quad\quad}\\
  &= O\Big(\frac{\lambda}{n^{1/3}}\Big)\!
      \left[ \Pr\Big(Y_{n^{2/3}/\lambda}\!\leq\! 0 \Big) +
             \Pr\Big(\exists t\!: Y_t \!=\! 0,\ n^{2/3}/\lambda \!\leq t \leq\! \lambda n^{2/3} \Big) \right]\!,
\label{n**-1/3*}
\end{align}
where in the last line we used $Y_{t+1}\geq Y_t-1$.

Let $X_t$ denote the event that $Y_t=0$.  Let $Z_t$ denote the event that $Y_t=0$ and $Y_s>0$ for $n^{2/3}/\lambda \leq s < t$.  Let $$ S = \sum_{n^{2/3}/\lambda \leq t \leq \lambda n^{2/3} + n^{2/3}/\lambda^2} \mathbb I_{X_t},$$
where as before $\mathbb I_A$ denotes the indicator of the event $A$.
We have
\newcommand{\sumt}{\sum_{n^{2/3}/\lambda \leq t \leq \lambda n^{2/3}}}
\newcommand{\sumd}{\sum_{0\leq\Delta\leq n^{2/3}/\lambda^2}}
\begin{align}
  S &\geq \sumt \mathbb I_{Z_t} \sumd \mathbb I_{X_{t+\Delta}}~, \notag\\
\intertext{so that}
E[S] &\geq \sumt\Pr[Z_t] E\left[\sumd \mathbb I_{X_{t+\Delta}} | Z_t \right] \notag\\
 &= \sumt\Pr[Z_t] \sumd \Pr \big( Y_{t+\Delta}=0 | Y_t=0 \big)  \label{yt-markov}\\
 &\geq \sumt\Pr[Z_t] \min_t \sumd \Pr \big( Y_{t+\Delta}=0 | Y_t=0 \big) \notag\\
 &\geq \Pr[\exists t : Y_t=0] \min_t \sumd \Pr \big( Y_{t+\Delta}=0 | Y_t=0 \big) \label{ES>Pr*min},
\end{align}
where in \eqref{yt-markov} we used the fact that the $Y_t$'s are Markovian, and in the last line the range of $t$ is given by $n^{2/3}/\lambda \leq t \leq \lambda n^{2/3}$.

Next we estimate $\Pr[Y_{t+\Delta}=0|Y_t=0]$ and $\Pr[Y_t=0|Y_0=1]$.  Since
$$Y_{s+\Delta} = Y_s - \Delta + \Bin(n-s-Y_s,1-(1-\p)^{\Delta}),$$
we seek
$$\Pr\Big(\Bin(m,r)=E[\Bin(m,r)]+x\Big)$$
where $m=n-s-Y_s$, $r=1-(1-\p)^\Delta$, and $x=\Delta-Y_s-mr$.

We have $m=(1+o(1))n$ and $r=\Delta \p(1+O(\Delta \p)) = \Theta(\Delta \p)$.  Thus $E[\Bin(m,r)]=mr=\Theta(\Delta)$ and $\Var[\Bin(m,r)] = \Theta(\Delta)$.  We approximate $x$ by
\begin{align*}
x &= \Delta - Y_s - (n-s-Y_s)(\Delta \p + O(\Delta^2 \p^2)) \\
  &= -Y_s + \Delta(1-\p n + \p s + \p Y_s + O(\Delta n\p^2)). \\
\intertext{Now assuming $Y_s=0$ or $1$, and $s=O(\lambda n^{2/3})$, and writing $\p$ as $(1+\widetilde\varepsilon)/n$, we further approximate}
x &= -Y_s + \Delta(-\widetilde\varepsilon+O(\lambda n^{-1/3}) +O(1/n) + O(\lambda n^{-1/3})) \\
  &= -Y_s + \Delta O(\varepsilon) \\
  &= \Delta O(\varepsilon) \ \ \text{if $Y_s=0$}.
\end{align*}

The normal approximation to the binomial is valid when
\[
x\ll \Var[\Bin(m,r)]^{2/3};
\]
see Feller \cite[Volume I, chapter VII.3]{feller-1}.
In particular, if $x\!=\!O(\Var[\Bin(m,r)]^{1/2}) \!=\! O(\Delta^{1/2})$, which happens when $\Delta = O(1/\varepsilon^2)=O(n^{2/3}/\lambda^2)$, we have
\begin{align}
\Pr[\Bin(m,r)=E[\Bin(m,r)]+x] &= \Theta(\Var[\Bin(m,r)]^{-1/2}) = \Theta(\Delta^{-1/2}) \notag\\
\intertext{which implies}
\Pr \big( Y_{s+\Delta}=0 | Y_s=0 \big) &= \Theta(\Delta^{-1/2}) \notag\\
\intertext{and hence}
\sumd \Pr \big( Y_{t+\Delta}=0 | Y_t=0 \big) &= \Theta(\Delta^{1/2})\Big|^{n^{2/3}/\lambda^2}_0 = \Theta(n^{1/3}/\lambda). \label{stutter}
\end{align}

Next we estimate $E[S]$.
If $s=0$ and $t=\Delta$, we bound $x$ by
\begin{align*}
x &= t - 1 - (n-1)(1-(1-\p)^t)
   \leq t-n + n (1-\p)^t
   \leq t - n t \p + n \binom{t}{2} \p^2 \\
  &\leq t - n t (2p-p^2) + 2 n t (t-1) p^2 \leq t - n t 2p + 2 n t^2 p^2
   \leq -t\varepsilon + t^2 (1+\varepsilon)^2/(2n) \\
  &\leq -t\varepsilon/9,
\end{align*}
provided $t\leq (\lambda+1/\lambda^2) n^{2/3} \leq n$ and $\varepsilon
< 1/3$.  The deviations $x$ from the mean are too large for the normal
approximation to be valid.  In this case one would typically use the
Chernoff bound, but we also need the $1/\sqrt{\Var}$ term not found in
the standard Chernoff bound.  Thus we note that the following variation
of the bound is easily deduced from the derivation given by Feller
of the normal approximation to the binomial distribution.
\begin{align}
\Pr[Y_t=0] &= O\big(\frac{1}{\sqrt{t}} e^{-\Theta(x^2/t)}\big) \notag\\
 &= O\big(\frac{1}{\sqrt{t}} e^{-\Theta(t\varepsilon^2)}\big)
 \ \ \text{recalling $|x|\geq t\varepsilon/9$}\notag\\
 &=  O\big(\frac{1}{\sqrt{n^{2/3}/\lambda}} e^{-\Theta(\lambda)}\big)
 \ \ \text{since $t\geq n^{2/3}/\lambda$} \notag\\
 &= O\big(n^{-1/3} e^{-\Theta(\lambda)}\big).\notag\\
\end{align}
{Now summing over $t$, we get}
\begin{align}
E[S] = \sum_{n^{2/3}/\lambda \leq t \leq \lambda n^{2/3}
 + n^{2/3}/\lambda^2}\Pr[Y_t=0] &= O\big(n^{1/3} e^{-\Theta(\lambda)}\big). \label{ES}
\end{align}

Combining the bounds \eqref{ES>Pr*min}, \eqref{stutter},
 and \eqref{ES}, we see that
\begin{equation}\Pr[\text{$Y_t=0$ for some $t$ such that $n^{2/3}/\lambda \leq t \leq \lambda n^{2/3}$}] = O\big(\exp(-\Theta(\lambda))\big). \label{pr-=0}
\end{equation}
We also need
\begin{equation}
\Pr[Y_{n^{2/3}/\lambda}\leq0] \leq \exp(-\Theta(\lambda)), \label{pr-<=0}
\end{equation}
 which follows from the straight Chernoff bound.
Substituting \eqref{pr-=0} and \eqref{pr-<=0} into \eqref{n**-1/3*} we obtain
the desired inequality
$$\Pr[n^{2/3}/\lambda \leq |C(0)| \leq \lambda n^{2/3}]
= O\big(n^{-1/3} \exp(-\Theta(\lambda))\big).$$

\bigskip

\noindent
{\bf Acknowledgements:}  The authors thank Dimitris Achlioptas  and
Riccardo Zecchina for
many useful discussions; in particular, it was Riccardo Zecchina who
suggested to us the significance of the backbone density.
We also thank the Institute for Advanced
Study in Princeton, where this work was begun.  Finally, one of us
(B. Bollob\'as) thanks Microsoft Research, where most of this work
was done.

\bigskip

\bibliographystyle{amsplain}

\begin{thebibliography}{{[MZKST99]}}

\bibitem[Ach00]{Ach00}
D. Achlioptas.
Setting 2 variables at a time yields a new lower bound for random 3-{SAT}
(extended abstract),
\textit{Proc.\ 32nd ACM Symposium on Theory of Computing}, 28--37
(2000).

\bibitem[Ald90]{Aldous90}
D.J. Aldous.
 A random walk construction of uniform spanning trees and uniform
  labelled trees,
 {\em SIAM J.\ on Discrete Mathematics}, \textbf{3}(4):450--465 (1990).

\bibitem[AM98]{AM98}
D. Achlioptas and M. Molloy.
Personal communication (1998).

\bibitem[APT79] {APT79} B. Aspvall, M.F. Plass and R.E. Tarjan.
A linear-time algorithm for testing the truth of certain
quantified Boolean formulas,
\textit{Inf.\ Process.\ Lett.\ } \textbf{8}:121--123 (1979).

\bibitem[AS92]{alon-spencer}
N. Alon and J. Spencer.
\textit{The Probabilistic Method},
John Wiley \& Sons, xiii + 254 pp (1992).

\bibitem[AS00]{AS00}
D. Achlioptas and G.B. Sorkin.
Optimal myopic algorithms for random 3-SAT,
\textit{Proc.\ 41st Symposium on the Foundations of Computer Science},
590--600 (2000).

\bibitem[BBCK98] {BBCK98} B. Bollob\'as, C. Borgs, J.T. Chayes
and J.H. Kim.  Lecture at the Workshop on the Interface between
Statistical Physics and Computer Science, Torino, Italy,
unpublished (1998).

\bibitem[BBCKW00] {BBCKW00} B. Bollob\'as, C. Borgs, J.T. Chayes,
J.H. Kim and D.B. Wilson.  Critical exponents of the 2-SAT transition,
in preparation (2000).

\bibitem[BCKS98a] {BCKS98a} C. Borgs, J.T. Chayes, H. Kesten and
J. Spencer. Uniform boundedness of critical crossing probabilities
implies hyperscaling, to appear in \textit {Rand.\ Struc.\ Alg.\ }

\bibitem[BCKS98b] {BCKS98b} C. Borgs, J.T. Chayes, H. Kesten and
J. Spencer. Birth of the infinite cluster:  Finite-size scaling
in percolation, preprint (1998).

\bibitem[BCM90] {BCM90} E.A. Bender, E.R. Canfield and B.D. McKay.
The asymptotic number of labeled connected graphs with a given number
of vertices and edges,
\textit{Rand.\ Struc.\ Alg.\ }
\textbf{1}:127--169 (1990).

\bibitem[BFU93] {BFU93} A. Broder, A. Frieze and E. Upfal. On the
satisfiability and maximum satisfiablity of random 3-CNF
formulas, \textit{Proc.\ 4th ACM-SIAM Symposium on Discrete
Algorithms}, 322--330 (1993).

\bibitem[Bol84] {Bol84} B. Bollob\'as.  The evolution of random
graphs,
\textit{Trans.\ Amer.\ Math.\ Soc.\ }  \textbf{286}:257--274 (1984).

\bibitem[Bol85] {Bol85} B. Bollob\'as. \textit{Random Graphs},
Academic Press,
London, xvi + 447 pp (1985).

\bibitem[Bri89] {Bri89} V.E. Britikov.
O strukture slucha\u\i nogo grafa vblizi kritichesko\u\i\ tochki,
\textit{Diskretn\t{a\i}a Matematika}, \textbf{1}:121--128 (1989).
English translation, On the random graph structure near the critical point,
\textit{Discrete Math.\ Appl.}, \textbf{1}:301--309 (1991).

\bibitem[CA93] {CA93} J.M. Crawford and L.D. Auton.
Experimental results on the crossover point in satisfiability problems,
\textit{Proc.\ 11th Natl.\ Conf.\ on Artificial Intelligence},
21--27 (1993).

\bibitem[CF86] {CF86} M.T. Chao and J. Franco.   Probabilistic
analysis of two heuristics for the 3-satisfiability problem,
\textit{SIAM J.\ on Computing} \textbf{5}:1106--1118
(1986).

\bibitem[CF90] {CF90} M.T. Chao and J. Franco. Probabilistic
analysis of a generalization of the unit-clause literal
selection heuristics for the $k$ satisfiable problem,
\textit{Information Science} \textbf{51}:289--314 (1990).

\bibitem[Cha98]{Cha98} J.T. Chayes. Finite-size scaling in
percolation, \textit{Proc.\ of the International Congress of
Mathematicians} Vol.\ III (Berlin, 1998), \textit{Doc.\ Math.\ }
Extra Vol.\ III, 113--122 (1998).

\bibitem[Che52]{Che52} H. Chernoff. A measure of asymptotic efficiency
for tests of a hypothesis based on the sum of observations,
\textit{Ann.\ Math.\ Stat.\ } \textbf{23}:493--507 (1952).

\bibitem[Coo71] {Coo71} S.A. Cook. The complexity of theorem-proving
procedures, \textit{Proc.\ 3rd ACM Symposium on Theory of
Computing}, 151--158 (1971).

\bibitem[CPS99] {CPS99} J.T. Chayes, A. Puha and T. Sweet. Independent
and dependent percolation, Probability theory and applications (Princeton,
NJ, 1996), 49--166 in \textit{IAS/Park City Math.\ Ser.\ Vol. 6},
 Amer.\ Math.\ Soc.,
Providence, RI (1999).

\bibitem[CR92] {CR92} V. Chv\'atal and B. Reed. Mick gets some
(the odds are on his side), \textit{Proc.\ 33rd Symposium
on the Foundations of Computer Science},
620--627 (1992).

\bibitem[CS88] {CS88} V. Chv\'atal and E. Szemer\'edi. Many hard
examples for resolution, \textit{J.\ ACM}
\textbf{35}:759--768 (1988).

\bibitem[DB97]{DB97} O. Dubois and Y. Boufkhad. A general
upper bound for the
satisfiablity threshold of random $k$-SAT formulas, \textit{J.\ Algorithms}
\textbf{24}:395--420 (1997).

\bibitem[DBM99]{DBM99}
O. Dubois, Y. Boufkhad, and J. Mandler.
Typical random 3-{SAT} formulae and the satisfiability threshold.
Research announcement at ICTP, Sept.\ 1999.  Two-page abstract
  appears in \textit{Proc.\ 11th ACM-SIAM Symposium on
  Discrete Algorithms}, 126--127 (2000).

\bibitem[Dub91] {Dub91} O. Dubois. Counting the number of solutions
for instances of satisfiability, \textit{Theoretical
Computer Science} \textbf{81}:49--64 (1991).

\bibitem[EF95]{EF95} A. El Maftouhi and W. Fernandez de la Vega.
On random $3$-sat.
\textit{Combin.\ Probab.\ Comput.\ } \textbf{4}:189--195 (1995).

\bibitem[ER60] {ER60} P. Erd\H os and A. R\'enyi.
On the evolution of random graphs,
\textit{Magyar Tud.\ Akad.\ Mat.\ Kutat\'o Int.\ K\"ozl.\ }
\textbf{5}:17--61 (1960).

\bibitem[ER61] {ER61} P. Erd\H os and A. R\'enyi.
On the evolution of random graphs,
\textit{Bull.\ Inst.\ Internat.\ Statist.\ }
\textbf{38}:343--347 (1961).

\bibitem[FB99] {FB99} E. Friedgut, with appendix by J. Bourgain.
Sharp thresholds of graph properties, and the $k$-sat problem,
\textit{J.\ Amer.\ Math.\ Soc.\ } \textbf{12}:1017--1054 (1999).

\bibitem[Fel68] {feller-1} W. Feller. \textit{An Introduction to
Probability Theory and Its Applications},
Volume I, 3rd Edition,
John Wiley \& Sons,
London, xviii + 509 pp (1968).

\bibitem[Fer92] {Fer92} W. Fernandez de la Vega. On random 2-SAT,
unpublished manuscript
(1992).

\bibitem[Fer98] {Fer98} W. Fernandez de la Vega. On random 2-SAT
(revised version),
preprint (1998).

\bibitem[FKG71]{FKG71} C.M. Fortuin, P.W. Kasteleyn and J. Ginibre.
Correlation inequalities on some partially ordered sets,
\textit{Commun.\ Math.\ Phys.\ } \textbf{22}:89--103 (1971).

\bibitem[FP83] {FP83} J. Franco and M. Paul. Probabilistic
analysis of the Davis-Putnam procedure for solving the
satisfiability problem, \textit{Discrete Applied Mathematics}
\textbf{5}:77--87 (1983).

\bibitem[FS96] {FS96} A. Frieze and S. Suen. Analysis of two simple
heuristics for a random instance of $K$-SAT,
\textit{J.\ Algorithms} \textbf{20}:312--335 (1996).

\bibitem[GJ79] {GJ79} M.R. Garey and D.S. Johnson.
\textit{Computers and Intractability:  A Guide to the
Theory of NP-Completeness}, New York, (1979).

\bibitem[GJS76] {GJS76} M.R. Garey, D.S. Johnson
and L. Stockmeyer. Some simplified NP-complete graph
problems, \textit{Theor.\ Comp.\ Sci.\ } \textbf{1}:237--267
(1976).

\bibitem[Goe92] {Goe92} A. Goerdt. A threshold for unsatisfiability,
\textit{Mathematical Foundations of Computer Science,
 17th Intl.\ Symposium},
I.M. Havel and V. Koubek, Eds.,
Lecture Notes in Computer Science \#629, Springer Verlag, 264--274 (1992).

\bibitem[Goe96] {Goe96} A. Goerdt. A threshold for unsatisfiability,
\textit{J.\ Computer and
System Sciences} \textbf{53}:469--486 (1996).

\bibitem[Goe99] {Goe99} A. Goerdt. A remark on random 2-SAT,
\textit{Discrete Applied Mathematics} \textbf{96--97}:107--110 (1999).

\bibitem[Har60]{Har60} T.E. Harris. A lower bound for the
critical probability in certain percolation processes,
\textit{Proc.\ Camb.\ Phil.\ Soc.\ } \textbf{56}:13--20 (1960).

\bibitem[H{\aa}s97]{Has97} J. H\aa stad.
Some optimal in-approximability results,
\textit{Proc.\ 29th ACM Symposium on Theory of Computation}, 1--10 (1997).

\bibitem[JSV00] {JSV00}
S. Janson, Y.C. Stamatiou, and M. Vamvakari.
 Bounding the unsatisfiability threshold of random 3-{SAT},
\textit{Rand.\ Struc.\ Alg.\ }
\textbf{17}:103--116 (2000).

\bibitem[JK{\L}P94] {JKLP94} S. Janson,  D. Knuth, T. \L uczak,
and B. Pittel.
The birth of the giant component,
\textit{Rand.\ Struc.\ Alg.\ }
\textbf{4}:231--358 (1994).

\bibitem[Kar90] {Kar90} R.M. Karp.
The transitive closure of a random digraph,
\textit{Rand.\ Struc.\ Alg.\ }
\textbf{1}:73--93 (1990).

\bibitem[KKK96] {KKK96} L. Kirousis, E. Kranakis and D. Krizanc.
Approximating the unsatisfiability threshold of random
formulas, \textit{Proc.\ 4th European Symposium
on Algorithms}, 27--38 (1996).

\bibitem[Kle66]{Kle66} D.J. Kleitman.
Families of non-disjoint subsets,
\textit{J.\ Combinatorial Theory} \textbf{1}:153--155 (1966).

\bibitem[KMPS95] {KMPS95} A. Kamath, R. Motwani, K. Palem
and P. Spirakis. Tail bounds for occupancy and the
satisfiability threshold conjecture, \textit{Rand.\ Struc.\ Alg.\ }
\textbf{7}:59--89 (1995).

\bibitem[KS94] {KS94} S. Kirkpatrick and B. Selman. Critical
behavior in the satisfiability of random Boolean expressions,
\textit{Science} \textbf{264}:1297--1301 (1994).

\bibitem[{\L}PW94] {LPW94} T. \L uczak, B. Pittel and J.C. Wierman.
The structure of a random graph at the
point of the phase trasnsition,
\textit{Trans.\ Amer.\ Math.\ Soc.\ }
\textbf{341}:721--748 (1994).

\bibitem[LT93] {LT93} T. Larrabee and Y. Tsuji.
Evidence for satisfiability threshold for random 3CNF formulas,
\textit{Proc.\ AAAI Symposium on
Artificial Intelligence and NP-Hard Problems},
112 (1993).

\bibitem[{\L}uc90]{Luc90} T. \L uczak.
Component behavior near the critical
point of the random graph process,
\textit{Rand.\ Struc.\ Alg.\ }
\textbf{1}:287--310 (1990).

\bibitem[McD89]{McDiarmid89}
C. McDiarmid.
On the method of bounded differences.
In {\em Surveys in Combinatorics, 1989}, pages 148--188 (1989).

\bibitem[MPV87] {MPV87} M. M\'ezard, G. Parisi and M.A. Virasoro.
\textit{Spin Glass Theory and Beyond}, World Scientific, Singapore
(1987).

\bibitem[MSL92] {MSL92} D. Mitchell, B. Selman and H. Levesque.
Hard and easy distributions of SAT problems,
\textit{Proc.\ 10th Natl.\ Conf.\ on Artificial Intelligence},
459--465 (1992).

\bibitem[MZ96] {MZ96} R. Monasson and R. Zecchina. The
entropy of the
$K$-satisfiability problem, \textit{Phys.\ Rev.\ Lett.\ }
\textbf{76}:3881 (1996).

\bibitem[MZ97] {MZ97} R. Monasson and R. Zecchina. Statistical
mechanics of the random $K$-SAT model, \textit{Phys.\ Rev.\ E}
\textbf{56}:1357--1370 (1997).

\bibitem[MZKST99] {MZKST99} R. Monasson, R. Zecchina,
S. Kirkpatrick, B. Selman and L. Troyansky.
$2+p$-{SAT}: Relation of typical-case complexity to the nature of the
  phase transition,
\textit{Rand.\ Struc.\ Alg.\ }
\textbf{15}:414--435 (1999).

\bibitem[SK96] {SK96} B. Selman and S. Kirkpatrick.
Critical behavior in the computational cost of satisfiability testing,
\textit{Artificial Intelligence} \textbf{81}:273--295 (1996).

\bibitem[Ver99] {Ver99} Y. Verhoeven. Random 2-SAT and unsatisfiability,
\textit{Inf.\ Process.\ Lett.\ } \textbf{72}:119--123 (1999).

\bibitem[Wil98] {Wil98} D.B. Wilson. \texttt{http://dbwilson.com/2sat-data/} (1998).

\bibitem[Wil00] {Wil00} D.B. Wilson.  The empirical values of the critical $k$-SAT exponents are wrong.  Preprint \texttt{math.PR/0005136} (2000).

\bibitem[Wri77] {Wri77} E.M. Wright.
The number of connected sparsely edged graphs,
\textit{J.\ Graph Theory}
\textbf{1}:317--330 (1977).

\bibitem[Wri80] {Wri80} E.M. Wright.
The number of connected sparsely edged
graphs III.  Asymptotic results,
\textit{J.\ Graph Theory}
\textbf{4}:393--407 (1980).

\end{thebibliography}

\newpage
\end{document}